\documentclass[a4paper,11pt]{preprint}

\usepackage{a4wide, amsmath, amsfonts, amssymb, array, mathrsfs, amsthm, bbm, bm, appendix, dsfont}
\usepackage[T1]{fontenc}
\usepackage[utf8]{inputenc} % only needed for pdflatex
\usepackage[osf]{newtxtext}
\usepackage{mathrsfs}
\usepackage{dsfont}
\usepackage{microtype}
\usepackage{newtxmath}
\usepackage{microtype}
\usepackage{tikz}
\usepackage{booktabs}
\usepackage[margin=1in]{geometry}
\definecolor{colorLink}{RGB}{0,100,162}
\definecolor{colorCite}{RGB}{8,124,100}

\usepackage[
    colorlinks=true,
    linkcolor=colorLink,
    citecolor=colorCite,
    urlcolor=colorLink,
    linktocpage=true
]{hyperref}

\usepackage{comment}

\usepackage{graphicx}
\graphicspath{{figures/}}

\makeatletter
\let\preprint@thebibliography\thebibliography
\let\preprint@endthebibliography\endthebibliography

\renewenvironment{thebibliography}[1]
  {\preprint@thebibliography{99}}
  {\preprint@endthebibliography}
\makeatother

\usepackage[colorinlistoftodos,prependcaption,color=yellow,textsize=tiny,textwidth=2cm]{todonotes}

\usepackage{xargs}
\usepackage{multirow}
\usepackage[pdftex,dvipsnames]{xcolor}  % Coloured text etc.
\usepackage[shortlabels]{enumitem}

\numberwithin{equation}{section}

\newtheorem{theorem}{Theorem}[section]
\newtheorem{lemma}[theorem]{Lemma}
\newtheorem{corollary}[theorem]{Corollary}
\newtheorem{remark}[theorem]{Remark}
\newtheorem{proposition}[theorem]{Proposition}

\newtheorem{example}[theorem]{Example}
\newtheorem{assumption}[theorem]{Assumption}

\allowdisplaybreaks

\renewcommand{\epsilon}{\varepsilon}
\newcommand{\R}{\mathbb{R}}
\newcommand{\C}{\mathbb{C}}
\newcommand{\T}{\mathbb{T}}
\newcommand{\Z}{\mathbb{Z}}
\newcommand{\N}{\mathbb{N}}

\newcommand{\indicator}[1]{\mathbbm{1}_{#1}}
\newcommand{\Id}{\,\mathrm{d}}
\newcommand{\norm}[1]{\left\lVert#1\right\rVert}

\newcommand{\cP}{\mathcal{P}}

\definecolor{phasegreen}{RGB}{170,205,170}
\definecolor{phasegray}{RGB}{225,225,225}
\definecolor{phasewhite}{RGB}{248,248,248}

\renewcommand{\hat}{\widehat}

\author{Benjamin Gess$^{1}$, Paul Nikolaev$^{2}$}
\institute{
	Technische Universit\"{a}t Berlin \& Max-Planck-Institut f\"{u}r Mathematik in den Naturwissenschaften,
	Email: \href{mailto:benjamin.gess@tu-berlin.de}{\color{black}\texttt{benjamin.gess@tu-berlin.de}} 
	\and 
	Technische Universit\"{a}t Berlin, Email: \href{mailto:paul.nikolaev@tu-berlin.de}{\color{black} \texttt{paul.nikolaev@tu-berlin.de}}
	}

\begin{document} 

\title{Long-time behaviour and bifurcation analysis of a  non-reciprocal two-species
aggregation-diffusion system on the torus
\footnote{\textbf{Acknowledgements:} The first author acknowledges support by the Max Planck Society through the Research Group “Stochastic Analysis in the Sciences (SAiS)" and by the ERC grant FluCo (Grant Agreement No. 101088488). Funded by the European Union. Views and opinions expressed are, however, those of the author(s) only and do not necessarily reflect those of the European Union or the European Research Council Executive Agency. Neither the European Union nor the granting authority can be held responsible for them. The second author acknowledges support by the DFG CRC/TRR 388 “Rough Analysis, Stochastic Dynamics and Related Fields”, Project A11.}}

%\date{\today}

\maketitle

\begin{abstract}
A two-species aggregation–diffusion system on the one-dimensional torus with antisymmetric cross-interactions modeling pursuit–evasion dynamics is studied. The self-interaction strength is used as a bifurcation parameter, and the linear stability of the homogeneous state is characterized by Fourier analysis. It is shown that non-reciprocal coupling enlarges the stability region. In a symmetric diffusion setting, exponential convergence to the homogeneous state is established for sufficiently weak self-interactions by means of an entropy argument, with a decay condition independent of the cross-interaction strength. Under suitable nondegeneracy and spectral assumptions, the first loss of stability is shown to occur through a Hopf bifurcation associated with a single spatial Fourier mode. A branch of time-periodic patterns is thereby obtained, along which the two species oscillate with a fixed phase difference. An explicit formula for the cubic Hopf coefficient is derived, and an essential dependence of the direction and stability of the bifurcating branch on the second spatial harmonic is identified. In the case of matching diffusion coefficients and interaction kernels, sufficiently strong non-reciprocity is shown to transform a subcritical stationary instability of the corresponding reciprocal system into a supercritical Hopf bifurcation with stable periodic solutions. An exchange-of-stability result is proved, and the theory is illustrated for Kuramoto, two-mode cosine, and top-hat interaction kernels.


\medskip

\end{abstract}

\tableofcontents

\section{Introduction}

We study the two-species aggregation--diffusion system
\begin{align} \label{eq: 2species_problem}
\begin{cases}
    \partial_t u_{1}= \sigma_1 \Delta u_{1} + \nabla \cdot \big(  u_{1}  (\alpha \nabla (W_1 * u_1 ) + \gamma \nabla( W* u_2))   \big) \\
    \partial_t u_{2}= \sigma_2 \Delta u_{2} + \nabla \cdot \big(  u_{2}  (\alpha \nabla (W_2 * u_2 ) - \gamma \nabla( W* u_1)) \big)
\end{cases}
\end{align}
posed on the one-dimensional torus \(\T = [-\pi,\pi]\), where \(\sigma_1, \sigma_2 > 0\) are the
diffusivities, \(\alpha \in \R\) is the common self-interaction strength, \(\gamma > 0\) measures the
strength of the cross-interaction, and \(W_1, W_2, W\) are even interaction kernels. 
Throughout,
\(\alpha\) plays the role of the bifurcation parameter from the homogeneous steady state, whereas \(\gamma\) is regarded as a fixed parameter.

The defining feature of \eqref{eq: 2species_problem} is the opposite sign of the cross-interactions: the first species is attracted towards the second, while the second is driven away from the first when \(W\) is attractive. Such interactions are called \emph{non-reciprocal}. Although microscopic forces obey Newton's third law, effective interactions between mesoscopic agents need not be reciprocal when mediated by non-equilibrium environments or arising from the agents' own activity~\cite{Ivlev2015,YouBaskaranMarchetti2020,MartinEtAl2025}. Examples range from chemically active colloids and droplets~\cite{SotoGolestanian2014,Meredith2020} and colloids with a cone of vision~\cite{Lavergne2019} to the chase-and-run migration of cell populations~\cite{Theveneau2013} and predator--prey systems~\cite{Tsyganov2003}.

Non-reciprocity can fundamentally alter the nature of instabilities in aggregation--diffusion systems. In contrast to the stationary bifurcations arising in reciprocal settings, non-reciprocal couplings can induce a Hopf bifurcation, leading to time-periodic patterns in which the two densities oscillate at a common frequency with a fixed phase lag. The onset of oscillations, however, does not by itself determine the dynamics that emerge beyond the instability threshold.

Despite progress based on linear stability analysis, single Fourier mode amplitude equations~\cite{YouBaskaranMarchetti2020}, and numerical simulations, a rigorous analysis of the nonlinear bifurcation in the full aggregation--diffusion system remains challenging. Linear stability analysis alone cannot determine the criticality of the Hopf bifurcation, which is governed by the cubic normal-form coefficient. Moreover, non-reciprocal couplings remove the standard gradient-flow structure that underlies many analytical results for reciprocal aggregation--diffusion systems~\cite{Carrillo2020,Carrillo2026}. This raises two central questions: to what extent can a single-mode truncation capture the full bifurcation picture, and how does the strength of non-reciprocity \(\gamma\) affect not only the onset of oscillations but also the criticality and local stability of the emerging periodic branch?

In this work, we address these questions through a rigorous bifurcation analysis of the non-reciprocal system~\eqref{eq: 2species_problem}. We first establish a Hopf bifurcation and derive an explicit formula for its cubic normal-form coefficient. This allows us to go beyond the linear instability analysis and characterize the nonlinear bifurcation in the full system. In particular, we show that the criticality cannot be determined from the critical Fourier mode \(k^*\) alone: the second harmonic \(2k^*\) also plays a decisive role through its contribution to the normal-form coefficient. This demonstrates that single-mode approximation cannot capture the full effects of the bifurcation. We then investigate the role of non-reciprocity in the matching regime \(\sigma_1=\sigma_2\) and \(W_1=W_2=W\). In this setting, we identify an explicit threshold \(\gamma^*>0\) such that the Hopf bifurcation is subcritical for \(\gamma<\gamma^*\) and supercritical for \(\gamma>\gamma^*\). This behaviour has no analogue in the reciprocal setting \footnote{The corresponding \emph{reciprocal} system is obtained by replacing \(-\gamma\) with \(+\gamma\) in the second equation (see~\eqref{eq: 2species_problem_reciprocal}), while retaining \(\alpha\) as the bifurcation parameter.} studied in~\cite{Carrillo2026}, where the corresponding bifurcation is stationary and the coupling strength does not influence the transition between subcritical and supercritical behaviour. Hence, the results of this work demonstrate that the nonlinear outcome of the instability depends both on interactions between Fourier modes and on the strength of the non-reciprocal coupling.

\subsection{Related literature}

For systems of interacting populations, the structure of solutions to~\eqref{eq: nspecies_problem} has been studied predominantly in the reciprocal setting, and even more extensively in the scalar case. Early works focused on the linear stability of homogeneous states, with the resulting pattern-selection mechanisms depending strongly on the structure of the interaction matrix \((\alpha_{ij})\). We refer to~\cite{pottslewis2019} for biologically motivated nonlocal models with top-hat interaction kernels, to~\cite{giunta2022detecting} for a combination of energy methods, numerical simulations, and asymptotic analysis in the local limit, to~\cite{giunta2024weakly} for a weakly nonlinear analysis near the first instability threshold, and to~\cite{JKME23} for the influence of the spatial dimension. The scalar aggregation--diffusion equation was analysed in the seminal work~\cite{Carrillo2020}. More recently, a complete rigorous description of bifurcations from the homogeneous state, based on simple eigenvalues and Crandall--Rabinowitz theory, was obtained for a two-species system with symmetric coupling in~\cite{Carrillo2026}, extending the scalar analysis of~\cite{Carrillo2020}. In addition, the clustering behaviour of weakly interacting diffusions subject to sufficiently localized attractive interactions has been investigated in~\cite{Gerber2026}.

A common feature of all these works is that the coupling is reciprocal. Consequently, the linearized operator is self-adjoint and the bifurcating solutions are stationary. In contrast, non-reciprocal interactions have attracted considerable attention in the physics literature, where the violation of action--reaction symmetry is known to provide a robust mechanism for the emergence of time-dependent collective behaviour. We refer to~\cite{FruchartVitelli2021} for an overview of non-reciprocal phase transitions and oscillatory phases, to~\cite{YouBaskaranMarchetti2020} for non-reciprocity as a generic mechanism leading to travelling states, to~\cite{Kreienkamp2022} for collective dynamics in chiral active matter, and to~\cite{SahaGolestanian2020,FrohoffThiele2023,BraunsMarchetti2024} for pattern formation in non-reciprocal Cahn--Hilliard models. These studies rely primarily on amplitude equations, normal-form reductions, and numerical continuation techniques.

A particularly influential contribution is~\cite{YouBaskaranMarchetti2020}, where a one-mode approximation is introduced by truncating the dynamics to the dominant Fourier mode. As we demonstrate in Section~\ref{sec: Sub-vs. Supercriticality}, such a reduction is insufficient to capture non-linear effects for the non-reciprocal aggregation--diffusion system~\eqref{eq: 2species_problem}. In our setting, resonant interactions with higher harmonics, and in particular with the \(2k^*\)-mode, play a decisive role in determining the criticality and stability of the bifurcating branch.

Closest to the present work, a Hopf transition leading to standing and travelling waves was reported for a non-reciprocal McKean--Vlasov model in~\cite{ArjunEtAl2026} using weakly nonlinear analysis and numerical simulations. On the other hand, a rigorous Crandall--Rabinowitz analysis of phase transitions in McKean--Vlasov systems was carried out in~\cite{HuLiu2026}; however, the interaction considered there is reciprocal and only stationary bifurcations arise.

\subsection{Main contribution}

This work provides, to the best of our knowledge, the first rigorous proof of a Hopf bifurcation from a spatially homogeneous state in a non-reciprocal aggregation–diffusion system. In contrast to the stationary pitchfork bifurcation of the corresponding reciprocal system, we obtain a branch of genuinely time-periodic solutions. We derive an explicit formula for the cubic Hopf coefficient \(\Lambda\) and show that the strength of non-reciprocity can change the bifurcation from subcritical to supercritical. Combined with an exchange-of-stability result at the first instability, this identifies a mechanism by which non-reciprocity generates stable periodic oscillations.

We first establish the global and linear stability theory of the homogeneous state. Fourier analysis gives the maximal interval \((\alpha_-,\alpha_+)\) containing \(\alpha=0\) on which this state is linearly stable (Proposition~\ref{prop: linear_stability}). Increasing \(\gamma\) preserves stability at every previously stable value of \(\alpha\), so that the stability interval can only expand (Remark~\ref{rem:gamma_dependence}).

For sufficiently weak self-interaction, we also prove exponential convergence to the homogeneous state in relative entropy for smooth kernels and normalized diffusivities \(\sigma_1=\sigma_2=1\) (Proposition~\ref{prop: entropy}). The cross-interaction terms cancel exactly in the entropy dissipation, yielding a decay estimate uniform in \(\gamma\).

\begin{figure}[h]
    \centering
    \includegraphics[width=\textwidth]{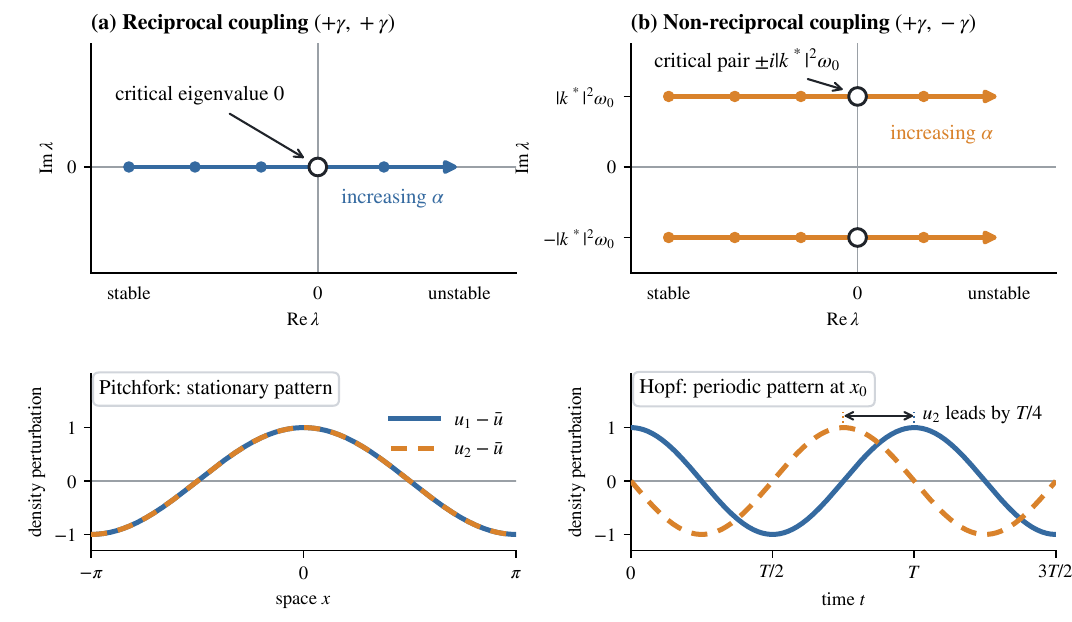}
    \caption{Different mechanisms for the first loss of stability.
    In the reciprocal system, a real eigenvalue crosses zero and a
    stationary branch is produced. In the non-reciprocal system, a
    complex-conjugate pair crosses the imaginary axis at
    \(\pm i|k^*|^2\omega_0\), producing a time-periodic branch. The lower
    panels display the matching regime: the reciprocal critical
    eigenvector has components in phase, whereas the non-reciprocal
    components have a quarter-period phase lag.}
    \label{fig:reciprocal_nonreciprocal_mechanism}
\end{figure}

We next prove a local Hopf bifurcation result at a critical Fourier
mode \(k^*\), which is not assumed to be the first mode to become unstable. Under suitable nondegeneracy and non-resonance assumptions, sufficiently strong non-reciprocity ensures that the eigenvalues associated with this mode cross the imaginary axis transversally at \(\alpha=\alpha_{k^*}\). We construct a branch of periodic solutions and obtain explicit leading-order expressions for their spatial profiles and oscillation frequency (Proposition~\ref{prop: hopf_bifurcation} and Corollary~\ref{cor: hopf_profile}). A key observation is that restricting the analysis to even functions removes the translational degeneracy and makes each eigenvalue of the critical pair simple, allowing the application of the Hopf theorem from simple eigenvalues~\cite[Theorem~I.8.2]{kielh2004}. Figure~\ref{fig:reciprocal_nonreciprocal_mechanism} illustrates the stationary and oscillatory bifurcation mechanisms and the corresponding leading-order patterns.

\begin{figure}[h]
    \centering
    \includegraphics[width=\textwidth]{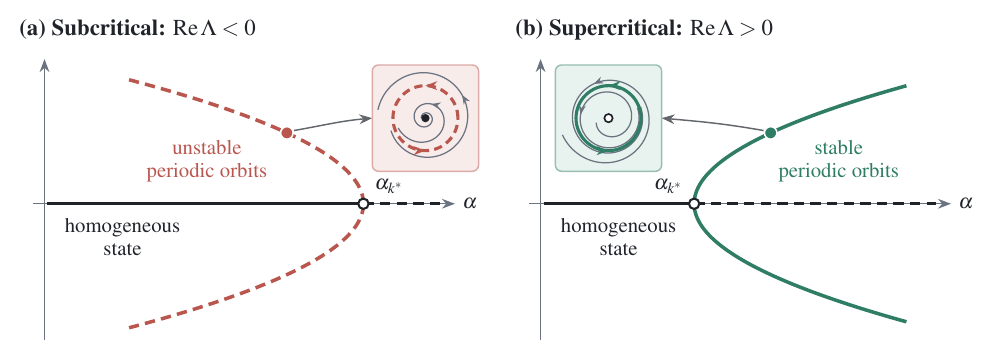}
  \caption{Direction and stability of the bifurcating branch at the first bifurcation point
\(\alpha_{k^*}=\alpha_+\), for positive crossing speed
\(\partial_\alpha\operatorname{Re}\lambda_+(\alpha_{k^*})>0\). The horizontal line is the homogeneous
state \(\bar{\boldsymbol u}\), the curve the branch of periodic solutions; solid lines are stable,
dashed lines unstable. The insets show nearby solutions in the
critical Fourier mode. (a) If \(\operatorname{Re}\Lambda<0\), the branch lies on the side where
\(\bar{\boldsymbol u}\) is stable and consists of unstable periodic orbits. (b) If
\(\operatorname{Re}\Lambda>0\), the branch lies on the side where \(\bar{\boldsymbol u}\) is unstable
and the periodic orbits inherit its stability. At a later bifurcation point the side of the branch is still determined by the sign of \(\operatorname{Re}\Lambda\), but all orbits are unstable.}
    \label{fig: sub_super}
\end{figure}

To determine on which side of \(\alpha_{k^*}\) the periodic solutions exist, we compute the cubic normal-form coefficient \(\Lambda\) explicitly (Proposition~\ref{prop: sub_supercritical}). The sign of its real part distinguishes subcritical from supercritical bifurcation and, together with the direction of the eigenvalue crossing, determines the direction of the branch (see Figure~\ref{fig: sub_super}). 

For a fixed mode \(k^*\), the critical value \(\alpha_{k^*}\) is independent of \(\gamma\), whereas \(\Lambda\) depends on it. The formula involves only the diffusivities, the non-reciprocity strength, and Fourier coefficients of the interaction kernels at \(k^*\) and its second harmonic \(2k^*\). It therefore captures the nonlinear interaction between these modes, which is absent from a truncation to the critical mode alone. We derive explicit criteria when the second-harmonic coefficients vanish or are small, an asymptotic formula for strong non-reciprocity, and a closed-form expression in the matching regime \(\sigma_1=\sigma_2\) and \(W=W_1=W_2\) (Corollaries~\ref{cor: real_part_vanishing_case}--\ref{cor: lambda_matching}).

\begin{figure}[h] 
    \centering
    \includegraphics[width=\textwidth]{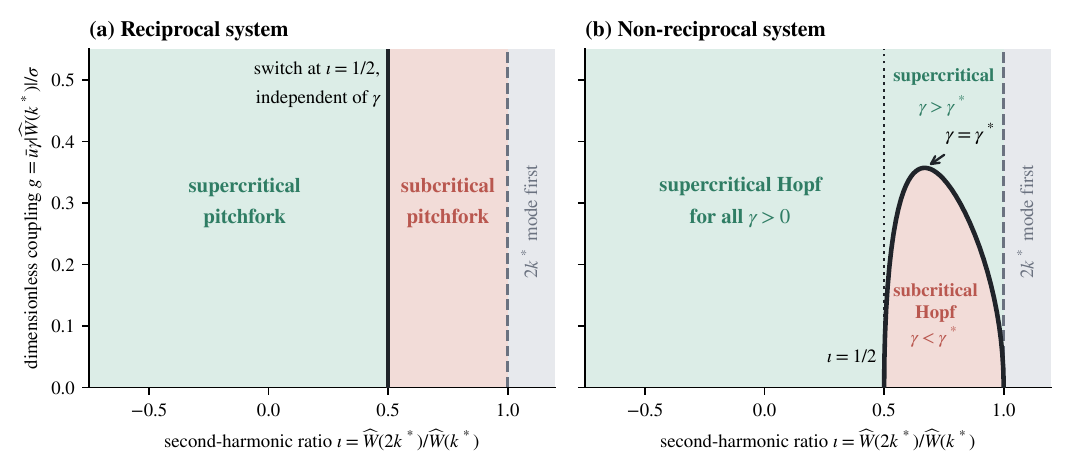}
\caption{Criticality diagrams in the matching regime \(\sigma_1=\sigma_2=\sigma\), \(W_1=W_2=W\),
            \(\widehat W(k^*)<0\). The vertical coordinate is the dimensionless coupling \(g=\bar u\gamma|\widehat W(k^*)|/\sigma\), which equals \(\omega_0/\sigma\) in the non-reciprocal system. (a) In the reciprocal system the pitchfork changes criticality at \(\iota=\tfrac12\), independently of \(\gamma\). (b) In the non-reciprocal system the black curve is the threshold \(\gamma=\gamma^*(\iota)\) of Corollary~\ref{cor: nonreciprocal_transition}, i.e.\ \(g=\sqrt{x^*}/\sigma\). For \(\tfrac12<\iota<1\) the Hopf bifurcation is subcritical below the curve and supercritical above it, and it is supercritical for every \(\gamma>0\) otherwise. The dotted line marks \(\iota=\tfrac12\). In the grey region \(\iota\ge1\) the mode \(2k^*\) destabilizes before \(k^*\).}
    \label{fig:reciprocal_nonreciprocal_regimes} 
\end{figure}

The dependence on \(\gamma\) is particularly explicit in the matching regime. For \(\widehat W(k^*)<0\), set
\[
    \iota:=\frac{\widehat W(2k^*)}{\widehat W(k^*)}.
\]
In the regime \(\iota<1\), the reciprocal bifurcation is supercritical for \(\iota<\tfrac12\) and subcritical for \(\tfrac12<\iota<1\), independently of \(\gamma\) [Appendix~\ref{app: reciprocal_caluclations}]. In the non-reciprocal system, the bifurcation is likewise supercritical for every \(\gamma>0\) when \(\iota<\tfrac12\). For \(\tfrac12<\iota<1\), however, there is a threshold \(\gamma^*>0\) below which the Hopf bifurcation is subcritical and above which it is supercritical [Corollary~\ref{cor: lambda_matching} and Remark~\ref{rem: matching_criticality}]. Thus, with the diffusivities and kernels fixed, increasing the non-reciprocity can change the criticality of the bifurcation. Figure~\ref{fig:reciprocal_nonreciprocal_regimes} displays this dependence: the vertical criticality boundary of the reciprocal problem becomes a threshold curve in the non-reciprocal problem.

Finally, we combine these results to identify the first instability and determine the stability of the emerging periodic orbit. Under explicit conditions on the Fourier coefficients, the local Hopf point satisfies \(\alpha_{k^*}=\alpha_+\), and the remainder of the spectrum lies strictly in the left half-plane (Lemma~\ref{lem: alpha_plus}). The first loss of stability is therefore a Hopf bifurcation. By the Principle of Exchange of Stability, a supercritical branch consists of locally orbitally asymptotically stable periodic solutions (Theorem~\ref{theorem: first_bifurcation}). In the matching regime with \(\tfrac12<\iota<1\), sufficiently strong non-reciprocity therefore produces stable oscillations at the first instability, although the corresponding reciprocal bifurcation is subcritical. Non-reciprocity thus determines both the emergence of time dependence and the stability of the bifurcating oscillations.

The bifurcation results require only bounded variation of the interaction kernels and therefore include top-hat kernels used in ecological modelling. To handle the non-self-adjoint linearized operator and its generally non-normal Fourier blocks, we also establish an invertibility criterion for block-diagonal operators in terms of singular values. This justifies the mode-by-mode inversions used in the computation of \(\Lambda\) (Theorem~\ref{theorem: block_inverse}).

Before stating the results precisely, we summarize the qualitative changes induced by the loss of reciprocity. Throughout, the comparison is with the reciprocal system obtained by replacing \(-\gamma\) with \(+\gamma\) in~\eqref{eq: 2species_problem} (see~\eqref{eq: 2species_problem_reciprocal}) while bifurcating in the same parameter \(\alpha\).
Figure~\ref{fig:reciprocal_nonreciprocal_mechanism} illustrates the corresponding instability mechanisms and the resulting leading-order patterns.

\underline{Reciprocal versus non-reciprocal interactions.}
\begin{itemize}
    \item \textbf{Type of bifurcation and emergent state.} In the reciprocal case the spectrum is real,
    the homogeneous state loses stability through a simple real eigenvalue, and the bifurcation is a
    steady pitchfork producing stationary patterns. In the non-reciprocal case a conjugate pair crosses
    the imaginary axis and the bifurcation is of Hopf type: the emergent states are time-periodic, the
    two densities oscillating at a common frequency.

    \item \textbf{Existence of the bifurcation.} In the reciprocal case a bifurcation point exists as
    soon as some Fourier mode of the kernel is negative, for every \(\gamma \ge 0\). In the
    non-reciprocal case an oscillatory bifurcation requires the non-reciprocity to dominate the mismatch
    between the two species, an explicit condition on \(\gamma\) relating the diffusivities and the
    self-interaction kernels. When the two species have the same intraspecific interaction \(W_1=W_2\) and same diffusion \(\sigma_1=\sigma_2\), then it holds for every \(\gamma > 0\);
    otherwise there is a genuine threshold given by~\eqref{eq: imaginary_condition} below which the crossing is through real eigenvalues.

    \item \textbf{Role of \(\gamma\) in the stability of the homogeneous state.} Here the two systems
    behave in opposite ways. In the reciprocal case increasing \(\gamma\) shrinks the region of linear
    stability in the interaction parameters. In
    the non-reciprocal case the interval \((\alpha_-, \alpha_+)\) widens monotonically in \(\gamma\), and
    converges as \(\gamma \to \infty\) to a nonempty interval determined by the self-interactions alone.
    Non-reciprocity thus stabilizes the homogeneous state at the linear level.

    \item \textbf{What determines criticality.} In the reciprocal case the direction of the branch is
    independent of \(\gamma\) and is decided by the ratio \(\iota\). In
    the non-reciprocal case criticality is governed by the coefficient \(\Lambda\), which depends on
    \(\gamma\). The ratio \(\iota\) no longer decides the criticality alone. 

    \item \textbf{The effect of strong non-reciprocity.} In the matching regime \(\sigma_1=\sigma_2 , W=W_1=W_2\), consider the interval
    \(\iota \in (\tfrac12, 1)\), which is precisely the parameter range where the reciprocal system undergoes a subcritical bifurcation. In contrast, the corresponding non-reciprocal Hopf bifurcation is subcritical only when the non-reciprocity strength \(\gamma\) is below the threshold \(\gamma^*\) and supercritical above it. By the Principle of Exchange of Stability, increasing the degree of non-reciprocity therefore stabilizes the bifurcating branch: beyond the threshold \(\gamma^*\), where the reciprocal system bifurcates subcritically to stationary patterns, the non-reciprocal system bifurcates supercritically to stable oscillations. 
\end{itemize}

\subsection{Organization of the paper}

Section~\ref{sec: setting} fixes the notation, the function spaces, and the standing
assumptions on the interaction kernels. Section~\ref{sec: stability} is devoted to the
homogeneous state: we prove exponential decay in relative entropy
(Proposition~\ref{prop: entropy}), determine the maximal interval of linear stability
(Proposition~\ref{prop: linear_stability}), and analyse the spectrum of the linearized
operator, which is neither self-adjoint nor generally normal. In Section~\ref{sec: local_bifurcation}
we verify the hypotheses of the Hopf bifurcation theorem from a simple eigenvalue and
construct the branch of time-periodic solutions (Proposition~\ref{prop: hopf_bifurcation}).

The criticality of this branch is discussed in Section~\ref{sec: criticallity_hopf}.
Section~\ref{sec: Sub-vs. Supercriticality} derives the general formula for the normal-form
coefficient \(\Lambda\) (Proposition~\ref{prop: sub_supercritical}).
Section~\ref{sec: vanishing_case} evaluates it when the kernels have vanishing second harmonic at
\(2k^*\), and perturbatively when that harmonic is small. In Section~\ref{sec: matching_case} we derive a closed-form expression in the matching case \(\sigma_1=\sigma_2\), \(W=W_1=W_2\), from which
a threshold \(\gamma^*\) is obtained at which an exchange subcritical to supercritical happens.
Section~\ref{sec: high_gamma} treats the regime of strong non-reciprocity and
Section~\ref{sec: exchange} identifies the first
bifurcation point and establishes the exchange of stability, yielding the main theorem
(Theorem~\ref{theorem: first_bifurcation}). Section~\ref{sec: application} applies the
results to the non-reciprocal Kuramoto model, a two-mode cosine kernel, and to the top-hat kernel.

The appendix collects auxiliary material: the algebraic identities used repeatedly in the
computation of \(\Lambda\), an abstract invertibility criterion for block-diagonal operators
in terms of singular values (Theorem~\ref{theorem: block_inverse}) together with its
application to the linearized two-species operator, and the corresponding computation for the
reciprocal system used for comparison (Appendix~\ref{app: reciprocal_caluclations}).

\section{Setting} \label{sec: setting}

Let \(\T^{d} := (\R/2\pi\Z)^{d}\) denote the flat torus of side \(2\pi\), which we
identify with \([-\pi,\pi]^{d}\). We denote by \(\cP(\T^{d})\) the class of Borel
probability measures on \(\T^{d}\); by \(\cP_{ac}(\T^{d})\subset\cP(\T^{d})\) those
absolutely continuous with respect to the Lebesgue measure; and by
\(\cP_{ac}^{+}(\T^{d})\subset\cP_{ac}(\T^{d})\) those having strictly positive
densities almost everywhere in \(\T^{d}\). We write \(C^{m}(\T^{d})\) for the
\(2\pi\)-periodic, \(m\)-times continuously differentiable functions. For an even
kernel \(W\) and \(f\in L^{2}(\T^{d})\), \(W*f\) denotes the spatial convolution
\begin{equation}\label{eq: convolution}
    W*f(x) := \int_{\T^{d}} W(x-y)\,f(y)\,\Id y .
\end{equation}
For \(1\le p\le\infty\), \(L^{p}(\T^{d})\) denotes the usual Lebesgue space;
\(W^{m,p}(\T^{d})\) the periodic Sobolev space of \(m\)-times weakly differentiable
functions with derivatives in \(L^{p}(\T^{d})\); and
\(H^{m}(\T^{d}) := W^{m,2}(\T^{d})\). We denote by \(\mathrm{BV}(\T^{d})\) the space
of functions of bounded variation, i.e.\ those \(f\in L^{1}(\T^{d})\) with a
well-defined distributional derivative \(Df\), where \(Df\) is a finite Radon measure. Throughout the paper we make the following assumption \[ W, W_1, W_2 \in \mathrm{BV}(\T^d). \]

We denote by \(L_{s}^{2}(\T^{d})\subset L^{2}(\T^{d})\) the closed subspace of
even, square-integrable functions. An orthonormal basis of \(L_{s}^{2}(\T^{d})\)
is given by \(\{e_{k}\}_{k\ge 0}\), where
\begin{equation}\label{eq: cosine_basis}
    e_{k}(x) :=
    \begin{cases}
        \sqrt{\tfrac{1}{\pi}}\,\cos(k\cdot x), & k\ge 1,\\[2pt]
        \sqrt{\tfrac{1}{2\pi}}, & k = 0 .
    \end{cases}
\end{equation}
We use the standard scalar product
\(\langle u,v\rangle  = \langle u,v\rangle_{L^{2}} := \int_{\T^{d}} u\,v\,\Id x\). Given \(f\in L_{s}^{2}(\T^{d}) \)
we write \(\hat f(k) := \langle f,e_{k} \rangle \) for its cosine coefficients. However, for the kernels \(W_1, W_2, W\) we make the following convention 
\begin{equation*}
    \hat W_i(k) = \sqrt{\pi} \langle W_i, e_k \rangle, \quad \hat W(k) = \sqrt{\pi} \langle W, e_k \rangle
\end{equation*}
such that the convolution/multiplication formula 
\(\widehat{W*f}(k) =  \hat W(k)\,\hat f(k)\) holds for \(k \ge 1\).
 Let 
\begin{equation*}
    \dot H_s^2(\T) := \{ u \in H^2(\T^d) \cap L^2_s(\T^d): \int_{\T^d} u(x) \Id x = 0\}
\end{equation*}
and 
\begin{equation*}
    \dot L_s^2 (\T) := \{u \in  L^2_s(\T^d): \int_{\T^d} u(x) \Id x = 0 \}. 
\end{equation*}
For the two-species system we consider elements of the product space
\(\boldsymbol f := (f_{1},f_{2})\in L_{s}^{2}(\T^{d})\times L_{s}^{2}(\T^{d})\),
which we equip with the inner product
\begin{equation}\label{eq: product_inner_product}
    \langle \boldsymbol f,\boldsymbol g\rangle_{H}
      := \langle f_{1},g_{1} \rangle  + \langle f_{2},g_{2} \rangle . 
\end{equation}
 From the orthonormal basis \(\{e_{k}\}_{k\ge 0}\) of
\(L_{s}^{2}(\T^{d})\) we obtain a natural orthonormal basis of \( L_{s}^{2}(\T^{d})\times L_{s}^{2}(\T^{d})\),
\begin{equation}\label{eq: product_basis}
    \{(e_{k},0),(0,e_{k})\}_{k\ge 0}.
\end{equation}

For a Banach
space \(\mathcal Y\), an exponent \(\theta\in(0,1)\) and a period \(\tau>0\) we denote by
\(C^{\theta}_{\tau}(\R,\mathcal Y)\) the space of \(\tau\)-periodic,
\(\theta\)-Hölder continuous functions \(u \colon\R\to\mathcal Y\), equipped with the norm
\begin{equation*}
  \norm{u}_{\mathcal Y,\theta}
  := \sup_{t\in\R}\norm{u(t)}_{\mathcal Y}
   + \sup_{\substack{t,s\in\R\\ t\neq s}}
     \frac{\norm{u(t)-u(s)}_{\mathcal Y}}{|t-s|^{\theta}},
\end{equation*}
and by
\(C^{1+\theta}_{\tau}(\R,\mathcal Y)
   :=\{u \in C^{\theta}_{\tau}(\R,\mathcal Y): \dot u \in C^{\theta}_{\tau}(\R,\mathcal Y)\}\),
normed by \(\norm{u}_{\mathcal Y,\theta}+\norm{\dot u}_{\mathcal Y,\theta}\).

We linearise around the homogeneous steady state
\begin{equation*}
u_1^* = u_2^*  =\bar u
\end{equation*}
from which we are going to bifurcate. On the torus $\T$, these constants are given by
\begin{equation*}
\bar u = \frac{1}{2\pi}.
\end{equation*}
Finally, we introduce the relative entropy with respect to the homogeneous steady state via 
\begin{equation*}
\mathcal H( \boldsymbol{u} \mid \bar {\boldsymbol{u}})
= \sum\limits_{i=1}^2
\int_{\mathbb T}
u_i(x)\log\!\left(\frac{u_i(x)}{\bar u}\right)\,dx.
\end{equation*}
for some function \(\boldsymbol{u}=(u_1, u_2)\).

\section{Stability of the homogeneous state } \label{sec: stability}

The goal of this section is to provide a stability analysis of the system~\eqref{eq: 2species_problem}. We provide an asymptotic convergence result in relative entropy for small \(\alpha \ge 0\). Additionally, we perform a linear stability analysis. In the process we characterize the spectrum of the linearized operator as well as properties of the associated semigroup and resolvent.

\subsection{Global asymptotic analysis}\label{sec: global_asymptotic_analysis}

We prove an exponential convergence result to the stationary state in relative entropy. By Carrillo, Salmaniw, Skrzeczkowski~\cite[Theorem~2.4]{Jakub2026} there exists a classical component-wise strictly positive solution \(\boldsymbol{u} \in ( \mathcal{P}_{\mathrm{ac}}(\T) \cap C^2(\T) ) \times (\mathcal{P}_{\mathrm{ac}}(\T) \cap C^2(\T) ) \) with finite entropy under the initial condition assumption \(\boldsymbol{u}_0 \in (\mathcal{P}_{\mathrm{ac}}(\T) \cap H^4(\T) ) \times (\mathcal{P}_{\mathrm{ac}}(\T) \cap H^4(\T) )\), $\boldsymbol{u}_0$ strictly positive as well as \(W, W_i \in L^\infty(\T) \cap \mathrm{BV}(\T)\). Therefore, the relative entropy is well-defined. Before we state the result we introduce \(W_{iu}\), which denotes the unstable components of the kernel \(W_{i}\), i.e. the portion with the negative Fourier coefficients,
\begin{equation*}
    W_{iu} :=\sum\limits_{k \ge 0}  \min( \langle W_i , e_k \rangle, 0) \;  e_k , \qquad i =1, 2. 
\end{equation*} 
We obtain the following convergence result.

\begin{proposition} \label{prop: entropy}
    Let \(\alpha \ge 0\), \(W,W_1,W_2 \in C^2(\T)\) with \(\norm{\Delta W_{iu}}_{L^\infty(\T)} < \infty\) for \(i=1,2\), \(\sigma_1=\sigma_2=1\), \(\boldsymbol{u}_0 \in (\mathcal{P}_{\mathrm{ac}}(\T) \cap H^4(\T) ) \times (\mathcal{P}_{\mathrm{ac}}(\T) \cap H^4(\T) )\) be strictly positive and 
    \begin{equation*}
    \sum\limits_{i=1}^2 \int_{\T} u_{i}(0,x) \log(u_i(0,x))\Id x < \infty. 
    \end{equation*}
Then, the classical solution \(\boldsymbol{u}\) solving~\eqref{eq: 2species_problem} is exponentially stable in relative entropy for small enough \(\alpha\). More precisely, 
\begin{equation*}
    \mathcal H (\boldsymbol u_t \vert  \boldsymbol{\bar u} ) \le \exp\Big(- (2-2 \alpha \max\limits_{i=1,2} \norm{\Delta W_{iu}}_{L^\infty(\T)} ) t  \Big)  \mathcal H (\boldsymbol u_0 \vert  \boldsymbol{\bar u} ) \quad \textrm{ for all } t \ge 0.
\end{equation*}
\end{proposition}
\begin{proof}
    The proof is a direct adaptation of \cite[Proposition~4.1]{Carrillo2026} and \cite[Theorem~1.1(b)]{Carrillo2020}. The only additional point that must be addressed is the non-reciprocal structure of system~\eqref{eq: 2species_problem}. A straightforward computation of the time derivative of the relative entropy yields
    $$
    \frac{\mathrm d}{\mathrm dt}
    \mathcal H (\boldsymbol u_t \mid \boldsymbol{\bar u})
    =
    -\sum_{i=1}^2
    \int_{\T}
    \frac{|\nabla u_i(x)|^2}{u_i(x)}
    \,\mathrm dx
    +
    \sum_{i=1}^2
    \alpha
    \int_{\T}
    u_i(x)\,
    (\Delta W_i * u_i)(x)
    \,\mathrm dx .
    $$
    The cross-interaction terms vanish identically. Indeed, by the evenness of the kernel \(W\), the function \(\Delta W\) is also even, which implies that the two cross contributions cancel after integration by parts. Consequently, the entropy dissipation estimate is independent of the non-reciprocity parameter \(\gamma\).
Applying the logarithmic Sobolev inequality with constant \(1/2\) to the first term and decomposing each kernel \(W_i\) into its stable and unstable components, the remainder of the argument is identical to that of \cite[Theorem~1.1(b)]{Carrillo2020}. This yields the claimed exponential decay estimate.
\end{proof}

Most methods for aggregation--diffusion equations rest on their gradient-flow structure. We discuss why 
the system~\eqref{eq: 2species_problem} admits no such structure. 
\begin{remark}
System~\eqref{eq: 2species_problem} is a particular case of the \(n\)-species aggregation--diffusion
system
\begin{equation} \label{eq: nspecies_problem}
      \partial_t u_{i}= \sigma_i \Delta u_{i} + \nabla \cdot \Big(  u_{i} \sum\limits_{j=1}^n \alpha_{ij} \nabla(W_{ij} * u_j)  \Big) ,
      \qquad i = 1, \ldots, n ,
\end{equation}
obtained by taking \(n=2\), \(\alpha_{ii} W_{ii} = \alpha W_i\) and
\(\alpha_{12} W_{12} = \gamma W = - \alpha_{21} W_{21}\). We restrict ourselves to two species, which is
the smallest setting in which non-reciprocity can occur, and to one spatial dimension, where the Fourier
modes of an even kernel are simple; both restrictions keep the analysis transparent without affecting
the mechanism we describe.

For \(n \ge 2\) the system \eqref{eq: nspecies_problem} need not possess a gradient-flow structure. It
does so if and only if there exist weights \(\pi_1, \ldots, \pi_n > 0\) such that
\begin{equation} \label{eq: detailed_balance}
    \pi_i \alpha_{ij} W_{ij} (x-y) = \pi_j \alpha_{ji} W_{ji} (y-x) , \qquad
    i,j = 1, \ldots, n , \quad x,y \in \T .
\end{equation}
This is known as the \emph{detailed balance condition}, which was introduced by
Boltzmann~\cite{Boltzman}.
Heuristically, \eqref{eq: detailed_balance} says that the matrix of interactions \(\alpha_{ij} W_{ij}\)
is symmetrisable through multiplication by a positive vector \((\pi_1, \ldots, \pi_n)\). When
\(W_{ij} = W\) for all \(i,j\), it reduces to the requirement that the coefficient matrix
\((\alpha_{ij})\) be symmetric. This condition has been used to establish the well-posedness of
solutions~\cite{Jakub2026,Portisch2022}, and it is what makes the bifurcation analysis of the reciprocal
problem manageable. The stationary states are exactly the critical points of a free energy, equivalently
the fixed points of an explicit nonlinear map, and the linearized operator is self-adjoint, so that its
spectrum is real. This is the structure exploited in the scalar case~\cite{Carrillo2020} and, for two
species with symmetric coupling, in~\cite{Carrillo2026}.

The system \eqref{eq: 2species_problem} violates \eqref{eq: detailed_balance} for every
\(\gamma \neq 0\). 
Indeed, taking \(i=1\), \(j=2\) and using that \(W\) is even, condition
\eqref{eq: detailed_balance} would read
\begin{equation*}
   \pi_1 \gamma W(x-y) = - \pi_2 \gamma W(y-x) = - \pi_2 \gamma W(x-y) ,
\end{equation*}
so that \((\pi_1 + \pi_2) \gamma W \equiv 0\). Since \(\gamma \neq 0\) and \(W \not\equiv 0\), this
forces \(\pi_1 = -\pi_2\), which is incompatible with the positivity of the weights. 
Hence, no reweighting of the species restores a gradient-flow structure, however weak the
non-reciprocity \(\gamma\) may be. 
The antisymmetry of the system~\eqref{eq: 2species_problem} implies that the linearized operator is no longer self-adjoint, and its spectrum need not be real. Finally, its
Fourier blocks are non-normal (under the condition \(W_1\neq W_2\) or \(\sigma_1\neq \sigma_2\)), so that the spectral theorem does not apply. 
\end{remark}

\subsection{Spectral analysis of the linearized aggregation-diffusion equation}\label{subsec: spectral_analysis}

Before we prove the linear stability result (Proposition~\ref{prop: linear_stability}) we start by analyzing the spectrum of the linearized operator. 
Let us rewrite the problem into an infinite dimensional evolution equation.
Let \(\tilde G \colon \operatorname{D}(\tilde G) \to  L_s^2(\T) \times  L^2_s(\T)    \) be the parabolic problem~\eqref{eq: 2species_problem} with domain \(\operatorname{D}(\tilde G) = (H^2(\T) \cap L_s^2)\times (H^2(\T) \cap L_s^2) \). More precisely, for \(\boldsymbol{u} = (u_1, u_2)\) we define 
\begin{align*}
    \tilde G (\boldsymbol{u}, \alpha) := 
    \begin{pmatrix}
        \sigma_1 \Delta u_{1} + \nabla \cdot \big(  u_{1}  (\alpha \nabla (W_1 * u_1 ) + \gamma \nabla( W* u_2) )   \big) \\
        \sigma_2 \Delta u_{2} + \nabla \cdot \big(  u_{2}  (\alpha \nabla (W_2 * u_2 ) - \gamma \nabla( W* u_1)) \big) 
    \end{pmatrix}.
\end{align*}

We recenter the map and define 
\begin{equation*}
    G(\boldsymbol{u}, \alpha) =  \tilde G(\boldsymbol{u} + \boldsymbol{\bar u}, \alpha)
\end{equation*}
such that \(G(\boldsymbol{0}, \alpha) = 0 \) and \(\operatorname{D}(G) = \dot H_s^2(\T)\times \dot H_s^2(\T) \). Hence, we can rewrite the two-species Fokker--Planck equation into an abstract evolution equation of the form 
\begin{equation} \label{eq: abstract_problem}
    \partial_t \boldsymbol u  =  G (\boldsymbol{u}, \alpha). 
\end{equation}
Next, we introduce the perturbation ansatz 
\begin{equation*}
u_i = \bar u + \varepsilon \rho_i
\end{equation*}
with \(\rho_i \in \dot H_s^2(\T)\). 

Substituting into the first equation yields
\begin{align*}
\partial_t (\bar u+\varepsilon \rho_1) &=\sigma_1 \Delta (\bar u+\varepsilon \rho_1) + \nabla\cdot \Big[ (\bar u+\varepsilon \rho_1)
\Big( \alpha \nabla(W_1*(\bar u+\varepsilon\rho_1)) + \gamma \nabla(W*(\bar u+\varepsilon\rho_2)) \Big) \Big].
\end{align*}

Since convolutions of constants are constant, we have $\nabla(W_i * \bar u)=0, \nabla(W * \bar u)=0$ and thus
\begin{align*}
&\nabla\cdot \Big[ (\bar u+\varepsilon \rho_1) \Big( \alpha \nabla(W_1*(\bar u+\varepsilon\rho_1)) + \gamma \nabla(W*(\bar u+\varepsilon\rho_2)) \Big)\Big]
\\
&= \nabla\cdot \Big[ (\bar u+\varepsilon \rho_1) \Big( \varepsilon \alpha \nabla(W_1*\rho_1) + \varepsilon \gamma \nabla(W*\rho_2) \Big) \Big].
\end{align*}

Keeping only terms of order $\epsilon$ we find,
\begin{align*}
 \varepsilon \bar u \nabla\cdot \Big( \alpha \nabla(W_1*\rho_1) + \gamma \nabla(W*\rho_2) \Big) + O(\varepsilon^2)
=  \varepsilon \bar u \Big( \alpha \Delta(W_1*\rho_1) + \gamma \Delta(W*\rho_2)\Big) + O(\varepsilon^2).
\end{align*}
By the regularity of \(\rho_1\), we have \(\norm{\nabla \cdot (\rho_1 \nabla W_1*\rho_1)}_{L^2(\T)} \le C \norm{\rho_1}_{H^2(\T)}^2\), where \(C\) depends on the bounded variation norm of \(W_1\). Hence, for any bounded set \(B \subset H^2(\T)\) the term \(O(\epsilon^2)\) is uniform in \(\rho_1\).

Proceeding analogously for the second equation we demonstrated that \(\tilde G\) is G\^{a}teaux differentiable uniformly with respect to any bounded set \(B\).
Hence, we obtain Fr\'echet differentiability of \(\tilde G\), and consequently Fr\'echet differentiability of \(G\) with \(D_u G(\boldsymbol 0, \alpha) \colon \operatorname{D}(G) \to \dot  L^2_s(\T) \times \dot L^2_s(\T)  \) corresponding to the linear problem
\begin{align} \label{eq: 2species_linearized}
\begin{cases}
\partial_t \rho_1 = \sigma_1 \Delta \rho_1 + \bar u \alpha \Delta(W_1*\rho_1) +\bar u \gamma \Delta(W*\rho_2), \\
\partial_t \rho_2 = \sigma_2 \Delta \rho_2 + \bar u \alpha \Delta(W_2*\rho_2) - \bar u \gamma \Delta(W*\rho_1).
\end{cases}
\end{align}

By similar arguments we derive the second Fr\'echet derivative 
\begin{equation} \label{eq: second_derivative}
        D_{uu}^2 G(\boldsymbol 0, \alpha)[\boldsymbol{\rho}, \boldsymbol{\eta}]
        =
        \begin{pmatrix}
            \alpha \nabla \cdot \big(  \nabla W_1 * \eta_1 \rho_1 +  \nabla W_1*\rho_1 \eta_1 \big) 
            + \gamma \nabla \cdot \big(  \nabla W* \eta_2 \rho_1 + \nabla W*\rho_2 \eta_1 \big)  \\
              \alpha  \nabla \cdot \big(  \nabla W_2 * \eta_2 \rho_2 +  \nabla W_2*\rho_2 \eta_2 \big) 
            - \gamma \nabla \cdot  \big(  \nabla W* \eta_1 \rho_2 + \nabla W*\rho_1 \eta_2 \big) 
        \end{pmatrix}.
\end{equation}
In particular, the derivative is independent of \(\boldsymbol{u}\). Therefore, \(G\) is smooth.

Moreover, we can readily check the normal ellipticity condition of~\cite[Theorem~2.4]{amann1990}, which implies in combination with the perturbation result~\cite[Theorem~3.2.10]{engel2000} and~\cite[Lemma~A.1.1]{Carrillo2026} that the linearisation of the system~\eqref{eq: 2species_problem} \(D_u G(\boldsymbol{0},\alpha)\) generates an analytic semigroup for \(\sigma_i > 0, i =1, 2\). Since, by the Rellich--Kondrachov theorem the embedding \(\operatorname{D}(G) \hookrightarrow \dot L^2(\T) \times \dot L^2(\T)\) is compact, it follows that the resolvent is compact, which implies the compactness for \(t > 0\) of the analytic semigroup generated by \(D_u G(\boldsymbol{0},\alpha)\). 
Moreover, the spectrum of the operator \(D_u G(\boldsymbol{0},\alpha)\) consists only of the point spectrum by Lemma~\ref{lemma: block_spectrum}. 
We want to analyse the point spectrum of the operator \(D_u G(\boldsymbol{0},\alpha)\). Thus, we introduce the Fourier expansion
\begin{equation*}
\rho_i(t,x)=\sum_{k \ge 1 }\hat\rho_i(k,t)e_k (x) .
\end{equation*}

Using $\widehat{W*\rho}(k)=\widehat W(k)\hat\rho(k),$ and $\Delta e_k(x) =-|k|^2 e_k(x)$, we obtain for each Fourier mode \(k\),
\begin{align*}
\partial_t \hat\rho_1(k) &=-|k|^2 \Big( \sigma_1+\bar u\alpha \widehat W_1(k) \Big)\hat\rho_1(k) -|k|^2 \bar u\gamma \widehat W(k)\hat\rho_2(k), \\
\partial_t \hat\rho_2(k) &=-|k|^2\Big( \sigma_2+\bar u\alpha \widehat W_2(k) \Big)\hat\rho_2(k) + |k|^2\bar u\gamma \widehat W(k)\hat\rho_1(k).
\end{align*}

Hence
\begin{equation*}
\partial_t
\begin{pmatrix}
\hat\rho_1(k)\\
\hat\rho_2(k)
\end{pmatrix}
=
A(k)
\begin{pmatrix}
\hat\rho_1(k) \\
\hat\rho_2(k)
\end{pmatrix},
\end{equation*}
with
\begin{equation} \label{eq: def_A_k}
A(k) := -|k|^2
\begin{pmatrix}
\sigma_1+\bar u\alpha \widehat W_1(k)
&
\bar u\gamma \widehat W(k)
\\
-\bar u\gamma \widehat W(k)
&
\sigma_2+\bar u\alpha \widehat W_2(k)
\end{pmatrix}.
\end{equation}
Thus, we obtain 
\begin{equation*}
    D_u G(\boldsymbol 0, \alpha) (c_1 e_k, c_2 e_k) = 
(A(k) c )  
    \begin{pmatrix}
        e_k \\
        e_k 
    \end{pmatrix} 
\end{equation*}
for \(c \in \C^2 \) and \(k \neq 0 \). 
Even though \(A(k)\) is in general not normal, we can 
apply Lemma~\ref{lemma: block_spectrum}. Hence, it is enough to compute the spectrum on each \(A(k)\) since we have 
\begin{equation*}
    \sigma(  D_u G(\boldsymbol 0, \alpha)) = \overline{\bigcup\limits_{k \ge 1} \sigma(A(k))}.
\end{equation*}

The eigenvalues are given by
\begin{equation} \label{eq: eigenvalues_abstract}
\lambda_\pm(k)=\frac{\operatorname{tr}(A(k))}{2}
\pm \sqrt{ \frac{\operatorname{tr}(A(k))^2}{4}- \det(A(k))
}.
\end{equation}

The trace is
\begin{equation*}
\operatorname{tr}(A(k)) = -|k|^2
\Big( \sigma_1+\sigma_2 +\bar u\alpha \widehat W_1(k) +\bar u\alpha \widehat W_2(k) \Big),
\end{equation*}
while the determinant equals
\begin{equation*}
\det(A(k))
= |k|^4\Big(\Big(\sigma_1+\bar u\alpha \widehat W_1(k)\Big)\Big( \sigma_2+\bar u\alpha \widehat W_2(k) \Big)+\bar u^2\gamma^2 \widehat W(k)^2 \Big).
\end{equation*}

\subsection{Linear stability analysis}

The homogeneous state is linearly stable precisely when, for every \(k\ge1\), both eigenvalues of \(A(k)\) have negative real part (see the proof of Proposition~\ref{prop: linear_stability}), which for a real \(2\times2\) matrix means \(\operatorname{tr}A(k)<0\) and \(\det A(k)>0\).
In order to track these conditions we introduce the following abbreviations: 
\begin{equation} \label{eq: S_i_of_k}
    S_1(k) := \sigma_1 + \bar u\alpha \widehat W_1(k),
    \qquad
    S_2(k) := \sigma_2 + \bar u\alpha \widehat W_2(k),
\end{equation}
one has \(\operatorname{tr} A(k) = -|k|^2\bigl(S_1(k)+S_2(k)\bigr)\), which is affine in \(\alpha\), and \(\det A(k) = |k|^4\bigl(S_1(k)S_2(k) + \bar u^2 \gamma^2 \widehat W(k)^2\bigr)\), which is quadratic.
The coefficient of \(\alpha\) in \(\operatorname{tr}A(k)\) is
\(-|k|^2\bar u(\widehat W_1(k)+\widehat W_2(k))\), whose sign decides the direction of the
constraint. Accordingly we split the modes into
\begin{equation*}
    \mathcal T^+ := \{k\ge 1 : \widehat W_1(k)+\widehat W_2(k)>0\},
    \qquad
    \mathcal T^- := \{k\ge 1 : \widehat W_1(k)+ \widehat W_2(k)<0\},
\end{equation*}
and for \(k\in\mathcal T^+\cup\mathcal T^-\) define the trace threshold
\begin{equation*}
    T(k)
    :=
    -\frac{\sigma_1+\sigma_2}
    {\bar u \widehat W_1(k)+\bar u\widehat W_2(k)} ,
\end{equation*}
so that \(\operatorname{tr}A(k)<0\) reads \(\alpha>T(k)\) on \(\mathcal T^+\) and \(\alpha<T(k)\) on \(\mathcal T^-\), while modes with \(\widehat W_1(k)+\widehat W_2(k)=0\) impose no trace constraint.

Moreover, \(\det A(k)=|k|^4D_k(\alpha)\) with the determinant polynomial
\begin{align*}
    D_k(\alpha)
    :&=
    S_1(k)S_2(k) + \bar u^2 \gamma^2 \widehat W(k)^2\\
    &= \bar u^2 \widehat W_1(k)\widehat W_2(k)\,\alpha^2
    +\bar u\bigl(\sigma_1\widehat W_2(k)+\sigma_2\widehat W_1(k)\bigr)\alpha
    +\sigma_1\sigma_2+\bar u^2 \gamma^2\widehat W(k)^2 .
\end{align*}
Since \(D_k(0)=\sigma_1\sigma_2+\bar u^2 \gamma^2\widehat W(k)^2  >0\), the origin lies in \(\{D_k>0\}\) and we may set
\begin{equation*}
    D_+(k):=\inf\{\alpha>0:D_k(\alpha)=0\},
    \qquad
    D_-(k):=\sup\{\alpha<0:D_k(\alpha)=0\},
\end{equation*}
with the conventions \(\inf\varnothing=+\infty\) and \(\sup\varnothing=-\infty\), so that
\((D_-(k),D_+(k))\) is the connected component of \(\{D_k>0\}\) containing \(0\).

Intersecting all constraints and retaining the component of \(\alpha=0\) leads to
\begin{equation*}
    \alpha_-
    :=
    \max\left\{
        \sup_{k\in\mathcal T^+} T(k),
        \ \sup_{k\ge 1} D_-(k)
    \right\},
    \qquad
    \alpha_+
    :=
    \min\left\{
        \inf_{k\in\mathcal T^-} T(k),
        \ \inf_{k\ge 1} D_+(k)
    \right\}.
\end{equation*}
\begin{lemma}\label{lem: tail}
Let \(a>0\). There exists \(R\in\N\), depending only on \(a\), \(\sigma_1,\sigma_2\) and the total variations of \(W_1,W_2\), such that
\[
S_1(k)\ge\tfrac12\sigma_1, \qquad S_2(k)\ge\tfrac12\sigma_2
\]
for all \(k\ge R\), all \(\alpha\in[-a,a]\) and all \(\gamma\ge0\). In particular, for \(k\ge R\) neither \(\operatorname{tr}A(k)\) nor \(\det A(k)\) vanishes on \([-a,a]\). Consequently, \(|T(k)|>a\) for \(k\in\mathcal T^+\cup\mathcal T^-\) with \(k\ge R\), and \(D_-(k)<-a\), \(D_+(k)>a\) for all \(k\ge R\) and all \(\gamma\ge0\).
Moreover, every eigenvalue \(\lambda(k)\) of \(A(k)\) with \(k\ge R\) satisfies
\[
\operatorname{Re}\lambda(k) \le - \frac{\sigma_1\sigma_2}{6(\sigma_1+\sigma_2)} |k|^2,
\]
uniformly in \(\alpha\in[-a,a]\) and \(\gamma\ge0\).
\end{lemma}

\begin{proof}
By~\eqref{eq: BV_decay} we have \(|\hat W_i(k)|\le V(W_i)/k\), where \(V( \cdot)\) denotes the total variation. Choosing \(R\ge 2\bar u a\max_i (V(W_i)/\sigma_i)\), we obtain for \(k\ge R\) and \(|\alpha|\le a\)
\[
S_i(k)\ge\sigma_i-\bar u a|\hat W_i(k)|\ge\tfrac12\sigma_i .
\]
Hence \(\operatorname{tr}A(k)=-|k|^2(S_1(k)+S_2(k))<0\) and \[\det A(k)\ge|k|^4S_1(k)S_2(k) \ge\tfrac14 k^4\sigma_1\sigma_2 >0\] on \([-a,a]\). Since \(T(k)\) is the zero of \(\operatorname{tr}A(k)\) and \(D_\pm(k)\) are zeros of \(\det A(k)\), the remaining assertions follow.

The same choice of \(R\) also gives \(S_i(k)\le\sigma_i+\bar u a|\hat W_i(k)|\le\tfrac32\sigma_i\), hence
\(|\operatorname{tr}A(k)|\le\tfrac32 k^2(\sigma_1+\sigma_2)\).  If the eigenvalues of \(A(k)\) are
complex, then \(\operatorname{Re}\lambda=\tfrac12\operatorname{tr}A(k)\le-\tfrac14k^2(\sigma_1+\sigma_2) \).
If they are real, both are negative. Denote them by \(\lambda_{\min}\le\lambda_{\max}<0\). Then
\(\lambda_{\min}=\operatorname{tr}A(k)-\lambda_{\max}>\operatorname{tr}A(k)\), so \(|\lambda_{\min}|<|\operatorname{tr}A(k)|\) and
\[
|\lambda_{\max}|=\frac{\det A(k)}{|\lambda_{\min}|}>\frac{\det A(k)}{|\operatorname{tr}A(k)|}
\ge\frac{\sigma_1\sigma_2}{6(\sigma_1+\sigma_2)}\,k^2 .
\]
\end{proof}
The next proposition shows that \((\alpha_-,\alpha_+)\) is exactly the maximal stability window.

\begin{proposition} \label{prop: linear_stability}

With the notation above, \(\alpha_- < 0 < \alpha_+\), the homogeneous equilibrium \((\bar u,\bar u)\) is linearly
stable for every \(\alpha\in(\alpha_-,\alpha_+)\), and \((\alpha_-,\alpha_+)\) is the maximal open
interval containing \(0\) on which it is linearly stable. In particular \(\alpha=0\) is always
linearly stable.
\end{proposition}

\begin{remark}[Structure of the determinant constraint]
\label{rem: partition_structure}
The determinant condition \(D_k(\alpha)>0\) can only fail when \(S_1(k)\) and \(S_2(k)\) have
opposite signs and their product outweighs the non-reciprocal cushion \(\bar u^2\gamma^2\widehat
W(k)^2\). Whenever \(S_1(k)S_2(k)\ge 0\) it holds automatically. The geometry of
\(D_k\) is governed by the sign of its leading coefficient \(\widehat W_1(k)\widehat W_2(k)\):
\begin{enumerate}
    \item If \(\widehat W_1(k)\widehat W_2(k)>0\), then \(D_k\) is an upward parabola with
    \(D_k(0)>0\). Its real roots, if any, lie on one side of \(0\) and carve out a forbidden band
    away from the origin, so the mode contributes at most a single one-sided bound
    (\(D_-(k)\) or \(D_+(k)\)). If \(D_k\) has no real root the mode imposes no determinant
    constraint at all.
    \item If \(\widehat W_1(k)\widehat W_2(k)<0\), then \(D_k\) is a downward parabola and
    \(D_k(0)>0\) forces two real roots straddling the origin, so the mode contributes a two-sided
    bound \(D_-(k)<0<D_+(k)\).
    \item If \(\widehat W_1(k)\widehat W_2(k)=0\), then \(D_k\) is affine (one bound), or, if
    additionally \(\widehat W_1(k)=\widehat W_2(k)=0\), constant and positive, imposing no constraint.
\end{enumerate}
Only in the first case can the set \(\{D_k>0\}\) be disconnected; the maximal stability interval
\((\alpha_-,\alpha_+)\) always refers to the component containing \(\alpha=0\), and additional stable
values of \(\alpha\) beyond a forbidden band, if present, are not connected to the trivial state.
\end{remark}

\begin{remark}[Widening of the stability interval under non-reciprocity]
\label{rem:gamma_dependence}
The non-reciprocity parameter \(\gamma\ge 0\) enters neither the trace thresholds \(T(k)\) nor the
leading and linear coefficients of \(D_k\). It appears only in the constant term
\(D_k(0)=\sigma_1\sigma_2+\bar u^2\gamma^2\widehat W(k)^2\). Raising \(\gamma\) therefore
lifts each determinant curve \(D_k\) upward, and in every regime of
Remark~\ref{rem: partition_structure} this pushes the nearest roots \(D_\pm(k)\) away from the origin,
or removes them altogether:
\begin{itemize}
    \item for \(\widehat W_1(k)\widehat W_2(k)>0\) the two roots move together and, once
    \(\bigl(\sigma_1\bar u \widehat W_2(k)+\sigma_2\bar u\widehat W_1(k)\bigr)^2<4\bar u^2 \widehat W_1(k)\widehat W_2(k)\,D_k(0)\),
    disappear, leaving \(D_-(k)=-\infty\), \(D_+(k)=+\infty\);
    \item for \(\widehat W_1(k)\widehat W_2(k)\le 0\) the roots spread apart, so \(D_-(k)\) decreases and
    \(D_+(k)\) increases.
\end{itemize}
Consequently, the window \((\alpha_-(\gamma),\alpha_+(\gamma))\) is non-decreasing in \(\gamma\), and the mechanism is independent of the signs of the self-interactions. Its limit as \(\gamma\to\infty\) is
\begin{equation} \label{eq: trace_interval}
    \lim_{\gamma\to\infty}\bigl(\alpha_-(\gamma),\alpha_+(\gamma)\bigr)
    = I_\infty\cap\bigcap_{k:\,\widehat W(k)=0}\bigl(D_-(k),D_+(k)\bigr),
    \qquad
    I_{\infty} := \Bigl( \sup_{k\in\mathcal T^+} T(k),\ \inf_{k\in\mathcal T^-} T(k) \Bigr),
\end{equation}
which contains \(0\) and depends only on \(\sigma_1,\sigma_2,W_1,W_2\). 
\end{remark}

\begin{proof}[Proof of Proposition~\ref{prop: linear_stability}.]
By Lemma~\ref{lemma: block_spectrum},
\(\sigma(D_uG(\mathbf 0,\alpha))=\bigcup_{k\ge1}\sigma(A(k))\), where \(\sigma(\cdot)\) denotes the spectrum of an operator. For each \(k\ge1\) both eigenvalues of the real
\(2\times2\) matrix \(A(k)\) have negative real part if and only if \(\operatorname{tr}A(k)<0\) and
\(\det A(k)>0\). We claim that \((\bar u,\bar u)\) is linearly stable if and only if these two
inequalities hold for all \(k\ge1\). Necessity is immediate, since \(\sigma(A(k))\subset\sigma(D_uG(\mathbf 0,\alpha))\).

For sufficiency, apply Lemma~\ref{lem: tail} with \(a:=\max(|\alpha|,1)\). Every eigenvalue of a mode
\(k\ge R(a)\) satisfies \(\operatorname{Re}\lambda\le-c\) for some \(c > 0\). Hence, 
\[
\sup\operatorname{Re}\sigma(D_uG(\mathbf 0,\alpha))\le\max\Bigl\{-c,\ \max_{1\le k<R(a)}\max\operatorname{Re}\sigma\bigl(A(k)\bigr)\Bigr\}<0 .
\]

\emph{Trace.} We have \(\operatorname{tr}A(k)<0\) if and only if
\begin{equation*}
    g_k(\alpha):=\sigma_1+\sigma_2+\alpha\bigl(\bar u\widehat W_1(k)+\bar u\widehat W_2(k)\bigr)>0 .
\end{equation*}
If \(\bar u\widehat W_1(k)+\bar u\widehat W_2(k)=0\) this reads \(\sigma_1+\sigma_2>0\) and holds for
every \(\alpha\). For \(k\in\mathcal T^+\), dividing by the positive coefficient gives
\(\alpha>T(k)\), whence \(\alpha>\sup_{k\in\mathcal T^+}T(k)\); for \(k\in\mathcal T^-\), dividing by the
negative coefficient reverses the inequality and gives \(\alpha<\inf_{k\in\mathcal T^-}T(k)\). Since
\(\sigma_1+\sigma_2>0\), one has \(T(k)<0\) on \(\mathcal T^+\) and \(T(k)>0\) on \(\mathcal T^-\); in
particular \(\alpha=0\) satisfies every trace condition.

\emph{Determinant.} We have \(\det A(k)>0\) if and only if \(D_k(\alpha)>0\). As \(D_k(0)=\sigma_1\sigma_2+\bar u^2\gamma^2\widehat W(k)^2>0\), the origin lies in the open set \(\{D_k>0\}\), whose connected component containing \(0\) is exactly the interval \((D_-(k),D_+(k))\) by definition of \(D_\pm(k)\). Hence, on the component of the origin, \(\det A(k)>0\) is equivalent to \(D_-(k)<\alpha<D_+(k)\); intersecting over \(k\) gives \(\sup_{k\ge 1}D_-(k)<\alpha<\inf_{k\ge 1}D_+(k)\).

\emph{Stability window.} Combining the trace and determinant conditions with the first step, the set of \(\alpha\) connected to \(\alpha=0\) at which \((\bar u,\bar u)\) is linearly stable is precisely
\((\alpha_-,\alpha_+)\). Applying Lemma~\ref{lem: tail} with \(a=1\), the modes \(k\ge R(1)\) satisfy \(|T(k)|>1\) and \(D_-(k)<-1<1<D_+(k)\). Each of the finitely many modes \(k<R(1)\) satisfies
\(T(k)\neq0\) and \(D_-(k)<0<D_+(k)\), since \(D_k(0)>0\). Hence \(\alpha_-<0<\alpha_+\).

\emph{Maximality.} Suppose \(\alpha_+<\infty\) and apply Lemma~\ref{lem: tail} with \(a:=\alpha_++1\).
Every mode \(k\ge R(a)\) has \(T(k)>a\) (if \(k\in\mathcal T^-\)) and \(D_+(k)>a\). Therefore
\[
\alpha_+=\min\Big\{\min_{k\in\mathcal T^-,\,k<R(a)}T(k),\ \min_{1\le k<R(a)}D_+(k)\Big\},
\] a minimum over finitely many values. So there exists \(k_0<R(a)\) with \(\alpha_+=T(k_0)\) or \(\alpha_+=D_+(k_0)\). In the first case
\(\operatorname{tr}A(k_0)=0\) at \(\alpha_+\), since \(T(k_0)\) is the unique zero of the affine function \(g_{k_0}\).  Since \(\det A(k_0)>0\) on
\([0,\alpha_+)\), continuity gives \(\det A(k_0) \ge 0\) at \(\alpha_+\). In the second case, \(D_{k_0}\) is a polynomial of degree at most two with \(D_{k_0}(0)>0\), so the set \(\{\alpha>0: D_{k_0}(\alpha)=0\}\) is finite and its infimum \(D_+(k_0)\) is attained. Hence, \(\det A(k_0)=0\) at \(\alpha_+\). In either case \(A(k_0)\) has an eigenvalue on the imaginary axis: \(\pm i\sqrt{\det A(k_0)}\) if the trace vanishes, and \(0\) if the determinant vanishes. Thus \((\bar u,\bar u)\) is not linearly stable at \(\alpha_+\). The same argument applies at \(\alpha_-\). Therefore \((\alpha_-,\alpha_+)\) is the
maximal open interval containing \(0\) on which the homogeneous state is linearly stable.

\end{proof}

\section{Local bifurcation theory} \label{sec: local_bifurcation}
The goal of this section is to describe the local bifurcation structure of system~\eqref{eq: 2species_problem}. In particular, we construct the branch of time-periodic solutions that bifurcates from the spatially homogeneous equilibrium and characterize its dependence on the Fourier coefficients of the interaction kernels \(W\), \(W_1\), and \(W_2\), as well as on the non-reciprocity parameter \(\gamma\). To leading order, the bifurcating branch is generated by a single Fourier mode \(e_{k^*}\) oscillating at the frequency \(\Im(\lambda_+(k^*))\) for a specific \(k^*\). The two species oscillate with the same amplitude, while their relative phase shift is determined by the non-reciprocal coupling.

The construction proceeds by verifying the Crandall--Rabinowitz Theorem for Hopf bifurcation~\cite[Theorem~I.8.2]{kielh2004}. Two major assumptions of this theory are the simplicity of eigenvalues \(\lambda_+\) of the map \(\operatorname{D}_u G(\mathbf{0} , \alpha)\), i.e. \(N(\lambda_+ - \operatorname{D}_u G(\mathbf{0} , \alpha)) = \operatorname{span}\{q\} \) and \(q \notin R(\lambda_+-\operatorname{D}_u G(\mathbf{0} , \alpha))\), and the non-resonance condition, which states that for all \(n \in \Z \setminus \{-1,1\}\), \(n \lambda_+\) is not an eigenvalue of \(\operatorname{D}_u G(\mathbf{0} , \alpha)\). 
Imposing specific conditions on the Fourier modes \(\hat W_1(k), \hat W_2(k), \hat W(k)\), we verify these conditions and thereby demonstrate a local Hopf bifurcation.

Throughout this work we fix a mode \(k^*\ge1\) with
\(\hat W_1(k^*)+\hat W_2(k^*)\neq0\) and set
\begin{equation} \label{eq: alpha_bifurcation_point}
    \alpha_{k^*} :=- \frac{\sigma_1 + \sigma_2}{\bar u (\hat W_1(k^*)+\hat W_2(k^*))} ,
\end{equation}
the unique parameter value at which \(\operatorname{tr}A(k^*)\) vanishes. Whether such a \(k^*\) exists, and which one is selected first, is addressed in
Theorem~\ref{theorem: first_bifurcation}.
Additionally, from now on we call a value of \(\gamma\) \emph{admissible} if it satisfies~\eqref{eq: imaginary_condition} and whenever the critical mode \(k^*\) is clear from the context, we write
\begin{equation} \label{eq: def_S_2}
S_2 := S_2(k^*) =  \frac{\sigma_2 \hat W_1(k^*) - \sigma_1 \hat W_2(k^*)}{\hat W_1(k^*) + \hat W_2(k^*) } \qquad \mathrm{and} \qquad \omega_0:=  \bigl(\bar u^2 \gamma^2 \widehat W(k^*)^2 - S_2^2(k^*)\bigr)^{1/2} 
\end{equation}
where \(S_2\) obviously depends on \(\alpha\) and should be viewed as evaluated at \(\alpha = \alpha_{k^*}\). 

\begin{assumption}[Interaction kernels]\label{ass: kernels}
The kernels satisfy
\begin{equation} \label{eq: cardinality_condition}
     \mathrm{card} \bigl\{k \ge 1 : \hat W_1(k) + \hat W_2(k) = \hat W_1(k^*) + \hat W_2(k^*)\bigr\} = 1
\end{equation}
and, at the parameter value \(\alpha=\alpha_{k^*}\),
\begin{equation} \label{eq: determinante_non_vanishing}
    \det A(k) \neq 0 \qquad \text{for all } k \ge 1 ,
\end{equation}
with \(A(k)\) as in~\eqref{eq: def_A_k}.
\end{assumption}

\begin{assumption}[Non-reciprocity]\label{ass: nonreciprocity}
At the parameter value \(\alpha=\alpha_{k^*}\) the coupling satisfies
\begin{equation} \label{eq: imaginary_condition}
       \bar u^2 \gamma^2 \widehat W(k^*)^2  >  S_2(k^*)^2 .
\end{equation}
\end{assumption}
Notice that~\eqref{eq: imaginary_condition} implies \(\omega_0 > 0\) and \(\hat W(k^*) \neq 0\). 
\begin{lemma}[Spectrum]\label{lem: critical_spectrum}
\; 
\begin{enumerate}[(i)]
  \item Let \(\alpha = \alpha_{k^*}\). Under Assumption~\ref{ass: nonreciprocity}, the block \(A(k^*)\) has the pair of purely
  imaginary eigenvalues
  \begin{equation} \label{eq: eigenvalue_bifurcation}
    \lambda_{\pm}(k^*) = \pm i |k^*|^2 \omega_0 =  \pm i |k^*|^2  \bigl(\bar u^2 \gamma^2 \widehat W(k^*)^2 - S_2^2(k^*)\bigr)^{1/2} ,
  \end{equation}
  and the eigenspace of \(\lambda_+(k^*)\) in the complexification of
  \(\dot L^2_s(\T)\times\dot L^2_s(\T)\) is spanned by \((e_{k^*},\chi_{k^*}e_{k^*})^{\mathrm{T}}\),
  where
  \begin{equation} \label{eq: def_chi}
    \chi_{k^*} := \frac{-\lambda_+(k^*)/|k^*|^2 + S_2(k^*)}{\bar u \gamma \hat W(k^*)} .
  \end{equation}
  Moreover \(|\chi_{k^*}|=1\).
  \item Let \(\alpha = \alpha_{k^*}\). Under Assumption~\ref{ass: kernels} in addition, no eigenvalue of
  \(\operatorname{D}_uG(\mathbf 0,\alpha_{k^*})\) outside the \(k^*\)-block is zero or purely
  imaginary. Consequently \(\lambda_\pm(k^*)\) are simple eigenvalues of
  \(\operatorname{D}_uG(\mathbf 0,\alpha_{k^*})\) and satisfy the non-resonance condition:
  \(n\lambda_+(k^*)\) lies in the resolvent set for every \(n\in\Z\setminus\{-1,1\}\).
  \item The crossing is transversal,
  \(\tfrac{\Id}{\Id\alpha}\Re\lambda_\pm(k^*)\vert_{\alpha=\alpha_{k^*}}
    =-\tfrac{|k^*|^2}{2}\bar u(\hat W_1(k^*)+\hat W_2(k^*))\neq0\).
\end{enumerate}
\end{lemma}

\begin{remark}
    The conditions~\eqref{eq: cardinality_condition} and~\eqref{eq: determinante_non_vanishing} ensure that the eigenvalues of the linearized operator \(  D_u G(\boldsymbol 0, \alpha)\) associated with the critical mode \(k^*\) are simple and satisfy the non-resonance condition.
    The inequality~\eqref{eq: imaginary_condition}  describes that in order to see Hopf bifurcation we require a strong enough non-reciprocity \(\gamma\) if \(\sigma_1 \hat W_2(k^*) \neq \sigma_2 \hat W_1(k^*)\). Otherwise, \(\gamma > 0 \) suffices. Hence, we require one Fourier mode such that the non-reciprocal strength \(\gamma>0\) is large enough in comparison to the diffusion of the system \(\sigma_1 , \sigma_2\). 
\end{remark}

\begin{proof}[Proof of Lemma~\ref{lem: critical_spectrum}]

In order for the eigenvalues to cross the imaginary axis outside the origin, we require the conditions \(\operatorname{tr}(A(k)) = 0\) and \(\det(A(k)) > 0\). At the Fourier mode \(k^*\) the first condition is equivalent to 
\begin{equation*}
    \sigma_1+\sigma_2 +\alpha \bar u \bigl( \hat W_1(k^*) + \hat W_2(k^*) \bigr) = 0. 
\end{equation*}
This condition is satisfied at the bifurcation point \(\alpha_{k^*}\) defined through~\eqref{eq: alpha_bifurcation_point}. 
The second condition \(\det(A(k^*))>0\) at the Fourier mode \(k^*\) is satisfied at the bifurcation point \(\alpha_{k^*}\) if and only if inequality~\eqref{eq: imaginary_condition} holds, i.e. 
\begin{equation*}
    \bar u^2 \gamma^2 \widehat W(k^*)^2  > S_2^2 =\Bigl(  \frac{\sigma_2 \hat W_1(k^*) - \sigma_1 \hat W_2(k^*)}{\hat W_1(k^*) + \hat W_2(k^*) }\Bigr)^2 .
\end{equation*}

Plugging in \(\alpha_{k^*}\) into the eigenvalue formula~\eqref{eq: eigenvalues_abstract} for the \(k^*\)-mode we find 
\begin{equation*} 
    \lambda_{\pm}(k^*) = \pm i |k^*|^2 \Big( - S_2^2 +   \bar u^2 \gamma^2 \widehat W(k^*)^2 \Big)^{\frac{1}{2}}, 
\end{equation*}
while for all other \( k \neq k^*\) we have
\begin{align*}
    \lambda_{\pm}( k) 
    & =-\frac{|k|^2}{2}
\Big( \sigma_1+\sigma_2 +\bar u\alpha_{k^*}\widehat W_1(k) +\bar u\alpha_{k^*}\widehat W_2(k) \Big)
 \pm \frac{i | k|^2}{2} \Big(  4\Big(\sigma_1+\bar u\alpha_{k^*}\widehat W_1(k)\Big)\Big( \sigma_2+\bar u\alpha_{k^*} \widehat W_2( k) \Big)\\
&\quad  + 4\bar u^2 \gamma^2 \widehat W( k)^2 - \Big( \sigma_1+\sigma_2 +\bar u\alpha_{k^*}\widehat W_1(k) +\bar u\alpha_{k^*}\widehat W_2(k) \Big)^2\Big)^{\frac{1}{2}}. 
\end{align*}
With this eigenvalue it is easy to verify that \((e_{k^*}, \chi_{k^*} e_{k^*})^{\mathrm{T}}\) is the corresponding eigenvector. Additionally, since \[\bigg|S_2(k^*) - \frac{\lambda_+(k^*)}{|k^*|^2} \bigg|^2 = S_2^2(k^*) + \omega_0^2 = \bar u^2 \gamma^2 \hat W^2(k^*) ,\]
it follows that \(|\chi_{k^*}| = 1 . \)

From the formula for \(\lambda_\pm(k^*)\) in a small interval around \(\alpha_{k^*}\), where the eigenvalues have non vanishing imaginary part, we have 
\begin{equation*}
    \Re ( \lambda_+(k^*) ) = \frac{\operatorname{tr}(A(k^*))}{2} = -\frac{|k^*|^2}{2}
\Big( \sigma_1+\sigma_2 +\bar u\alpha \widehat W_1(k^*) +\bar u\alpha \widehat W_2(k^*) \Big).
\end{equation*}
Consequently, we obtain 
\begin{equation} \label{eq: real_derivative_eigenvalue}
    \frac{\Id }{\Id \alpha}  \Re ( \lambda_\pm(k^*) )_{|\alpha=\alpha_{k^*}}   = - \frac{|k^*|^2}{2} \bar u (\hat W_1(k^*)+\hat W_2(k^*)) ,
\end{equation}
and 
\begin{equation}\label{eq: imiginary_derivative_eigenvalue}
    \frac{\Id }{\Id \alpha} \Im (  \lambda_\pm(k^*) )_{|\alpha=\alpha_{k^*}}  = \pm
    \frac{|k^*|^2 \bar u S_2 (\hat W_1(k^*)- \hat W_2(k^*))  }{2 \omega_0 }
\end{equation}
which implies that the eigenvalue is crossing the imaginary axis with non-vanishing speed since \(\hat W_1(k^*) + \hat W_2(k^*) \neq 0 \) by the  assumptions of this section, which proves \((iii)\).

Additionally, at the bifurcation point \(\alpha_{k^*} \) we have \(\operatorname{tr}A(k )=0\) if and only if 
\begin{align*}
    &  \sigma_1 +  \bar u \alpha_{k^*} \hat W_1(k)  +  \sigma_2 +  \bar u \alpha_{k^*} \hat W_2(k) = 0 \\ 
\Leftrightarrow\quad  & \hat W_1(k) + \hat W_2(k) = \hat  W_1(k^*)+ \hat W_2(k^*).
\end{align*}
Hence, if conditions~\eqref{eq: cardinality_condition}
and~\eqref{eq: determinante_non_vanishing} hold, then at the bifurcation point \(\alpha_{k^*}\) we have
\[
\operatorname{Re}\big(\lambda_{\pm}(k)\big)
 \neq 0
\]
for all \(k \neq k^*\). Indeed, if \(\lambda_{\pm}\) is non-real, then 
\[
\operatorname{Re}\big(\lambda_{\pm}(k)\big)
=\frac{1}{2} \operatorname{tr}(A(k))
\neq 0.
\]
Otherwise, if \(\lambda_\pm\) is real, then \(\det(A(k))\neq 0\) implies \(\lambda_\pm \neq 0\) via formula~\eqref{eq: eigenvalues_abstract}.
With the spectral decomposition established in Lemma~\ref{lemma: block_spectrum}, we can transfer the analysis of the finite-dimensional matrices \(A(k)\) to the infinite-dimensional operator \(\operatorname{D}_uG\). Since $\omega_0>0$, the eigenvalues $\lambda_\pm(k^*)$ of the $2\times2$ block are distinct, hence algebraically simple in the block, and by the block decomposition (Theorem~\ref{theorem: block_inverse}) also for $\operatorname{D}_uG(\mathbf 0,\alpha_{k^*})$.
Consequently, by the above analysis this implies that the two eigenvalues defined in~\eqref{eq: eigenvalue_bifurcation} are simple. Additionally, the non-resonance condition holds, i.e. \( n \lambda_{\pm} (k^*) \) is not an eigenvalue for all \(n \in \Z  \setminus \{-1, 1\}\) since all other eigenvalues have a non-vanishing real part under the conditions~\eqref{eq: cardinality_condition},~\eqref{eq: determinante_non_vanishing}. This proves \((i)\) and \((ii)\).
\end{proof}

\begin{proposition}[Hopf bifurcation]\label{prop: hopf_bifurcation}
Let Assumptions~\ref{ass: kernels} and~\ref{ass: nonreciprocity} hold. Then a Hopf bifurcation occurs for~\eqref{eq: 2species_problem} at \(\alpha_{k^*}\). 
For every \(\theta\in(0,1)\) there is a continuously differentiable curve
\[
    r \mapsto (u(r),\alpha(r),\kappa(r))
    \in
    \Bigl(
        C^{1+\theta}_{2\pi/\kappa(r)}(\R, \dot L^2_s(\T) \times \dot L^2_s(\T))
        \cap
        C^{\theta}_{2\pi/\kappa(r)}(\R, \dot H^2_s(\T) \times \dot H^2_s(\T))
    \Bigr) \times \R \times \R ,
\]
defined for \(|r|\) small, such that \(u(r)\) is a real \(\tfrac{2\pi}{\kappa(r)}\)-periodic
solution of~\eqref{eq: abstract_problem} at the parameter value \(\alpha(r)\), and
\[
  (u(0),\alpha(0))=(0,\alpha_{k^*}),
  \qquad
  \kappa(0)=\Im\lambda_+(k^*)=|k^*|^2\omega_0 .
\]
\end{proposition}

\begin{proof}
Combining the spectral properties of \(\operatorname{D}_uG\) established in Lemma~\ref{lem: critical_spectrum} with the results in Section~\ref{subsec: spectral_analysis}, which show that \(\operatorname{D}_uG\) generates a compact analytic semigroup for \(t>0\), we conclude that all assumptions of the Hopf bifurcation theorem~\cite[Theorem~I.8.2]{kielh2004} are satisfied. Proposition~\ref{prop: hopf_bifurcation} therefore follows.
\end{proof}
 
An immediate consequence of the result~\cite[Theorem~I.8.2]{kielh2004}, the explicit bifurcation formulas~\cite[Theorem~I.9.1]{kielh2004} and Lemma~\ref{lem: critical_spectrum} is the Taylor expansion of the time periodic solution \(u\). 
\begin{corollary}\label{cor: hopf_profile}
In the setting of Proposition~\ref{prop: hopf_bifurcation}, rescale time by \(s=\kappa(r)t\) and
put \(\tilde u(r)(s):=u(r)(s/\kappa(r))\), so that all solutions on the branch are
\(2\pi\)-periodic. Then there is \(r_0>0\) such that for all \(|r|<r_0\)
\begin{align*}
    \tilde u(r)(s)&=
    2r\,e_{k^*}
    \begin{pmatrix}
        \cos(s) \\
        \Re(\chi_{k^*})\cos(s)-\Im(\chi_{k^*})\sin(s)
    \end{pmatrix}
    +O(r^{2}) \\
    &= 2r\,e_{k^*}
    \begin{pmatrix}
        \cos(s) \\
        \dfrac{S_2}{\bar u\gamma\hat W(k^*)}\cos(s)
        +\dfrac{\omega_0}{\bar u\gamma\hat W(k^*)}\sin(s)
    \end{pmatrix}
    +O(r^{2}),
\end{align*}
the remainder being estimated in
\(C^{1+\theta}_{2\pi}(\R,\dot L^2_s(\T)\times\dot L^2_s(\T))
 \cap C^{\theta}_{2\pi}(\R,\dot H^2_s(\T)\times\dot H^2_s(\T))\) with \(\theta\) as in Proposition~\ref{prop: hopf_bifurcation}.
\end{corollary}

\section{Criticality of the Hopf bifurcation }\label{sec: criticallity_hopf}

The aim of this section is to determine the criticality of the Hopf bifurcation established in Proposition~\ref{prop: hopf_bifurcation}, that is, the direction in which the branch of periodic solutions emerges from the critical parameter \(\alpha_{k^*}\). To explain the two possible scenarios, suppose that \(k^*\in\mathcal T^-\), so that the critical mode loses stability as \(\alpha\) increases through \(\alpha_{k^*}\). A supercritical bifurcation gives rise to periodic solutions for \(\alpha>\alpha_{k^*}\), so that oscillations emerge continuously as the critical mode becomes unstable. In contrast, a subcritical bifurcation produces unstable periodic solutions for \(\alpha<\alpha_{k^*}\), where the critical mode is still stable. In this case, no stable periodic solutions emerge beyond the instability threshold. For \(k^*\in\mathcal T^+\), the roles of the two sides of \(\alpha_{k^*}\) are interchanged. This distinction is particularly relevant for the stability of the emerging periodic solutions. Indeed, as we show in Section~\ref{sec: exchange}, a supercritical Hopf bifurcation gives rise to locally orbitally asymptotically stable periodic solutions, provided all remaining modes are stable at the bifurcation point.

To distinguish between these two scenarios, we rely on classical Hopf bifurcation theory, which relates the criticality to the sign of the real part of the cubic normal-form coefficient \(\Lambda\)~\cite[Section~I.9]{kielh2004}. We first recall the abstract formula for \(\Lambda\) and the precise mathematical definitions of supercritical and subcritical bifurcations in Section~\ref{sec: intro_Sub-vs. Supercriticality}. Since this abstract formula does not directly reveal the criticality for our system~\eqref{eq: 2species_problem}, we then compute \(\Lambda\) explicitly in Proposition~\ref{prop: sub_supercritical}, expressing it in terms of the diffusivities \(\sigma_1,\sigma_2\), the Fourier coefficients of the interaction kernels, the bifurcation point \(\alpha_{k^*}\), and the strength of the non-reciprocal coupling \(\gamma\).

Although the formula for \(\Lambda\) is explicit, its complexity makes it difficult to determine the sign of \(\Re(\Lambda)\) directly. We therefore investigate three particular regimes in which the expression simplifies or can be analysed further, allowing us to derive explicit conditions on the parameters of system~\eqref{eq: 2species_problem} for the bifurcation to be supercritical or subcritical.

Motivated by the Kuramoto interaction kernel (Example~\ref{example: kuramoto}), corresponding to \(W_1=W_2=W=-\cos\), we first consider the case of vanishing second-harmonic coefficients, that is, \(\hat W(2k^*)=\hat W_1(2k^*)=\hat W_2(2k^*)=0\). In this regime, we show in Corollary~\ref{cor: equal_diffusion_super} that, for equal diffusivities \(\sigma_1=\sigma_2\), the Hopf bifurcation is supercritical for every admissible value of \(\gamma\). Interestingly, this conclusion no longer holds when the diffusivities differ. In this case, the bifurcation can also be subcritical, with its criticality depending in particular on the ratio \(\hat W_2(k^*)/\hat W_1(k^*)\) of the Fourier coefficients at the critical mode.

In the second regime, we allow for non-vanishing second-harmonic coefficients but restrict ourselves to the matching case \(\sigma_1=\sigma_2\) and \(W_1=W_2=W\). Thus, the two species have identical diffusivities and interaction kernels, differing only in the sign of their cross-interactions. A particular example is given by kernels involving both the first and second Fourier modes, which describe systems with polar and nematic alignment interactions (see Example~\ref{example: cos}).

In this setting, Corollary~\ref{cor: lambda_matching} provides an explicit polynomial \(N_R\) such that the Hopf bifurcation is supercritical if and only if \(N_R(\omega_0^2)>0\). In particular, we derive a parameter regime where \(N_R\) has a unique positive root, which identifies a threshold \(\gamma^*>0\) such that the bifurcation is subcritical for \(\gamma<\gamma^*\) and supercritical for \(\gamma>\gamma^*\). Hence, sufficiently strong non-reciprocity can change the criticality of the Hopf bifurcation and lead to the emergence of stable periodic solutions, provided the remaining modes are stable at the bifurcation point.
Remarkably, this dependence on the coupling strength is absent in the corresponding reciprocal system, where the criticality is independent of \(\gamma\). Our result therefore shows that non-reciprocity does not merely give rise to oscillations, but fundamentally changes how the strength of the interaction determines the nonlinear dynamics of the emerging branch (see Remark~\ref{rem: matching_criticality}).

Motivated by this observation, we finally investigate the influence of \(\gamma\) on the criticality of the Hopf bifurcation in the general case, without assuming equal diffusivities or matching interaction kernels. 
Our analysis shows that sufficiently strong non-reciprocity can enforce supercriticality under an explicit condition on the self-interaction kernels. More precisely, Corollary~\ref{cor: high_gamma_cond} establishes that, for sufficiently large \(\gamma\), the Hopf bifurcation is supercritical whenever
\begin{equation*}
    \frac{\hat W_1(2k^*)+\hat W_2(2k^*)}
         {\hat W_1(k^*)+\hat W_2(k^*)}
    < \frac{6\sqrt{2}-5}{4},
\end{equation*}
independently of the ratio \(\iota=\hat W(2k^*)/\hat W(k^*) \notin \{0,\pm 1/2\}\) of the cross-interaction kernel. Thus, even beyond the matching regime, sufficiently strong non-reciprocity provides a mechanism for the emergence of a supercritical Hopf bifurcation.

\subsection{The cubic normal-form coefficient} \label{sec: intro_Sub-vs. Supercriticality}
Throughout this section we write
\begin{equation}\label{eq: def_phi0_qstar}
 q_{k^*}:=\begin{pmatrix} e_{k^*}\\ \chi_{k^*}e_{k^*}\end{pmatrix},
  \qquad
  q_{k^*}^*:=\frac{1}{1-\overline{\chi}_{k^*}^{\,2}}\begin{pmatrix} e_{k^*}\\ -\overline{\chi}_{k^*}\,e_{k^*}\end{pmatrix},
\end{equation}
and we use the complex pairing
\(\langle \boldsymbol f,\boldsymbol g\rangle:=\int_{\T}\bigl(f_1\overline{g_1}+f_2\overline{g_2}\bigr)\Id x\)
on the complexification of \(\dot L^2_s(\T)\times\dot L^2_s(\T)\). By Lemma~\ref{lem: critical_spectrum},
\(q_{k^*}\) spans the eigenspace of \(\operatorname{D}_uG(\mathbf 0,\alpha_{k^*})\) for \(\lambda_+(k^*)\).
Using~\eqref{eq: def_chi} one checks
\(A(k^*)^{\mathrm T}(1,-\overline{\chi}_{k^*})^{\mathrm T}=\overline{\lambda_+(k^*)}\,(1,-\overline{\chi}_{k^*})^{\mathrm T}\),
so \(q_{k^*}^*\) spans the eigenspace of the adjoint operator \(\operatorname{D}_uG(\mathbf 0,\alpha_{k^*})^*\) for
\(\overline{\lambda_+(k^*)}\). Since \(\omega_0>0\) we have \(\chi_{k^*}^2\neq1\), so \(q_{k^*}^*\) is well defined, and
\(|\chi_{k^*}|=1\) gives the normalisation \(\langle q_{k^*},q_{k^*}^*\rangle=1\).
Following~\cite[Section~I.9]{kielh2004} we define the cubic normal form coefficient
\begin{equation}\label{def: Lambda}
\begin{aligned}
\Lambda
&:= -\Biggl\langle D_{uu}^2 G(\mathbf{0}, \alpha_{k^*})\Bigl[
\bar q_{k^*} ,
\bigl(2\lambda_+(k^*)-D_uG(\mathbf{0},\alpha_{k^*})\bigr)^{-1}
\Bigr. \Bigl. D_{uu}^2G(\mathbf{0},\alpha_{k^*})\Bigl[
q_{k^*} ,
q_{k^*}
\Bigr]\Bigr]
,q_{k^*}^*
\Biggr\rangle \\
&\quad + 2\Biggl\langle D_{uu}^2 G(\mathbf{0}, \alpha_{k^*})\Bigl[
q_{k^*},
D_uG(\mathbf{0},\alpha_{k^*})^{-1}
\Bigr.\Bigl. D_{uu}^2G(\mathbf{0},\alpha_{k^*})\Bigl[
q_{k^*},
\bar q_{k^*}
\Bigr]\Bigr]
,q_{k^*}^*
\Biggr\rangle.
\end{aligned}
\end{equation}
Through the formula~\cite[Theorem~I.9.1]{kielh2004}, 
\begin{equation} \label{eq: alpha_second_derivative}
    \alpha''(0) = \frac{ \Re(\Lambda)}{ \frac{\Id }{\Id \alpha}  \Re ( \lambda_\pm(k^*) )_{|\alpha=\alpha_{k^*}} } = -
   \frac{2}{\bar u ( \hat W_1(k^*) + \hat W_2(k^*)) |k^*|^2}  \Re(\Lambda)
\end{equation}
for the second derivative of the bifurcation curve, the real part of \(\Lambda\) fixes the
direction in which the branch opens.

By~\eqref{eq: alpha_second_derivative} and the transversality condition~\eqref{eq: real_derivative_eigenvalue} we notice that the sign of \(\Re(\Lambda)\) determines whether the branch opens towards the side of
\(\alpha_{k^*}\) on which the homogeneous state is unstable, regardless of the direction in
which \(\alpha\) destabilizes. We therefore call the bifurcation supercritical if
\(\Re(\Lambda)>0\) and subcritical if \(\Re(\Lambda)<0\). Equivalently,
by~\eqref{eq: alpha_second_derivative}, the bifurcation is supercritical if and only if
\(\alpha''(0)>0\) for \(k^*\in\mathcal T^-\), and if and only if \(\alpha''(0)<0\) for
\(k^*\in\mathcal T^+\). At the first bifurcation point \(\alpha_+\) of
Theorem~\ref{theorem: first_bifurcation} one necessarily has \(k^*\in\mathcal T^-\), so there the
two formulations agree.

In the non-reciprocal setting, the parameter $\gamma$ enters $A(k)$ only through the off-diagonal entries and therefore not through its trace, so the critical parameter $\alpha_{k^*}$ itself is 
independent of $\gamma$: at linear order non-reciprocity is visible only in the frequency $\omega_0$ and the corresponding eigenvector. Its effect primarily appears at cubic order, in $\Lambda$. 
The computations of this section are consequently the main place where the strength of the non-reciprocal coupling can be compared quantitatively with the reciprocal theory~\cite{Carrillo2026}. 

Moreover, the cubic normal-form coefficient \(\Lambda\) governs also the second derivative of the frequency \(\kappa(r)\) via 
\begin{equation} \label{eq: def_frequency}
\kappa''(0) = \frac{  \frac{\Id }{\Id \alpha}  \Im ( \lambda_\pm(k^*) )_{|\alpha=\alpha_{k^*}} \Re(\Lambda)}{ \frac{\Id }{\Id \alpha}  \Re ( \lambda_\pm(k^*) )_{|\alpha=\alpha_{k^*}} } - \Im (\Lambda)= - \frac{S_2(k^*) (\hat W_1(k^*)- \hat W_2(k^*)) }{(\hat W_1(k^*) + \hat W_2(k^*))\omega_0 }  \Re(\Lambda)- \Im (\Lambda) 
\end{equation}
so that the period of the bifurcating orbits is not \(2\pi/(|k^*|^2\omega_0)\) but depends on their amplitude already at second order.

\subsection{The cubic normal-form coefficient for system~\eqref{eq: 2species_problem}} \label{sec: Sub-vs. Supercriticality}
We compute the abstract coefficient \(\Lambda\) given by~\eqref{def: Lambda} in terms of parameters of the system~\eqref{eq: 2species_problem}.

\begin{proposition}\label{prop: sub_supercritical}
Let the assumptions of Proposition~\ref{prop: hopf_bifurcation} hold true. 

Then \(\Lambda\) can be computed to equal 
\begin{align}
\begin{split} \label{eq: def_Lambda}
\Lambda&= 
-\frac{  |k^*|^2}{2 \sqrt \pi(1-\chi_{k^*}^2)}   \Bigg(
b_1 \frac{-\sigma_1+i \omega_0}{\bar u} -2  \alpha_{k^*}  b_1 \hat W_1(2k^*) -2 b_2 \gamma \hat W(2k^*) \\
&\quad -  b_2 \frac{-\sigma_2+i \omega_0}{\bar u} + 2 b_2 \alpha_{k^*}  \hat W_2(2k^*) -2 b_1\gamma \hat W(2k^*)   \\
&\quad \quad - 2   \Big(    - \tilde b_1 \frac{\sigma_1 + i \omega_0}{\bar u}    -2  \alpha_{k^*}   \tilde  b_1 \hat W_1(2k^*) -2 \tilde b_2 \gamma \hat W(2k^*) \\
&\quad +  \chi_{k^*}^2  \tilde{b}_2 \frac{\sigma_2+ i \omega_0}{\bar u} + 2 \chi_{k^*}^2  \tilde  b_2 \alpha_{k^*}  \hat W_2(2k^*) -2 \chi_{k^*}^2  \tilde b_1\gamma \hat W(2k^*) \Big)  \bigg) 
\end{split}
\end{align}
with \(b_i\), \(\tilde b_i\) given by 
\begin{align*}
    b_1&=  \frac{2|k^*|^2 \bigg( \Big( 2\lambda_+(k^*) + |2k^*|^2(\sigma_2+\bar u\alpha_{k^*} \widehat W_2(2k^*)) \Big) \frac{\sigma_1+i\omega_0}{\bar u} - |2k^*|^2 \gamma \chi_{k^*}^2 \widehat W(2k^*) (\sigma_2+i\omega_0) \bigg) }{\sqrt \pi \det(2\lambda_+(k^*)- A_{2k^*})} , \\
    b_2&=  \frac{2|k^*|^2 \bigg( |2k^*|^2 \gamma \widehat W(2k^*) (\sigma_1+i\omega_0) + \Big[ 2\lambda_+(k^*) + |2k^*|^2(\sigma_1+\bar u\alpha_{k^*} \widehat W_1(2k^*)) \Big] \chi_{k^*}^2 \frac{\sigma_2+i\omega_0}{\bar u} \bigg) }{\sqrt \pi \det(2\lambda_+(k^*)- A_{2k^*})}, \\
    \tilde b_1&=  - \frac{8|k^*|^4 }{\sqrt \pi \det(A_{2k^*})}\bigg( (\sigma_2+\bar u\alpha_{k^*} \widehat W_2(2k^*))  \frac{\sigma_1 }{\bar u} -  \bar u\gamma \widehat W(2k^*)  \frac{\sigma_2}{\bar u}\bigg) , \\
    \tilde b_2&=  - \frac{8|k^*|^4 }{\sqrt \pi \det(A_{2k^*})}\bigg( \bar u\gamma \widehat W(2k^*) \frac{\sigma_1}{\bar u} +  (\sigma_1+\bar u\alpha_{k^*} \widehat W_1(2k^*)) \frac{\sigma_2}{\bar u} \bigg) . 
\end{align*}

\end{proposition}

\begin{remark}
In this generality it is very hard to provide algebraic conditions for the sign of \(\Re(\Lambda)\) to be positive or negative. Numerical simulation of the coefficients show that the closer the \(2k^*\)-frequency is to the \(k^*\)-frequency the more likely it is for the sign of \(\Re(\Lambda)\) to switch (even multiple times). Hence, this means that if we have resonance of the \(k^*\) and \(2k^*\) mode numerical simulations suggest  that we can actually not make any guess whether the Hopf bifurcation is sub or supercritical. However, if we let \(|\hat W (2k^*)|\) be small we are going to verify in Corollary~\ref{cor: sub_subcritical_small_case} that if \(\sigma_1=\sigma_2\) the sign of \(\alpha''(0)\) corresponds to the sign of \(-(\hat W_1(k^*)+\hat W_2(k^*))\). Hence, in combination with transversality~\eqref{eq: real_derivative_eigenvalue} and the Principle of Exchange of Stability this will demonstrate that the Hopf bifurcation at the critical point always produces a stable orbit and will be proven later on in Theorem~\ref{theorem: first_bifurcation}. Additionally, notice that in the matching case \(\sigma=\sigma_1=\sigma_2\) and \(\hat W_1 = \hat W_2\) the second derivative of the frequency \(\kappa(r)\) only depends on the imaginary part of \(\Lambda\). 
\end{remark}

\begin{proof}[Proof of Proposition~\ref{prop: sub_supercritical}]
Based on the formula~\eqref{def: Lambda} we require the second derivative of our abstract map \(G\), which we computed in~\eqref{eq: second_derivative}. 
In particular, using \(V * e_k= \sqrt{\pi} \langle V, e_k \rangle  e_k\) for any even function \(V\), \(\nabla e_k e_k = \frac{1}{2} \nabla e_k^2 \) and \(\frac12 \Delta e_k^2 =  -\tfrac{|k|^2}{\sqrt \pi} e_{2k}\), we obtain for \(k \ge 1 \),
\begin{equation*}
        D_{uu}^2 G(\boldsymbol 0, \alpha)[\boldsymbol{e_k}, \boldsymbol{e_k}]
        =
        \begin{pmatrix}
             -\frac{2|k|^2}{\sqrt \pi} \big( \alpha \hat W_1(k) e_{2k} + \gamma \hat W(k) e_{2k} \big) \\
              -\frac{2|k|^2}{\sqrt \pi} \big( \alpha \hat W_2(k) e_{2k} - \gamma \hat W(k) e_{2k} \big) 
        \end{pmatrix}
\end{equation*}
In particular, recall
\begin{equation*}
    \chi_{k^*} = \frac{-\lambda_+(k^*)/|k^*|^2 - ( \sigma_1 + \bar u \alpha_{k^*}  \hat W_1(k^*)) }{\bar u \gamma \hat W(k^*)} =\frac{-\lambda_+(k^*)/|k^*|^2 + S_2}{\bar u \gamma \hat W(k^*)}  
\end{equation*}
and 
\begin{align*}
        D_{uu}^2 G(\boldsymbol 0, \alpha_{k^*})[q_{k^*},q_{k^*}]
        &= 
        \begin{pmatrix}
              -\frac{2|k^*|^2}{\sqrt \pi} \big( \alpha_{k^*} \hat W_1(k^*) e_{2k^*} + \chi_{k^*} \gamma \hat W(k^*) e_{2k^*} \big) \\
              -\frac{2|k^*|^2}{\sqrt \pi}\big( \chi_{k^*}^2 \alpha_{k^*} \hat W_2(k^*) e_{2k^*} - \chi_{k^*}\gamma \hat W(k^*) e_{2k^*} \big) 
        \end{pmatrix} \\
        &= \frac{2|k^*|^2}{\sqrt \pi}
          \begin{pmatrix}
               \frac{\sigma_1+i \omega_0}{\bar u} \\
               \chi_{k^*}^2 \frac{\sigma_2+i \omega_0}{\bar u}
        \end{pmatrix} e_{2k^*},
\end{align*}
where we applied \(|\chi_{k^*}| =1\) and formulas from Table~\ref{eq: aux_special_case_lyapunov}. 
Next, let us invert the operator \(2\lambda_+(k^*)- A(2k^*) \). Notice that the inverse exists by Lemma~\ref{lem: critical_spectrum}.
In comparison to the definition of \(\Lambda\) in~\eqref{def: Lambda} it is enough to compute the inverse on each block. Indeed, we just need to check the condition of Theorem~\ref{theorem: block_inverse}. By the proof of Corollary~\ref{cor: A0_invertible} we know that there exists a \(c > 0\) and \(R \in\N \) such that \(c|k|^2 \le \sigma_{min}(A(k)):= \min_{|v|=1} |A(k) v|  \) for all \(|k| \ge R\). Thus, 
\begin{equation*}
  \min_{|v|=1} |(2\lambda_+(k^*) - A(k))  v| \ge  \sigma_{min}(A(k)) - 2|\lambda_+(k^*)| \ge 2|\lambda_+(k^*)| > 0 
\end{equation*}
for all \(|k| \ge \max(R, 2  \sqrt{|\lambda_+(k^*)| / c } ) \). Moreover, by the non-resonance condition we know that \(2\lambda_+(k^*)- A(k) \) is invertible for finitely many \(|k| \le \max(R, 2  \sqrt{|\lambda_+(k^*)| / c })\). Consequently, we proved the condition of Theorem~\ref{theorem: block_inverse} and we can compute the inverse block-wise. 
We obtain
{\small
\begin{align*}
    &(2\lambda_+(k^*)- A(2k^*))^{-1}\\
&\quad = \frac{1}{\det(2\lambda_+(k^*)- A_{2k^*})}
    \begin{pmatrix}
    2\lambda_+(k^*) +   |2k^*|^2 (\sigma_2+\bar u\alpha_{k^*} \widehat W_2(2k^*) )
&
-|2k^*|^2 \bar u\gamma \widehat W(2k^*)
\\
|2k^*|^2\bar u\gamma \widehat W(2k^*)
&
 2\lambda_+(k^*)  +  |2k^*|^2( \sigma_1+\bar u\alpha_{k^*} \widehat W_1(2k^*))
    \end{pmatrix} .
\end{align*}
}
We compute 
{\small
\begin{align*}
  & (2\lambda_+(k^*)- A(2k^*))^{-1}  D_{uu}^2 G(\boldsymbol 0, \alpha_{k^*})[q_{k^*},q_{k^*}]
      =  \frac{2|k^*|^2 e_{2k^*}}{\sqrt \pi \det(2\lambda_+(k^*)- A_{2k^*})} \\
   &\quad \cdot  \begin{pmatrix}
\Big[ 2\lambda_+(k^*) + |2k^*|^2(\sigma_2+\bar u\alpha_{k^*} \widehat W_2(2k^*)) \Big] \frac{\sigma_1+i\omega_0}{\bar u} - |2k^*|^2 \bar u\gamma \chi_{k^*}^2 \widehat W(2k^*) \frac{\sigma_2+i\omega_0}{\bar u}\\
\\
|2k^*|^2\bar u\gamma \widehat W(2k^*) \frac{\sigma_1+i\omega_0}{\bar u} + \Big[ 2\lambda_+(k^*) + |2k^*|^2(\sigma_1+\bar u\alpha_{k^*}\widehat W_1(2k^*)) \Big] \chi_{k^*}^2 \frac{\sigma_2+i\omega_0}{\bar u}
\end{pmatrix}
\end{align*}
}

Hence, there exists some constants \(b_1, b_2 \in \C\) depending on all parameters such that 
\begin{equation*}
     (2\lambda_+(k^*)- A(2k^*))^{-1} D_{uu}^2 G(\boldsymbol 0, \alpha_{k^*})[q_{k^*}, q_{k^*}] = 
     \begin{pmatrix}
         b_1 e_{2k^*} \\
         b_2 e_{2k^*}
     \end{pmatrix}.
\end{equation*}
Utilizing \(F_1:= \nabla \cdot (\nabla e_{k^*} e_{2k^*} ) =  -\frac{3}{2\sqrt{\pi}} |k^*|^2 e_{3k^*} + \frac{1}{2\sqrt{\pi}}|k^*|^2 e_{k^*} \) and \(F_2:= \nabla \cdot (\nabla e_{2k^*} e_{k^*} ) = -\frac{3}{\sqrt \pi}|k^*|^2 e_{3k^*} - \frac{1}{\sqrt \pi}|k^*|^2 e_{k^*}  \), \(|\chi_{k^*}| =1\) and~\eqref{eq: aux_special_case_lyapunov} we find 
\begin{align*}
          &  D_{uu}^2 G(\boldsymbol 0, \alpha_{k^*})[ \begin{pmatrix}
            e_{k^*} \\ \bar \chi_{k^*} e_{k^*}
        \end{pmatrix},   \begin{pmatrix}
         b_1 e_{2k^*} \\
         b_2 e_{2k^*}
     \end{pmatrix}]\\
        &\quad = 
        \begin{pmatrix}
          \big( \alpha_{k^*} b_1 \hat W_1(k^*)  + b_1\bar \chi_{k^*} \gamma \hat W(k^*)  \big)  F_1  +\big(  \alpha_{k^*} b_1 \hat W_1(2k^*) + b_2 \gamma \hat W(2k^*)\big)  F_2 \\
        \big( \bar \chi_{k^*} b_2 \alpha_{k^*} \hat W_2(k^*)- b_2 \gamma \hat W(k^*) \big) F_1  +\big( \bar \chi_{k^*} b_2 \alpha_{k^*} \hat W_2(2k^*) - b_1 \bar \chi_{k^*} \gamma \hat W(2k^*)  \big) F_2
        \end{pmatrix} \\
        &\quad = 
           \begin{pmatrix}
           b_1 \frac{-\sigma_1+i \omega_0}{\bar u} F_1  + \big(  \alpha_{k^*} b_1 \hat W_1(2k^*) + b_2 \gamma \hat W(2k^*)\big)  F_2 \\
        \bar \chi_{k^*} b_2 \frac{-\sigma_2+i \omega_0}{\bar u} F_1  +\big( \bar \chi_{k^*} b_2 \alpha_{k^*} \hat W_2(2k^*) - b_1 \bar \chi_{k^*} \gamma \hat W(2k^*)  \big) F_2
        \end{pmatrix}.
\end{align*}
We pair this vector with \(q_{k^*}^*\) from~\eqref{eq: def_phi0_qstar}. Since \(q_{k^*}^*\) is supported on
\(e_{k^*}\) and \(\langle e_{3k^*},e_{k^*}\rangle=0\), the \(e_{3k^*}\)-terms do not contribute.
The eigenvector \(q_{k^*}^*\) is supported on \(e_{k^*}\). Hence, we can disregard \(e_{3k*}\) since \(\langle e_{3k^*}, e_{k*} \rangle =0\).
Utilizing again \(\chi_{k^*} \bar \chi_{k^*} = 1 \), we obtain 
\begin{align*}
& \left \langle D_{uu}^2 G(\boldsymbol 0, \alpha_{k^*})[ \begin{pmatrix}
            e_{k^*} \\ \bar \chi_{k^*} e_{k^*}
        \end{pmatrix},   \begin{pmatrix}
         b_1 e_{2k^*} \\
         b_2 e_{2k^*}
     \end{pmatrix}] ,  q_{k^*}^* 
     \right \rangle  \\
&= 
\frac{|k^*|^2}{2 \sqrt \pi(1-\chi_{k^*}^2)}    
\bigg(  b_1 \frac{-\sigma_1+i \omega_0}{\bar u} -2  \alpha_{k^*} b_1 \hat W_1(2k^*) -2 b_2 \gamma \hat W(2k^*) \\
&\quad - \big( b_2 \frac{-\sigma_2+i \omega_0}{\bar u} - 2 b_2 \alpha_{k^*} \hat W_2(2k^*) +2 b_1\gamma \hat W(2k^*) \big)  
\bigg)  \\
&= 
\frac{|k^*|^2}{2 \sqrt \pi(1-\chi_{k^*}^2)}    
\bigg(  b_1 \frac{-\sigma_1+i \omega_0}{\bar u} -2  \alpha_{k^*} b_1 \hat W_1(2k^*) -2 b_2 \gamma \hat W(2k^*) \\
&\quad -  b_2 \frac{-\sigma_2+i \omega_0}{\bar u} + 2 b_2 \alpha_{k^*} \hat W_2(2k^*) -2 b_1\gamma \hat W(2k^*) \big)  
\bigg) 
\end{align*}

Similarly, we find 
\begin{equation*}
        D_{uu}^2 G(\boldsymbol 0, \alpha_{k^*})[q_{k^*}
        ,\bar q_{k^*}]
        = 
        \begin{pmatrix}
              -\frac{2|k^*|^2}{\sqrt \pi} \big( \alpha_{k^*} \hat W_1(k^*) e_{2k^*} + \Re( \chi_{k^*} ) \gamma \hat W(k^*) e_{2k^*} \big) \\
              -\frac{2|k^*|^2}{\sqrt \pi}\big(  \alpha_{k^*} \hat W_2(k^*) e_{2k^*} - \Re(\chi_{k^*})\gamma \hat W(k^*) e_{2k^*} \big) .
        \end{pmatrix}
\end{equation*}
As previously, \(A(2k^*)^{-1}\) exists by Lemma~\ref{lem: critical_spectrum} and can be computed on each block by Corollary~\ref{cor: A0_invertible},
which implies 
{\small
\begin{align*}
  &  A(2k^*)^{-1}  D_{uu}^2 G(\boldsymbol 0, \alpha_{k^*})[q_{k^*}, \bar q_{k^*}]
      =   \frac{8|k^*|^4 e_{2k^*}}{\sqrt \pi \det(A_{2k^*})} \\
   &\quad \cdot  \begin{pmatrix}
(\sigma_2+\bar u\alpha_{k^*}\widehat W_2(2k^*)) \big(\alpha_{k^*} \widehat W_1(k^*) + \Re( \chi_{k^*})  \gamma \widehat W(k^*)\big) -  \bar u\gamma \widehat W(2k^*) \big( \alpha_{k^*} \widehat W_2(k^*) - \Re( \chi_{k^*}) \gamma \widehat W(k^*)\big) \\
\\
\bar u\gamma \widehat W(2k^*) \big(\alpha_{k^*} \widehat W_1(k^*) + \Re(\chi_{k^*}) \gamma \widehat W(k^*)\big) +  (\sigma_1+\bar u\alpha_{k^*} \widehat W_1(2k^*)) \big(  \alpha_{k^*} \widehat W_2(k^*) - \Re( \chi_{k^*})\gamma \widehat W(k^*)\big)
\end{pmatrix} \\
&= - \frac{8|k^*|^4 e_{2k^*}}{\sqrt \pi \det(A_{2k^*})}
\begin{pmatrix}
(\sigma_2+\bar u\alpha_{k^*} \widehat W_2(2k^*))  \frac{\sigma_1 }{\bar u} -  \sigma_2 \gamma \widehat W(2k^*) 
\\
    \sigma_1 \gamma \widehat W(2k^*)  +  (\sigma_1+\bar u\alpha_{k^*} \widehat W_1(2k^*)) \frac{\sigma_2}{\bar u} 
\end{pmatrix}
\end{align*}
}
Here we used~\eqref{eq: aux_real_part_chi_gamma2} and~\eqref{eq: aux_real_part_chi_gamma3} in the last step.
Consequently, we have found the \(\tilde b_1, \tilde b_2\). 
Further, using~\eqref{eq: aux_special_case_lyapunov} and \(|\chi_{k^*}|= 1\), we have 
\begin{align*}
          &  D_{uu}^2 G(\boldsymbol 0, \alpha_{k^*})[ \begin{pmatrix}
            e_{k^*} \\ \chi_{k^*} e_{k^*}
        \end{pmatrix},   \begin{pmatrix}
         \tilde{b}_1 e_{2k^*} \\
         \tilde{b}_2 e_{2k^*}
     \end{pmatrix}]\\
        &\quad = 
        \begin{pmatrix}
          \big( \alpha_{k^*} \tilde{b}_1 \hat W_1(k^*)  + \tilde{b}_1 \chi_{k^*} \gamma \hat W(k^*)  \big) F_1 +\big(  \alpha_{k^*} \tilde{b}_1 \hat W_1(2k^*) + \tilde{b}_2 \gamma \hat W(2k^*)\big) F_2 \\
        \big( \chi_{k^*} \tilde{b}_2 \alpha_{k^*} \hat W_2(k^*)- \tilde{b}_2 \gamma \hat W(k^*) \big) F_1  +\big(  \chi_{k^*} \tilde{b}_2 \alpha_{k^*} \hat W_2(2k^*) - \tilde{b}_1  \chi_{k^*} \gamma \hat W(2k^*)  \big)  F_2
        \end{pmatrix}\\
    &\quad = 
    \begin{pmatrix}
          - \tilde b_1 \frac{\sigma_1 + i \omega_0}{\bar u}   F_1 +\big(  \alpha_{k^*} \tilde{b}_1 \hat W_1(2k^*) + \tilde{b}_2 \gamma \hat W(2k^*)\big)  F_2 \\
        - \chi_{k^*} \tilde{b}_2 \frac{\sigma_2+ i \omega_0}{\bar u} F_1  +\big(  \chi_{k^*} \tilde{b}_2 \alpha_{k^*} \hat W_2(2k^*) - \tilde{b}_1  \chi_{k^*} \gamma \hat W(2k^*)  \big) F_2
        \end{pmatrix}
\end{align*}
and
\begin{align*}
&\left \langle D_{uu}^2 G(\boldsymbol 0, \alpha_{k^*})[ \begin{pmatrix}
            e_{k^*} \\ \chi_{k^*} e_{k^*}
        \end{pmatrix},   \begin{pmatrix}
         \tilde b_1 e_{2k^*} \\
         \tilde b_2 e_{2k^*}
     \end{pmatrix}] ,  q_{k^*}^* 
     \right \rangle  \\
&= 
\frac{|k^*|^2}{2 \sqrt \pi(1-\chi_{k^*}^2)} 
\bigg(  - \tilde b_1 \frac{\sigma_1 + i \omega_0}{\bar u}    -2  \alpha_{k^*}  \tilde  b_1 \hat W_1(2k^*) -2 \tilde b_2 \gamma \hat W(2k^*) \\
&\quad - \big( - \chi_{k^*}^2  \tilde{b}_2 \frac{\sigma_2+ i \omega_0}{\bar u} -2 \chi_{k^*}^2  \tilde  b_2 \alpha_{k^*} \hat W_2(2k^*) +2 \chi_{k^*}^2  \tilde b_1\gamma \hat W(2k^*) \big)  
\bigg) \\
&= 
\frac{|k^*|^2}{2 \sqrt \pi(1-\chi_{k^*}^2)} 
\bigg(  - \tilde b_1 \frac{\sigma_1 + i \omega_0}{\bar u}    -2  \alpha_{k^*}  \tilde  b_1 \hat W_1(2k^*) -2 \tilde b_2 \gamma \hat W(2k^*) \\
&\quad +  \chi_{k^*}^2  \tilde{b}_2 \frac{\sigma_2+ i \omega_0}{\bar u} + 2 \chi_{k^*}^2  \tilde  b_2 \alpha_{k^*} \hat W_2(2k^*) -2 \chi_{k^*}^2  \tilde b_1\gamma \hat W(2k^*)   
\bigg) 
\end{align*}

Combining the above equations we arrive at the formula of Proposition~\ref{prop: sub_supercritical}. 
\end{proof}

\subsection{Kernels with vanishing or small second harmonic}\label{sec: vanishing_case}

Although Proposition~\ref{prop: sub_supercritical} provides an explicit formula for the cubic Hopf normal-form coefficient \(\Lambda\), its dependence on the Fourier coefficients at the second harmonic \(2k^*\) makes the sign of \(\Re(\Lambda)\) difficult to characterize. In this subsection, we consider the case of vanishing second-harmonic coefficients, that is,
\(\hat W_1(2k^*)=\hat W_2(2k^*)=\hat W(2k^*)=0\).
Under this assumption, Corollary~\ref{cor: real_part_vanishing_case} shows that the sign of \(\Re(\Lambda)\) is determined by an explicit quadratic polynomial \(P(\omega_0^2)\). Corollary~\ref{cor: sub_subcritical_vanishing_case} then translates this criterion into conditions involving the ratio \(\hat W_2(k^*)/\hat W_1(k^*)\) of the Fourier coefficients at the critical mode. Finally, Corollary~\ref{cor: sub_subcritical_small_case} shows that the resulting classification persists for sufficiently small second-harmonic coefficients.

An important consequence is that, for equal diffusivities \(\sigma_1=\sigma_2\), the Hopf bifurcation is supercritical for every admissible value of \(\gamma\). This is consistent with the corresponding result for the reciprocal system~\cite[Theorem~2.6]{Carrillo2026}. However, as we shall see in Section~\ref{sec: high_gamma}, this conclusion does not extend to general interaction kernels: non-vanishing second-harmonic coefficients can introduce a dependence of the criticality on the strength of the non-reciprocal coupling \(\gamma\).

\begin{corollary}\label{cor: real_part_vanishing_case}
Let the assumptions of Proposition~\ref{prop: hopf_bifurcation} hold and let \(\hat W_1(2k^*) = \hat W_2(2k^*) =\hat W(2k^*)=0 \). Then 
\begin{equation}\label{eq: ReLambda_two_term}
\Re(\Lambda)
=
\frac{3|k^*|^2}{4\pi\bar u^2}
\left(
\frac{\sigma_1\bigl(\omega_0^2+2\sigma_1^2+\sigma_1 S_2\bigr)}{\omega_0^2+4\sigma_1^2}
+
\frac{\sigma_2\bigl(\omega_0^2+2\sigma_2^2-\sigma_2 S_2\bigr)}{\omega_0^2+4\sigma_2^2}
\right),
\end{equation}
and after combination over a common denominator,
\begin{equation}\label{eq: ReLambda_P_form}
\Re(\Lambda)
=
\frac{3|k^*|^2(\sigma_1+\sigma_2)\,P(\omega_0^2)}{4\pi\bar u^2\,(\omega_0^2+4\sigma_1^2)(\omega_0^2+4\sigma_2^2)},
\qquad
P(x):=x^2+c_1x+8\sigma_1^2\sigma_2^2,
\end{equation}
with
\begin{equation}\label{eq: c1_def}
c_1 := 2(\sigma_1^2 +\sigma_1\sigma_2+\sigma_2^2) + S_2 (\sigma_1 - \sigma_2) .
\end{equation}
Moreover, we have 
\begin{equation}\label{eq: ImLambda_final}
\Im(\Lambda)=\frac{3|k^*|^2}{4\pi\bar u^2\, \omega_0 (\omega_0^2+4\sigma_1^2)(\omega_0^2+4\sigma_2^2)}
\Bigl(\,\omega_0^2\bigl((\sigma_1^2+\sigma_2^2)\omega_0^2+8\sigma_1^2\sigma_2^2\bigr)+S_2(\sigma_2-\sigma_1)\,Q(\omega_0^2)\Bigr)
\end{equation}
with 
\begin{equation}\label{eq: def_Q}
Q(x):=x^2+2\bigl(\sigma_1^2-\sigma_1\sigma_2+\sigma_2^2\bigr)x+8\sigma_1^2\sigma_2^2.
\end{equation}
\end{corollary}

\begin{proof}
Recall that so far we computed general formulas for any symmetric interaction potential \(W_1, W_2,W\). In this step we reduce the problem to the setting in Corollary~\ref{cor: real_part_vanishing_case}, i.e. \(\hat W(2k^*) = \hat W_1(2k^*) =  \hat W_2(2k^*) = 0 \). 
Hence, utilizing~\eqref{eq: aux_real_part_chi_gamma2},~\eqref{eq: aux_real_part_chi_gamma3}  
we find 
\begin{equation*}
 D_{uu}^2 G(\boldsymbol 0, \alpha_{k^*})[q_{k^*}, q_{k^*}]
        = \frac{2|k^*|^2}{\sqrt \pi \bar u } \begin{pmatrix}
               \sigma_1+i \omega_0 \\
               \chi_{k^*}^2 (\sigma_2+i \omega_0)
        \end{pmatrix} e_{2k^*}, 
        \quad 
        D_{uu}^2 G(\boldsymbol 0, \alpha_{k^*})[q_{k^*},\bar q_{k^*}]
        = \frac{2|k^*|^2}{\sqrt \pi \bar u}
        \begin{pmatrix}
               \sigma_1 \\
              \sigma_2 
        \end{pmatrix}e_{2k^*} .
\end{equation*}

In the next step the crucial assumption \(\hat W(2k^*) = 0\) implies that the matrix \(A_{2k^*}\) is diagonal and hence we find 
\begin{equation*}
    (2\lambda_+(k^*) - A_{2k^*})^{-1} = 
    \begin{pmatrix}
        \frac{1}{2\lambda_+(k^*) + |2k^*|^2 \sigma_1} & 0 \\
        0 & \frac{1}{2\lambda_+(k^*) + |2k^*|^2 \sigma_2}
    \end{pmatrix}
    , \quad 
    A_{2k^*}^{-1} = -
    \begin{pmatrix}
        \frac{1}{|2k^*|^2 \sigma_1} & 0 \\
        0 & \frac{1}{|2k^*|^2 \sigma_2}
    \end{pmatrix}. 
\end{equation*}
This implies 
\begin{equation*}
     (2\lambda_+(k^*) - A_{2k^*})^{-1}     D_{uu}^2 G(\boldsymbol 0, \alpha_{k^*} )[q_{k^*}, q_{k^*}]
        = \frac{1}{\sqrt \pi \bar u }
        \begin{pmatrix}
           \frac{\sigma_1+i\omega_0}{i \omega_0 + 2\sigma_1 } \\
          \chi_{k^*}^2 \frac{\sigma_2+i \omega_0}{i \omega_0 + 2 \sigma_2}
        \end{pmatrix} e_{2k^*}
\end{equation*}

and 
\begin{equation*}
      A_{2k^*}^{-1}     D_{uu}^2 G(\boldsymbol 0, \alpha_{k^*} )[q_{k^*},\bar q_{k^*}]
        =- \frac{1}{2 \sqrt \pi \bar u }
        \begin{pmatrix}
           1 \\ 1 
        \end{pmatrix} e_{2k^*}.
\end{equation*}

Utilizing the above formulas and for the moment we ignore the \(e_{3k^*}\) modes, since they will vanish anyway by pairing the vector with the adjoint, we find 
\begin{align*}
    &  D_{uu}^2 G(\boldsymbol 0, \alpha_{k^*} )[ \bar q_{k^*},  
        (2\lambda_+(k^*) - A_{2k^*})^{-1}     D_{uu}^2 G(\boldsymbol 0, \alpha_{k^*} )[q_{k^*},  q_{k^*} ]]\\
        &\quad =\frac{|k^*|^2}{2 \pi \bar u ^2}
            \begin{pmatrix}
           \frac{(\sigma_1+i\omega_0 )(i \omega_0-\sigma_1)}{i \omega_0 + 2\sigma_1 }    \\
        \frac{\chi_{k^*}(i\omega_0- \sigma_2)(i\omega_0+ \sigma_2)}{(i\omega_0+2\sigma_2) }  
        \end{pmatrix}e_{k^*} .
\end{align*}

Pairing it with the adjoint eigenvector \(q^*\) we obtain
\begin{align*}
      & \Big\langle  D_{uu}^2 G(\boldsymbol 0, \alpha_{k^*} )[ \bar q_{k^*},  
        (2\lambda_+(k^*) - A_{2k^*})^{-1}     D_{uu}^2 G(\boldsymbol 0, \alpha_{k^*} )[q_{k^*}, q_{k^*}]]  , q_{k^*}^* \Big \rangle \\
&\quad = \frac{|k^*|^2}{2(1- \chi_{k^*}^2) \pi \bar u^2 } \Big(   \frac{\sigma_1+i\omega_0}{i \omega_0 + 2\sigma_1 }  (i \omega_0-\sigma_1) - \frac{\chi_{k^*}^2 (i\omega_0- \sigma_2)(i\omega_0+ \sigma_2)}{(i\omega_0+2\sigma_2) } \Big) 
\end{align*}

For the second term we again ignore the \(3k^* \) modes, which yields
\begin{equation*}
            D_{uu}^2 G(\boldsymbol 0, \alpha_{k^*} )[ q_{k^*},  A_{2k^*}^{-1}     D_{uu}^2 G(\boldsymbol 0, \alpha_{k^*} )[q_{k^*} , \bar q_{k^*}]] = \frac{|k^*|^2}{4 \pi \bar u^2 }
            \begin{pmatrix}
          \sigma_1 +i \omega_0  \\
          (i\omega_0+\sigma_2) \chi_{k^*}
        \end{pmatrix} e_{k^*}.
\end{equation*}
Pairing it again with the adjoint eigenvector \(q_{k^*}^*\) we find
\begin{align*}
           2 &\Big\langle    D_{uu}^2 G(\boldsymbol 0, \alpha_{k^*} )[ q_{k^*},  A_{2k^*}^{-1}     D_{uu}^2 G(\boldsymbol 0, \alpha_{k^*} )[q_{k^*}, \bar q_{k^*}]] , q_{k^*}^* \Big\rangle \\
        &\quad = \frac{|k^*|^2}{2 (1-\chi_{k^*}^2) \pi \bar u^2 } \Big( \sigma_1 + i \omega_0 -  \chi_{k^*}^2   (i\omega_0+\sigma_2) \Big)  . 
\end{align*}

Subtracting the first term from the second we find 
\begin{align*}
    &\frac{3|k^*|^2}{2 (1-\chi_{k^*}^2) \pi \bar u^2 } \Big( 
    \frac{\sigma_1 ( \sigma_1 + i \omega_0) }{i\omega_0 + 2\sigma_1} - 
    \chi_{k^*}^2 \frac{\sigma_2 (i \omega_0 + \sigma_2)}{i\omega_0 + 2 \sigma_2 } \Big) \\
    &\quad =
     \frac{3|k^*|^2}{2  (1-\chi_{k^*}^2) \pi \bar u^2 } \Big( 
    \frac{\sigma_1 ( \sigma_1 + i \omega_0) (2\sigma_1 -i \omega_0)  }{\omega_0^2 + 4\sigma_1^2} - 
    \chi_{k^*}^2 \frac{\sigma_2 (i \omega_0 + \sigma_2) (2\sigma_2 -i \omega_0) }{\omega_0^2 + 4 \sigma_2^2 } \Big) \\
    &\quad = \frac{3|k^*|^2}{2  (1-\chi_{k^*}^2) \pi \bar u^2 } \Big( 
    \frac{\sigma_1 ( 2\sigma_1^2 + \omega_0^2 + i \omega_0 \sigma_1)  }{\omega_0^2 + 4\sigma_1^2} - 
    \chi_{k^*}^2 \frac{\sigma_2 ( 2\sigma_2^2 + \omega_0^2 + i \omega_0 \sigma_2)}{\omega_0^2 + 4 \sigma_2^2 } \Big) . 
\end{align*}

In the next step let us simplify the expression \(1/(1-\chi_{k^*}^2) \). Using the definition of \(\chi_{k^*} \) and~\eqref{eq: square_root_relation} we obtain
\begin{align} \label{eq: chi_quadrat_formula}
\begin{split}
    &1-\chi_{k^*}^2 = (1-\chi_{k^*}) (1+\chi_{k^*}) \\
    &\quad = \frac{\big( \gamma \bar u \hat W(k^*) -(\sigma_2+\alpha_{k^*}  \bar u \hat W_2(k^*))  + i \omega_0  \big) \big( \gamma \bar u \hat W(k^*) +(\sigma_2+\alpha_{k^*}  \bar u \hat W_2(k^*))  - i \omega_0  \big)  }{\gamma^2 \bar u^2 \hat W^2(k^*)} \\
    &\quad = \frac{\big( \gamma^2 \bar u^2 \hat W^2(k^*)- ( \sigma_2+\alpha_{k^*} \bar u \hat W_2(k^*)-i \omega_0)^2  }{\gamma^2 \bar u^2 \hat W^2(k^*)}\\
    &\quad  = \frac{\big( \gamma^2 \bar u^2 \hat W^2(k^*)- ( \sigma_2+\alpha_{k^*} \bar u \hat W_2(k^*))^2 +2 i ( \sigma_2+\alpha_{k^*} \bar u \hat W_2(k^*)) \omega_0 + \omega_0^2  }{\gamma^2 \bar u^2 \hat W^2(k^*)} \\
    &= 2 \frac{\omega_0  (i( \sigma_2+\alpha_{k^*} \bar u \hat W_2(k^*))+ \omega_0) }{\gamma^2 \bar u^2 \hat W^2(k^*)}
    \end{split}
\end{align}
This implies 
\begin{align*}
    \frac{1}{1-\chi_{k^*}^2 } &=  \frac{ \gamma^2 \bar u^2 \hat W^2(k^*) }{ 2 \omega_0  (i( \sigma_2+\alpha_{k^*} \bar u \hat W_2(k^*))+ \omega_0)} \\ 
    &=   \frac{ \gamma^2 \bar u^2 \hat W^2(k^*) (\omega_0 - i( \sigma_2+\alpha_{k^*} \bar u \hat W_2(k^*))) }{ 2 \omega_0  (( \sigma_2+\alpha_{k^*} \bar u \hat W_2(k^*))^2+ \omega_0^2) }\\
    &=   \frac{ \gamma^2 \bar u^2 \hat W^2(k^*) (\omega_0 - i( \sigma_2+\alpha_{k^*} \bar u \hat W_2(k^*))) }{2  \omega_0  \gamma^2 \bar u^2 \hat W^2(k^*)}\\
    &=   \frac{  (\omega_0 - i( \sigma_2+\alpha_{k^*} \bar u \hat W_2(k^*))) }{2 \omega_0  }, 
\end{align*}
where we used again relation~\eqref{eq: square_root_relation}. Additionally, we have 
\begin{equation} \label{eq: aux_realpart_computations_identitiy}
    (\omega_0 - i( \sigma_2+\alpha_{k^*}  \bar u \hat W_2(k^*)))\chi_{k^*}^2 = - (\omega_0 + i ( \sigma_2+\alpha_{k^*}  \bar u \hat W_2(k^*))).
\end{equation}
Indeed, 
\begin{align*}
    &(\omega_0 - i( \sigma_2+\alpha_{k^*} \bar u \hat W_2(k^*)))\chi_{k^*}^2 
    = \frac{(\omega_0 - i( \sigma_2+\alpha_{k^*} \bar u \hat W_2(k^*))) (-i\omega_0 + ( \sigma_2+\alpha_{k^*} \bar u \hat W_2(k^*)))^2}{\gamma^2 \bar u^2 \hat W^2(k^*)} \\
    &= \frac{- i (\omega_0 - i( \sigma_2+\alpha_{k^*} \bar u \hat W_2(k^*)))(\omega_0 + i( \sigma_2+\alpha_{k^*} \bar u \hat W_2(k^*))) (-i\omega_0 + ( \sigma_2+\alpha_{k^*} \bar u \hat W_2(k^*)))}{\gamma^2 \bar u^2 \hat W^2(k^*)} \\
    &= \frac{- i (\omega_0^2 +( \sigma_2+\alpha_{k^*} \bar u \hat W_2(k^*))^2 ) (-i\omega_0 + ( \sigma_2+\alpha_{k^*} \bar u \hat W_2(k^*)))}{\gamma^2 \bar u^2 \hat W^2(k^*)}
\\
&= -i  (-i\omega_0 + ( \sigma_2+\alpha_{k^*} \bar u \hat W_2(k^*))) = -(\omega_0 +i ( \sigma_2+\alpha_{k^*} \bar u \hat W_2(k^*)))
\end{align*}

Plugging first our expression for \(1/(1-\chi_{k^*}^2)\) into our formula and then using~\eqref{eq: aux_realpart_computations_identitiy} we find 
\begin{align*}
    &\frac{3|k^*|^2}{2  (1-\chi_{k^*}^2) \pi \bar u^2 } \Big( 
    \frac{\sigma_1 ( 2\sigma_1^2 + \omega_0^2 + i \omega_0 \sigma_1)  }{\omega_0^2 + 4\sigma_1^2} - 
    \chi_{k^*}^2 \frac{\sigma_2 ( 2\sigma_2^2 + \omega_0^2 + i \omega_0 \sigma_2)}{\omega_0^2 + 4 \sigma_2^2 } \Big) \\
    &\quad =   \frac{3|k^*|^2}{4  \pi \bar u^2 \omega_0 } \Big( 
    \frac{ (\omega_0 - i( \sigma_2+\alpha_{k^*} \bar u \hat W_2(k^*))) \sigma_1  ( 2\sigma_1^2 + \omega_0^2 + i \omega_0 \sigma_1)  }{\omega_0^2 + 4\sigma_1^2} \\
    &\quad \quad - 
     (\omega_0 - i( \sigma_2+\alpha_{k^*} \bar u \hat W_2(k^*))) \chi_{k^*}^2 \frac{\sigma_2 ( 2\sigma_2^2 + \omega_0^2 + i \omega_0 \sigma_2)}{\omega_0^2 + 4 \sigma_2^2 } \Big)\\ 
    &\quad =   \frac{3|k^*|^2}{4  \pi \bar u^2 \omega_0 } \Big( 
    \frac{ (\omega_0 - i( \sigma_2+\alpha_{k^*} \bar u \hat W_2(k^*))) \sigma_1 ( 2\sigma_1^2 + \omega_0^2 + i \omega_0 \sigma_1)  }{\omega_0^2 + 4\sigma_1^2} \\
    &\quad \quad + 
     \frac{ (\omega_0 + i( \sigma_2+\alpha_{k^*} \bar u \hat W_2(k^*)))\sigma_2 ( 2\sigma_2^2 + \omega_0^2 + i \omega_0 \sigma_2)}{\omega_0^2 + 4 \sigma_2^2 } \Big) 
\end{align*}
Taking the real part of the above equation, we can find for the first term in the brackets
\begin{align*}
&\Re\Big(\frac{ (\omega_0 - i( \sigma_2+\alpha_{k^*} \bar u \hat W_2(k^*))) \sigma_1 ( 2\sigma_1^2 + \omega_0^2 + i \omega_0 \sigma_1)  }{\omega_0^2 + 4\sigma_1^2} \Big) \\
&\quad = \frac{\sigma_1}{\omega_0^2 + 4\sigma_1^2} \big(  \omega_0 (2\sigma_1^2+\omega_0^2) + ( \sigma_2+\alpha_{k^*} \bar u \hat W_2(k^*))\omega_0 \sigma_1      \big)  
\end{align*}
and for the imaginary part we derive
\begin{align*}
    &\Im \Big( \frac{ (\omega_0 - i( \sigma_2+\alpha_{k^*} \bar u \hat W_2(k^*))) \sigma_1 ( 2\sigma_1^2 + \omega_0^2 + i \omega_0 \sigma_1)  }{\omega_0^2 + 4\sigma_1^2} \Big) \\
&\quad = \frac{\sigma_1}{\omega_0^2+4\sigma_1^2} \big(
\sigma_1\omega_0^2
-\bigl(\sigma_2+\alpha_{k^*}\bar u\,\widehat W_2(k^*)\bigr)
\bigl(2\sigma_1^2+\omega_0^2\bigr) \big)
.
\end{align*}
Analogously 
\begin{align*}
&\Re\Big(\frac{ (\omega_0 +  i( \sigma_2+\alpha_{k^*} \bar u \hat W_2(k^*))) \sigma_2 ( 2\sigma_2^2 + \omega_0^2 + i \omega_0 \sigma_2)  }{\omega_0^2 + 4\sigma_2^2} \Big) \\
&\quad = \frac{\sigma_2}{\omega_0^2 + 4\sigma_2^2} \big(  \omega_0 (2\sigma_2^2+\omega_0^2) - ( \sigma_2+\alpha_{k^*} \bar u \hat W_2(k^*))\omega_0 \sigma_2      \big)  
\end{align*}
and
\begin{align*}
   &\Im\Big(\frac{ (\omega_0  + i( \sigma_2+\alpha_{k^*} \bar u \hat W_2(k^*))) \sigma_2 ( 2\sigma_2^2 + \omega_0^2 + i \omega_0 \sigma_2)  }{\omega_0^2 + 4\sigma_2^2} \Big) \\
&\quad = 
\frac{\sigma_2}{\omega_0^2+4\sigma_2^2} \big( 
\sigma_2\omega_0^2
+\bigl(\sigma_2+\alpha_{k^*} \bar u\,\widehat W_2(k^*)\bigr)
\bigl(2\sigma_2^2+\omega_0^2\bigr)\big) .
\end{align*}

Hence, we find 
\begin{align*}
    &\Re \bigg(   -  D_{uu}^2 G(\boldsymbol 0, \alpha_{k^*})[ \begin{pmatrix}
            e_{k^*} \\ \bar \chi_{k^*} e_{k^*}
        \end{pmatrix},  
        (2\lambda_+(k^*) - A_{2k^*})^{-1}     D_{uu}^2 G(\boldsymbol 0, \alpha_{k^*})[\begin{pmatrix}
            e_{k^*} \\ \chi_{k^*} e_{k^*}
        \end{pmatrix}, \begin{pmatrix}
            e_{k^*} \\ \chi_{k^*} e_{k^*}
        \end{pmatrix}]]\\      
         &\quad +2 \Big[   D_{uu}^2 G(\boldsymbol 0, \alpha_{k^*})[ \begin{pmatrix}
            e_{k^*} \\ \chi_{k^*} e_{k^*}
        \end{pmatrix},  A_{2k^*}^{-1}     D_{uu}^2 G(\boldsymbol 0, \alpha_{k^*})[\begin{pmatrix}
            e_{k^*} \\ \chi_{k^*} e_{k^*}
        \end{pmatrix}, \begin{pmatrix}
            e_{k^*} \\ \bar \chi_{k^*} e_{k^*}
        \end{pmatrix}]] , q^* \Big] \bigg)  \\
        &\quad = \frac{3|k^*|^2}{4  \pi \bar u^2  } \bigg( \frac{\sigma_1 \big( 2\sigma_1^2+\omega_0^2 + \sigma_1( \sigma_2+\alpha_{k^*} \bar u \hat W_2(k^*))     \big)  }{\omega_0^2 + 4\sigma_1^2} + \frac{\sigma_2\big( 2\sigma_2^2+\omega_0^2-\sigma_2  ( \sigma_2+\alpha_{k^*} \bar u \hat W_2(k^*))     \big)  }{\omega_0^2 + 4\sigma_2^2}  \bigg) \\
          &\quad = \frac{3|k^*|^2}{4  \pi \bar u^2  } \bigg( \frac{\sigma_1 \big( (2\sigma_1^2+\omega_0^2) + \sigma_1 S_2    \big)  }{\omega_0^2 + 4\sigma_1^2} + \frac{\sigma_2\big( (2\sigma_2^2+\omega_0^2) -\sigma_2  S_2     \big)  }{\omega_0^2 + 4\sigma_2^2}  \bigg) 
\end{align*}
Let us write the above fraction on a common denominator and recall \(S_2= \sigma_2+\alpha_{k^*} \bar u \hat W_2(k^*)\). 
Expanding the two products,
\begin{align*}
\sigma_1(\omega_0^2+2\sigma_1^2+\sigma_1S_2)(\omega_0^2+4\sigma_2^2)
&=\sigma_1\omega_0^4+\sigma_1\omega_0^2\bigl(4\sigma_2^2+2\sigma_1^2+\sigma_1S_2\bigr)+4\sigma_1\sigma_2^2\bigl(2\sigma_1^2+\sigma_1S_2\bigr),\\
\sigma_2(\omega_0^2+2\sigma_2^2-\sigma_2S_2)(\omega_0^2+4\sigma_1^2)
&=\sigma_2\omega_0^4+\sigma_2\omega_0^2\bigl(4\sigma_1^2+2\sigma_2^2-\sigma_2S_2\bigr)+4\sigma_1^2\sigma_2\bigl(2\sigma_2^2-\sigma_2S_2\bigr).
\end{align*}
Collecting powers of $\omega_0^2$:
 The quadratic term has coefficient $\sigma_1+\sigma_2$.
 For the constant term the last bracket in each line expands to $8\sigma_1^3\sigma_2^2+4\sigma_1^2\sigma_2^2S_2$ and $8\sigma_1^2\sigma_2^3-4\sigma_1^2\sigma_2^2S_2$ respectively; the two $S_2$-terms are equal and opposite and cancel exactly, leaving
\[ 8\sigma_1^3\sigma_2^2+8\sigma_1^2\sigma_2^3=8\sigma_1^2\sigma_2^2(\sigma_1+\sigma_2).
\]
For the linear term, the coefficient in front of $\omega_0^2$ is
\[ \sigma_1\bigl(4\sigma_2^2+2\sigma_1^2+\sigma_1S_2\bigr)+\sigma_2\bigl(4\sigma_1^2+2\sigma_2^2-\sigma_2S_2\bigr)
=2(\sigma_1^3+\sigma_2^3)+4\sigma_1\sigma_2(\sigma_1+\sigma_2)+S_2(\sigma_1^2-\sigma_2^2).
\]
Using $\sigma_1^3+\sigma_2^3=(\sigma_1+\sigma_2)(\sigma_1^2-\sigma_1\sigma_2+\sigma_2^2)$ and $\sigma_1^2-\sigma_2^2=(\sigma_1-\sigma_2)(\sigma_1+\sigma_2)$, this factors as
\[
(\sigma_1+\sigma_2)\Bigl[2(\sigma_1^2-\sigma_1\sigma_2+\sigma_2^2)+4\sigma_1\sigma_2+S_2(\sigma_1-\sigma_2)\Bigr]
=(\sigma_1+\sigma_2)\Bigl[2(\sigma_1^2+\sigma_1\sigma_2+\sigma_2^2)+S_2(\sigma_1-\sigma_2)\Bigr].
\]
Every coefficient therefore carries the common factor $\sigma_1+\sigma_2$, and dividing it out identifies the bracket exactly with $\omega_0^4+c_1\omega_0^2+8\sigma_1^2 \sigma_2^2=P (\omega_0^2) $, where \( P(x) = x^2 + c_1 x + 8 \sigma_1^2 \sigma_2^2\) and 
\(c_1\) given by~\eqref{eq: c1_def}. 

For the imaginary part we find 
\begin{equation*}
    \Im (\Lambda) = \frac{3|k^*|^2}{4 \omega_0 \pi \bar u^2 } \bigg( \frac{\sigma_1 \big(
\sigma_1\omega_0^2
-S_2
\bigl(2\sigma_1^2+\omega_0^2\bigr) \big)}{\omega_0^2+4\sigma_1^2}  + \frac{\sigma_2\big( 
\sigma_2\omega_0^2
+S_2
\bigl(2\sigma_2^2+\omega_0^2\bigr)\big)}{\omega_0^2+4\sigma_2^2}  \bigg).
\end{equation*}
Expanding the two products,
\begin{align*}
\sigma_1\bigl(\sigma_1\omega_0^2-S_2(2\sigma_1^2+\omega_0^2)\bigr)(\omega_0^2+4\sigma_2^2)
&=\sigma_1(\sigma_1-S_2)\omega_0^4+\Bigl[4\sigma_1\sigma_2^2(\sigma_1-S_2)-2\sigma_1^3S_2\Bigr]\omega_0^2-8\sigma_1^3\sigma_2^2S_2,\\
\sigma_2\bigl(\sigma_2\omega_0^2+S_2(2\sigma_2^2+\omega_0^2)\bigr)(\omega_0^2+4\sigma_1^2)
&=\sigma_2(\sigma_2+S_2)\omega_0^4+\Bigl[4\sigma_1^2\sigma_2(\sigma_2+S_2)+2\sigma_2^3S_2\Bigr]\omega_0^2+8\sigma_1^2\sigma_2^3S_2.
\end{align*}
Collecting powers of \(\omega_0^2\):
The quadratic term coefficient is given by 
\[\sigma_1(\sigma_1-S_2)+\sigma_2(\sigma_2+S_2)=\sigma_1^2+\sigma_2^2+S_2(\sigma_2-\sigma_1).
\]
For the constant term coefficient the two boundary contributions are
\[-8\sigma_1^3\sigma_2^2S_2+8\sigma_1^2\sigma_2^3S_2=8\sigma_1^2\sigma_2^2S_2(\sigma_2-\sigma_1).
\]
Finally, the linear term coefficient in front of \(\omega_0^2\) is
\[4\sigma_1\sigma_2^2(\sigma_1-S_2)-2\sigma_1^3S_2+4\sigma_1^2\sigma_2(\sigma_2+S_2)+2\sigma_2^3S_2
=8\sigma_1^2\sigma_2^2+S_2\bigl(4\sigma_1^2\sigma_2-4\sigma_1\sigma_2^2-2\sigma_1^3+2\sigma_2^3\bigr).
\]
Using \(\sigma_2^3-\sigma_1^3=(\sigma_2-\sigma_1)(\sigma_1^2+\sigma_1\sigma_2+\sigma_2^2)\) and \(\sigma_1^2\sigma_2-\sigma_1\sigma_2^2=\sigma_1\sigma_2(\sigma_1-\sigma_2)\), the bracket factors as
\[4\sigma_1\sigma_2(\sigma_1-\sigma_2)+2(\sigma_2-\sigma_1)(\sigma_1^2+\sigma_1\sigma_2+\sigma_2^2)
=2(\sigma_2-\sigma_1)\bigl[\sigma_1^2-\sigma_1\sigma_2+\sigma_2^2\bigr],
\]
so that the linear-term coefficient equals \(8\sigma_1^2\sigma_2^2+2S_2(\sigma_2-\sigma_1)(\sigma_1^2-\sigma_1\sigma_2+\sigma_2^2)\).
Altogether, denoting \(x:=\omega_0^2\), the numerator splits into an \(S_2\)-independent and an \(S_2\)-proportional part
\[N(x)=x\bigl[(\sigma_1^2+\sigma_2^2)x+8\sigma_1^2\sigma_2^2\bigr]+S_2(\sigma_2-\sigma_1)\,Q(x),
\]
where
\begin{equation*}
Q(x):=x^2+2\bigl(\sigma_1^2-\sigma_1\sigma_2+\sigma_2^2\bigr)x+8\sigma_1^2\sigma_2^2.
\end{equation*}
Hence
\begin{equation*}
\Im(\Lambda)=\frac{3|k^*|^2}{4\pi\bar u^2 \omega_0\,(\omega_0^2+4\sigma_1^2)(\omega_0^2+4\sigma_2^2)}
\Bigl[\,\omega_0^2\bigl((\sigma_1^2+\sigma_2^2)\omega_0^2+8\sigma_1^2\sigma_2^2\bigr)+S_2(\sigma_2-\sigma_1)\,Q(\omega_0^2)\Bigr].
\end{equation*}

\end{proof}

A simple consequence of Corollary~\ref{cor: real_part_vanishing_case} is the characterization of sub- vs. supercriticality depending on the polynomial \(P\). To state it, we introduce the constants
\begin{equation}\label{eq: AB_def}
    B_1 := 2\sigma_1^2 + (3+4\sqrt 2)\,\sigma_1\sigma_2 + \sigma_2^2,
    \qquad
    B_2 := \sigma_1^2 + (3+4\sqrt 2)\,\sigma_1\sigma_2 + 2\sigma_2^2 ,
\end{equation}
which are strictly positive and satisfy \(B_2(\sigma_1,\sigma_2)=B_1(\sigma_2,\sigma_1)\) and
\(B_1-B_2=\sigma_1^2-\sigma_2^2\), together with the quantity
\begin{equation}\label{eq: Q_criterion}
    Z(k^*)
    :=
    \bigl(\hat W_1(k^*)+\hat W_2(k^*)\bigr)
    \bigl(B_1\,\hat W_1(k^*)+B_2\,\hat W_2(k^*)\bigr) ,
\end{equation}
and, whenever \(\hat W_1(k^*)\neq0\), the ratio \(r:=\hat W_2(k^*)/\hat W_1(k^*)\), which satisfies
\(r\neq-1\).

\begin{corollary}[Criticality when \(\hat W_1(2k^*)=\hat W_2(2k^*)=\hat W(2k^*)=0\)]
\label{cor: sub_subcritical_vanishing_case}
Let the assumptions of Corollary~\ref{cor: real_part_vanishing_case} hold.
Then the Hopf bifurcation at \(\alpha_{k^*}\) is supercritical if and only
if \(P(\omega_0^2)>0\), and subcritical if and only if \(P(\omega_0^2)<0\).
\end{corollary}
\begin{proof}
By definition the bifurcation is supercritical if and only if \(\Re(\Lambda)>0\). In
\eqref{eq: ReLambda_P_form} the factor \(3|k^*|^2(\sigma_1+\sigma_2)/\bigl(4\pi\bar u^2
(\omega_0^2+4\sigma_1^2)(\omega_0^2+4\sigma_2^2)\bigr)\) is strictly positive, so
\(\Re(\Lambda)\) and \(P(\omega_0^2)\) carry the same sign, which proves the claim.
\end{proof}

\begin{corollary}\label{cor: equivalent_cond}

Let the assumptions of Corollary~\ref{cor: real_part_vanishing_case} hold. The following statements
are equivalent:
\begin{enumerate}[(i)]
    \item\label{item: super_all_gamma}
    the Hopf bifurcation is supercritical for every admissible \(\gamma\);
    \item\label{item: c1_criterion}
    the coefficient \(c_1\) of~\eqref{eq: c1_def} satisfies \(c_1>-4\sqrt2\,\sigma_1\sigma_2\);
    \item\label{item: Z_criterion}
    \(Z(k^*)>0\);
    \item\label{item: componentwise_criterion}
    \(\hat W_2(k^*)>-\tfrac{B_1}{B_2}\hat W_1(k^*)\) if \(\hat W_1(k^*)+\hat W_2(k^*)>0\), and
    \(\hat W_2(k^*)<-\tfrac{B_1}{B_2}\hat W_1(k^*)\) if \(\hat W_1(k^*)+\hat W_2(k^*)<0\).
\end{enumerate}
\end{corollary}
\begin{proof}
Since \(S_2\) does not depend on \(\gamma\), the quantity
\(\omega_0^2=\bar u^2\gamma^2\widehat W(k^*)^2-S_2^2\) ranges over all of \((0,\infty)\) as \(\gamma\)
ranges over the admissible values. By Corollary~\ref{cor: sub_subcritical_vanishing_case}, statement~\ref{item: super_all_gamma} is therefore equivalent to \(P>0\) on \((0,\infty)\). Since \(P(x)=x^2+c_1x+8\sigma_1^2\sigma_2^2\) has positive leading coefficient and \(P(0)>0\), this holds if and only if \(P\) has no root in \((0,\infty)\), that is, if and only if
\(c_1>-4\sqrt2\,\sigma_1\sigma_2\). Hence,~\ref{item: super_all_gamma} is equivalent to~\ref{item: c1_criterion}. 

 Substituting \(S_2=(\sigma_2\hat W_1(k^*)-\sigma_1\hat W_2(k^*))/(\hat W_1(k^*)+\hat W_2(k^*))\)
into~\eqref{eq: c1_def} and collecting terms yields
\begin{equation}\label{eq: c1_Z_identity}
    c_1+4\sqrt2\,\sigma_1\sigma_2
    =
    \frac{B_1\,\hat W_1(k^*)+B_2\,\hat W_2(k^*)}{\hat W_1(k^*)+\hat W_2(k^*)} ,
\end{equation}
whose right-hand side is positive if and only if \(Z(k^*)>0\). Hence~\ref{item: c1_criterion} is equivalent to~\ref{item: Z_criterion}.
Dividing by the positive constant \(B_2\) and distinguishing the sign of \(\hat W_1(k^*)+\hat W_2(k^*)\) gives the equivalence to~\ref{item: componentwise_criterion}. 
\end{proof}

\begin{corollary} \label{cor: equal_diffusion_super}
    Let the assumptions of Corollary~\ref{cor: real_part_vanishing_case} hold. Additionally, let one of the two conditions hold 
    \[\sigma_1 = \sigma_2  \qquad \mathrm{or} \qquad r=1. \]
    Then, for every admissible \(\gamma\) the Hopf bifurcation is supercritical. 
\end{corollary}
\begin{proof}
    Equal diffusivities \(\sigma_1=\sigma_2\) force \(c_1=6\sigma_1^2>0\) independent of \(\hat W_1(k^*),\hat W_2(k^*)\). 
Identical self-interaction kernels, i.e. \(r=1\) implies \(c_1=\tfrac32(\sigma_1+\sigma_2)^2>0\). In either case~\ref{item: c1_criterion}
holds and the bifurcation is supercritical for every admissible \(\gamma\).
\end{proof}

\begin{corollary}\label{cor: subcritical_window}
Let the assumptions of Corollary~\ref{cor: real_part_vanishing_case} hold, and let \(\hat W_1(k^*)\neq0\). Then the Hopf bifurcation is subcritical for some admissible \(\gamma\) if and only if \(r=\hat W_2(k^*)/\hat W_1(k^*)\) satisfies
\begin{equation}\label{eq: r_interval}
    -\tfrac{B_1}{B_2}<r<-1 \ \text{ if } \sigma_1>\sigma_2,
    \qquad
    -1<r<-\tfrac{B_1}{B_2} \ \text{ if } \sigma_1<\sigma_2 .
\end{equation}
In this case \(P\) has two positive roots \(x_\pm=\tfrac12\bigl(-c_1\pm\sqrt{c_1^2-32\sigma_1^2\sigma_2^2}\bigr)\), and the bifurcation is
subcritical exactly for those admissible \(\gamma\) with \(\omega_0^2\in(x_-,x_+)\) and supercritical for all other admissible \(\gamma\) with \(\omega_0^2\notin\{x_-,x_+\}\). In particular this requires
\(\sigma_1\neq\sigma_2\).
\end{corollary}

\begin{proof}
By Corollary~\ref{cor: equivalent_cond}, the bifurcation is subcritical for some admissible \(\gamma\) if and only if \(Z(k^*)\le 0\) and \(P(\omega_0^2) < 0\) for some \(\omega_0^2> 0\). If \(Z(k^*)=0\), then \(c_1=-4\sqrt2\sigma_1\sigma_2\) and \(P(x)=(x-2\sqrt2\sigma_1\sigma_2)^2\ge0\), so no subcritical value exists. If \(Z(k^*)<0\), then \(c_1<-4\sqrt2\sigma_1\sigma_2\), so \(c_1^2>32\sigma_1^2\sigma_2^2\) and \(c_1<0\), and \(P\) has the two positive roots \(x_\pm\), with
\(P<0\) exactly on \((x_-,x_+)\). Since \(\omega_0^2\) attains every value in \((0,\infty)\), the claim follows once we characterize \(Z(k^*)<0\). Writing \(\hat W_2(k^*)=r\hat W_1(k^*)\) gives \(Z(k^*)=\hat W_1(k^*)^2(1+r)(B_1+B_2r)\), which is negative if and only if \(r\) lies strictly between \(-1\) and \(-B_1/B_2\). The order of these two roots is determined by \(B_1-B_2=\sigma_1^2-\sigma_2^2\), which gives~\eqref{eq: r_interval}. If \(\sigma_1=\sigma_2\), then \(B_1=B_2\) and the interval is empty.
\end{proof}

\begin{remark}
    Notice, that since \(\omega_0\to\infty\) as \(\gamma\to\infty\),
strong non-reciprocity always produces a supercritical bifurcation in the regime of this Section, which foreshadows results of Section~\ref{sec: high_gamma}. 
\end{remark}

\begin{corollary}[Sub Supercritical in the case \(|\hat W_i(2k^*)| \) small]\label{cor: sub_subcritical_small_case}
Let the assumptions of Proposition~\ref{prop: hopf_bifurcation} hold, fix the Fourier coefficients of
\(W_1,W_2,W\) at \(k^*\) and an admissible \(\gamma\). Then, there exists an \(\epsilon >0 \) such that if \(|\hat W(2k^*)| \le \epsilon\), \(|\hat W_i(2k^*)| \le \epsilon\) for \(i=1,2\), 
then the Hopf bifurcation at \(\alpha_{k^*}\) 

is supercritical if \(P(\omega_0^2)>0\) and subcritical if \(P(\omega_0^2)<0\).
\end{corollary}
\begin{proof}
Let us denote by \(\boldsymbol{W} =(W_1,W_2,W) \).
    Tracing the computation of Proposition~\ref{prop: sub_supercritical} carefully we notice that in the general non-vanishing case \(\hat {\boldsymbol{W}}(2k^*)  \neq 0 \) we can write the real part of \(\Lambda\) as follows 
    \begin{equation*}
        \Re (\Lambda) = \frac{3|k^*|^2(\sigma_1+\sigma_2)\,P(\omega_0^2)}{4\pi\bar u^2\,(\omega_0^2+4\sigma_1^2)(\omega_0^2+4\sigma_2^2)} + \Re ( f(\hat {\boldsymbol{W}}(2k^*), \hat {\boldsymbol{W}}(k^*)))
    \end{equation*}
    for some continuous complex valued function \(f\) with \(\Re(f(0,x)) = 0\) for all \(x \in \R^3\). Hence, the claim follows by Corollary~\ref{cor: sub_subcritical_vanishing_case} in combination with the continuity of \(f\) by first fixing a \(\hat {\boldsymbol W}(k^*)\) such that \(P(\omega_0^2) > 0\) and then choosing \(\hat W(2k^*)\) small enough.
    The case \(P(\omega_0^2) < 0\) is proven analogously. 
\end{proof}

\subsection{Identical species}\label{sec: matching_case}

In this subsection, we exploit the matching conditions \(\sigma_1=\sigma_2=\sigma\) and \(W_1=W_2=W\) to simplify the general expression for \(\Lambda\) obtained in Proposition~\ref{prop: sub_supercritical}. We derive explicit representations of \(\Re(\Lambda)\) and \(\Im(\Lambda)\) in terms of polynomials in \(\omega_0^2\), valid for every admissible value of \(\gamma\). In particular, Corollary~\ref{cor: lambda_matching} reduces the sign of \(\Re(\Lambda)\) to that of a cubic polynomial \(N_R\), whose coefficients depend only on \(\sigma\) and the ratio \(\iota=\hat W(2k^*)/\hat W(k^*)\). Since the leading coefficient of \(N_R\) is positive for \(\iota\neq 0\), this representation shows that sufficiently strong non-reciprocity leads to a supercritical Hopf bifurcation. This contrasts with the corresponding reciprocal system~\cite{Carrillo2026}, where the criticality is independent of the interaction strength.
We also compute this threshold \(\gamma^*\) explicitly. 

\begin{corollary}[Closed form of \(\Lambda\) in the matching case \(\sigma_1=\sigma_2=\sigma\), \(W_1=W_2=W\)]
\label{cor: lambda_matching}
Let the assumptions of Proposition~\ref{prop: hopf_bifurcation} hold, \(k^*  \in \mathcal T^- \) and assume in addition the
\emph{matching regime}
\[
  \sigma_1=\sigma_2=:\sigma,\qquad W_1=W_2=:W .
\]
Set
\[
  \iota:=\frac{\widehat W(2k^*)}{\widehat W(k^*)}. 
\]
Then \(S_2=0\), \(\chi_{k^*}=i\), \(\alpha_{k^*}=-\dfrac{\sigma}{\bar u\,\widehat W(k^*)}\), and \(\Lambda\) of Proposition~\ref{prop: sub_supercritical} reduces to
\[
  \Lambda
  =\frac{|k^*|^{2}\bigl(\sigma\,N_R(\omega_0^2)+i\,\omega_0\,N_I(\omega_0^2)\bigr)}{D_{b \tilde b}(\omega_0^2)}
\]
where the denominator is the strictly positive polynomial
\[
  D_{b \tilde b}(x)=2\pi\bar u^{2}\,
  \bigl(\iota^{2}x+\sigma^{2}(\iota-1)^{2}\bigr)
  \bigl((2\iota-1)^{2}x+4\sigma^{2}(\iota-1)^{2}\bigr)
  \bigl((2\iota+1)^{2}x+4\sigma^{2}(\iota-1)^{2}\bigr),
\]
and the numerators are the cubic in \(\omega_0^2\)
\[
\begin{aligned}
  N_R(x)&=4\iota^{2}\bigl(12\iota^{4}-4\iota^{3}+\iota^{2}-\iota+1\bigr)x^{3}
        +\sigma^{2}(\iota-1)\bigl(144\iota^{5}-152\iota^{4}+64\iota^{3}-42\iota^{2}+8\iota-3\bigr)x^{2}\\
        &\quad+2\sigma^{4}(\iota-1)^{3}\bigl(72\iota^{3}-56\iota^{2}+22\iota-9\bigr)x
        +24\sigma^{6}(\iota-1)^{5}(2\iota-1),\\
  N_I(x)&=\iota^{2}(2\iota-1)(2\iota+1)(4\iota^{2}+1)x^{3}
        +\iota\sigma^{2}(2\iota+1)\bigl(24\iota^{4}-56\iota^{3}+46\iota^{2}-16\iota+1\bigr)x^{2}\\
        &\quad+\sigma^{4}(\iota-1)^{2}\bigl(48\iota^{4}-80\iota^{3}+32\iota^{2}-2\iota+3\bigr)x
        +4\sigma^{6}(\iota-1)^{4}\bigl(4\iota^{2}-6\iota+3\bigr).
\end{aligned}
\]
In particular, since \(D_{b \tilde b}(x)>0\) for every \(x>0\), the real and imaginary parts are
\[
  \operatorname{Re}\Lambda=\frac{|k^*|^{2}\,\sigma\,N_R(\omega_0^{2})}{D_{b \tilde b}(\omega_0^{2})},
  \qquad
  \operatorname{Im}\Lambda=\frac{|k^*|^{2}\,\omega_0\,N_I(\omega_0^{2})}{D_{b \tilde b}(\omega_0^{2})} .
\]
\end{corollary}

\begin{proof}[Proof of Corollary~\ref{cor: lambda_matching}]
Throughout we use the auxiliary relations of Table~\ref{tab: aux_relations} specialised to the matching
regime, and abbreviate \( \Sigma :=\sigma+\bar u\alpha_{k^*}\widehat W(2k^*)= \sigma(1-\iota)\).
 
\medskip
\emph{Step 1: reduction of the linear data.}
Since \(\sigma_1=\sigma_2=\sigma\) and \(\widehat W_1=\widehat W_2=\widehat W\),
\[
  S_2=\frac{\sigma_2\widehat W_1(k^*)-\sigma_1\widehat W_2(k^*)}{\widehat W_1(k^*)+\widehat W_2(k^*)}
     =\frac{\sigma\widehat W(k^*)-\sigma\widehat W(k^*)}{2\widehat W(k^*)}=0,
\]
and \(\alpha_{k^*}=-\dfrac{\sigma_1+\sigma_2}{\bar u(\widehat W_1(k^*)+\widehat W_2(k^*))}
=-\dfrac{2\sigma}{2\bar u\widehat W(k^*)}=-\dfrac{\sigma}{\bar u\widehat W(k^*)}\).
With \(S_2=0\), the definition \(\omega_0^{2}=\bar u^{2}\gamma^{2}\widehat W(k^*)^{2}-S_2^{2}\) and the
square--root relation~\eqref{eq: square_root_relation} give \(\bar u^{2}\gamma^{2}\widehat W(k^*)^{2}=\omega_0^{2}\),
while \eqref{eq: aux_special_case_lyapunov} yields \(\chi_{k^*}=(S_2-i\omega_0)/ (- \omega_0) = i\), since \(k^* \in \mathcal T^- \), whence
\[
  \chi_{k^*}^{2}=-1,\qquad 1-\chi_{k^*}^{2}=2 .
\]
Moreover \eqref{eq: aux_special_case_lyapunov} gives the two identities we use repeatedly,
\[
  \bar u\,\alpha_{k^*}\widehat W(k^*)=-\sigma,
  \qquad
  \bar u\,\gamma\,\widehat W(k^*)= - \omega_0 ,
\]
and, multiplying by \(\iota=\widehat W(2k^*)/\widehat W(k^*)\),
\[
  \bar u\,\alpha_{k^*}\widehat W(2k^*)=-\sigma\iota,
  \qquad
  \bar u\,\gamma\,\widehat W(2k^*)=-\iota\,\omega_0. 
\]
 
\medskip
\emph{Step 2: the two determinants.}
In the matching regime the \(2k^*\)--block is
\[
  A_{2k^*}=-4|k^*|^{2}
  \begin{pmatrix} \Sigma & \bar u\gamma\widehat W(2k^*)\\[2pt] -\bar u\gamma\widehat W(2k^*) & \Sigma\end{pmatrix}
  =-4|k^*|^{2}\begin{pmatrix} \Sigma &- \iota\omega_0\\[2pt] \iota\omega_0 & \Sigma\end{pmatrix},
\]
using \(\bar u\gamma\widehat W(2k^*)=-\iota\omega_0\). Hence, we have 
\[
  \det A_{2k^*}=16|k^*|^{4}\bigl(\Sigma^{2}+\iota^{2}\omega_0^{2}\bigr)
              =16|k^*|^{4}\bigl(\iota^{2}\omega_0^{2}+\sigma^{2}(1-\iota)^{2}\bigr),
\]
and, since \(2\lambda_+(k^*)I-A_{2k^*}=2|k^*|^{2}\left(\begin{smallmatrix} i\omega_0+2\Sigma & -2\iota\omega_0\\ 2\iota\omega_0 & i\omega_0+2\Sigma\end{smallmatrix}\right)\),
\[
  \det\!\bigl(2\lambda_+(k^*)-A_{2k^*}\bigr)
  =4|k^*|^{4}\Bigl[(i\omega_0+2\Sigma)^{2}+4\iota^{2}\omega_0^{2}\Bigr].
\]

\medskip
\emph{Step 3: reduction of \(b_1,b_2,\tilde b_1,\tilde b_2\).}
Insert Steps 1--2 into the definitions of Proposition~\ref{prop: sub_supercritical}. For \(b_1\) we have 
\(2\lambda_+(k^*)+|2k^*|^{2}\Sigma=2|k^*|^{2}(i\omega_0+2\Sigma)\), \(\chi_{k^*}^{2}=-1\) and
\(\bar u\gamma\widehat W(2k^*)= - \iota\omega_0\), which implies 
\[
  b_1
  =\frac{2|k^*|^{2}\,\dfrac{\sigma+i\omega_0}{\bar u}\,
        \bigl[(2i|k^*|^{2}\omega_0+4|k^*|^{2}\Sigma)-4|k^*|^{2}\iota\omega_0\bigr]}
        {\sqrt\pi\,\det(2\lambda_+(k^*)-A_{2k^*})}
  =\frac{(\sigma+i\omega_0)\bigl(i\omega_0-2\iota\omega_0+2\sigma(1-\iota)\bigr)}
        {\sqrt\pi\,\bar u\,\bigl[(i\omega_0+2\Sigma)^{2}+4\iota^{2}\omega_0^{2}\bigr]} ,
\]
where the \(4|k^*|^{4}\) cancels against \(\det(2\lambda_+(k^*)-A_{2k^*})\). The same substitution gives
\[
\begin{aligned}
  b_2&=\frac{(\sigma+i\omega_0)\bigl(-2\iota\omega_0 -\,i\omega_0-2\sigma(1-\iota)\bigr)}
        {\sqrt\pi\,\bar u\,\bigl[(i\omega_0+2\Sigma)^{2}+4\iota^{2}\omega_0^{2}\bigr]},\\[4pt]
  \tilde b_1&=\frac{\sigma\bigl(-\iota\omega_0-\sigma(1-\iota)\bigr)}
        {2\sqrt\pi\,\bar u\,\bigl[\iota^{2}\omega_0^{2}+\sigma^{2}(\iota-1)^{2}\bigr]},
  \qquad
  \tilde b_2=\frac{-\,\sigma\bigl(-\iota\omega_0+\sigma(1-\iota)\bigr)}
        {2\sqrt\pi\,\bar u\,\bigl[\iota^{2}\omega_0^{2}+\sigma^{2}(\iota-1)^{2}\bigr]},
\end{aligned}
\]
where we used \(W_1= W_2= W\) and \(\bar u\alpha_{k^*}\widehat W(2k^*)=-\sigma\iota\). 
 
\medskip
\emph{Step 4: reduction of the two brackets and assembly.}
Denote by \(\mathrm{FP}\) and \(\mathrm{SP}\) the two bracketed groups of~\eqref{eq: def_Lambda} (the
\(b\)-group and the \(\tilde b\)-group). Substituting
\(\bar u\alpha_{k^*}\widehat W(2k^*)=-\sigma\iota\), \(\bar u\gamma\widehat W(2k^*)=-\iota\omega_0\) and
\(\chi_{k^*}^{2}=-1\), and collecting terms, they group as
\[
  \mathrm{FP}=(b_1-b_2)\,\frac{i\omega_0-\sigma(1-2\iota)}{\bar u}
             +\frac{2\iota\omega_0}{\bar u}\,(b_1+b_2),
  \qquad
  \mathrm{SP}=(\tilde b_1+\tilde b_2)\,\frac{-\sigma(1-2\iota)-i\omega_0}{\bar u}
             -\frac{2\iota\omega_0}{\bar u}\,(\tilde b_1-\tilde b_2).
\]
From Step 3, and writing \(D_b:=(i\omega_0+2\Sigma)^{2}+4\iota^{2}\omega_0^{2}\),
\(D_{\tilde b}:=\iota^{2}\omega_0^{2}+\sigma^{2}(\iota-1)^{2}\) and \(2\Sigma=2\sigma(1-\iota)\), the required
combinations are
\begin{align*}
  b_1-b_2&=\frac{2(\sigma+i\omega_0)(i\omega_0+2\Sigma)}{\sqrt\pi\,\bar u D_b },
  \quad
  b_1+b_2=- \frac{4\iota\omega_0(\sigma+i\omega_0)}{\sqrt\pi\,\bar u D_b},
\\
  \tilde b_1+\tilde b_2&= \frac{\sigma^{2}(\iota-1)}{\sqrt\pi\,\bar u D_{\tilde b}},
  \quad
  \tilde b_1-\tilde b_2=- \frac{\sigma\iota\omega_0}{\sqrt\pi\,\bar u D_{\tilde b}}.
\end{align*}

\emph{Step 5: Reduction of $\mathrm{FP}$}
Substituting $b_1 - b_2$ and $b_1 + b_2$ into $\mathrm{FP}$:
\begin{align*}
    \mathrm{FP} &= \left( \frac{2(\sigma+i\omega_0)(i\omega_0+2\Sigma)}{\sqrt\pi\,\bar u D_b} \right) \frac{i\omega_0-\sigma(1-2\iota)}{\bar u} - \frac{2\iota\omega_0}{\bar u} \left( \frac{4\iota\omega_0(\sigma+i\omega_0)}{\sqrt\pi\,\bar u D_b} \right) \nonumber \\
    &= \frac{2(\sigma+i\omega_0)}{\sqrt\pi\,\bar u^2 D_b} \left[ (i\omega_0+2\Sigma)(i\omega_0-\sigma(1-2\iota)) - 4\iota^2\omega_0^2 \right].
\end{align*}
Using the identity $-4\iota^2\omega_0^2 = (i\omega_0+2\Sigma)^2 - D_b$:
\begin{align*}
    N_{\mathrm{FP}} :&= (i\omega_0+2\Sigma)(i\omega_0-\sigma(1-2\iota)) + (i\omega_0+2\Sigma)^2 - D_b \nonumber \\
    &= (i\omega_0+2\Sigma) \left[ i\omega_0 - \sigma(1-2\iota) + i\omega_0 + 2\sigma(1-\iota) \right] - D_b \\
    &= (i\omega_0+2\Sigma)(2i\omega_0 + \sigma) - D_b.
\end{align*}
Thus:
\begin{equation*}
    \mathrm{FP} = \frac{2(\sigma+i\omega_0)}{\sqrt\pi\,\bar u^2} \left( \frac{(i\omega_0+2\Sigma)(2i\omega_0+\sigma)}{D_b} - 1 \right).
\end{equation*}

\emph{Step 6: Reduction of $\mathrm{SP}$}
Substituting $\tilde b_1 + \tilde b_2$ and $\tilde b_1 - \tilde b_2$ into $\mathrm{SP}$:
\begin{align*}
    \mathrm{SP} &= \left( \frac{\sigma^2(\iota-1)}{\sqrt\pi\,\bar u D_{\tilde b}} \right) \frac{-\sigma(1-2\iota)-i\omega_0}{\bar u} + \frac{2\iota\omega_0}{\bar u} \left( \frac{\sigma\iota\omega_0}{\sqrt\pi\,\bar u D_{\tilde b}} \right)  \\
    &= \frac{\sigma}{\sqrt\pi\,\bar u^2 D_{\tilde b}} \left[ \sigma(\iota-1)(-\sigma(1-2\iota)-i\omega_0) + 2\iota^2\omega_0^2 \right].
\end{align*}
Since $D_{\tilde b} = \iota^2\omega_0^2 + \sigma^2(1-\iota)^2$, we have $2\iota^2\omega_0^2 = 2D_{\tilde b} - 2\sigma^2(1-\iota)^2$:
\begin{align*}
    N_{\mathrm{SP}} &= 2D_{\tilde b} - 2\sigma^2(1-\iota)^2 + \sigma^2(1-\iota)(1-2\iota) + i\omega_0\sigma(1-\iota)  \\
    &= 2D_{\tilde b} + \sigma^2(1-\iota)\left[ -2(1-\iota) + (1-2\iota) \right] + i\omega_0\sigma(1-\iota)  \\
    &= 2D_{\tilde b} - \sigma^2(1-\iota) + i\omega_0\sigma(1-\iota) = 2D_{\tilde b} + \sigma(1-\iota)(i\omega_0 - \sigma).
\end{align*}
Thus:
\begin{equation*}
    \mathrm{SP} = \frac{\sigma}{\sqrt\pi\,\bar u^2} \left( 2 + \frac{\sigma(1-\iota)(i\omega_0 - \sigma)}{D_{\tilde b}} \right).
\end{equation*}

\emph{Step 7: Combining $\mathrm{FP} - 2\,\mathrm{SP}$}
Combining the simplified forms:
\begin{align*}
    \mathrm{FP} - 2\,\mathrm{SP} &= \frac{2}{\sqrt\pi\,\bar u^2} \left( (\sigma+i\omega_0)\left( \frac{(i\omega_0+2\Sigma)(2i\omega_0+\sigma)}{D_b} - 1 \right) - \sigma \left( 2 + \frac{\sigma(1-\iota)(i\omega_0-\sigma)}{D_{\tilde b}} \right) \right)  \\
    &= \frac{2}{\sqrt\pi\,\bar u^2} \left( \frac{(\sigma+i\omega_0)(i\omega_0+2\Sigma)(2i\omega_0+\sigma)}{D_b} - \frac{\sigma^2(1-\iota)(i\omega_0-\sigma)}{D_{\tilde b}} - (3\sigma+i\omega_0) \right).
\end{align*}

\emph{Step 8: Expression for $\Lambda$}
Substituting back into $\Lambda$:
\begin{equation*}
    \Lambda = \frac{|k^*|^2}{2\pi\bar u^2} \left( (3\sigma + i\omega_0) + \frac{\sigma^2(1-\iota)(i\omega_0 - \sigma)}{D_{\tilde b}} - \frac{(\sigma + i\omega_0)(i\omega_0 + 2\Sigma)(2i\omega_0 + \sigma)}{D_b} \right).
\end{equation*}

\emph{Step 9: Decomposition into Real and Imaginary Parts}
Let $\Lambda = \frac{|k^*|^2}{2\pi\bar u^2} (I_1 + I_2 - I_3)$.

\emph{First and Second Terms}
\begin{align*}
    \Re (I_1) &= 3\sigma, \quad \Im (I_1) = \omega_0, \\
    \Re (I_2) &= -\frac{\sigma^3(1-\iota)}{D_{\tilde b}}, \quad \Im (I_2) = \frac{\sigma^2(1-\iota)\omega_0}{D_{\tilde b}}.
\end{align*}

\emph{Third Term $I_3 = \frac{N_3}{D_b}$}
Expanding $N_3 = (\sigma + i\omega_0)(i\omega_0 + 2\Sigma)(2i\omega_0 + \sigma)$:
\begin{align*}
    N_3 &= (\sigma + i\omega_0) \left( (2\Sigma\sigma - 2\omega_0^2) + i\omega_0(4\Sigma + \sigma) \right)  \\
    &= \left( 2\Sigma\sigma^2 - (3\sigma + 4\Sigma)\omega_0^2 \right) + i\omega_0 \left( 6\Sigma\sigma + \sigma^2 - 2\omega_0^2 \right).
\end{align*}
With $\Sigma = \sigma(1-\iota)$, we find 
\begin{align*}
    \Re(N_3) = 2\sigma^3(1-\iota) - \sigma(7 - 4\iota)\omega_0^2, \\
    \Im (N_3) = \omega_0 \left( \sigma^2(7 - 6\iota) - 2\omega_0^2 \right).
\end{align*}

Expanding $D_b = (2\Sigma + i\omega_0)^2 + 4\iota^2\omega_0^2 = \Re (D_b)  + i \Im(D_b) $:
\begin{align*}
    \Re(D_b) = 4\sigma^2(1-\iota)^2 + (4\iota^2 - 1)\omega_0^2, \\
   \Im (D_b) = 4\sigma(1-\iota)\omega_0.
\end{align*}

Multiplying by the complex conjugate $\frac{\Re(D_b)  - i \Im(D_b) }{\Re(D_b)  - i  \Im(D_b)}$:
\begin{align*}
    \text{Re}(I_3) &= \frac{\Re(N_3)\Re(D_b)  +  \Im(N_3)  \Im(D_b)}{\Re(D_b) ^2 +  \Im(D_b)^2}, \\
    \text{Im}(I_3) &= \frac{\Im(N_3) \Re(D_b)  - \Re(N_3)  \Im(D_b)}{\Re(D_b) ^2 +  \Im(D_b)^2}.
\end{align*}

\emph{Final Decomposed Formula}
Collecting the components yields:
\begin{align*}
    \Re(\Lambda) &= \frac{|k^*|^2}{2\pi\bar u^2} \left( 3\sigma - \frac{\sigma^3(1-\iota)}{D_{\tilde b}} - \frac{\Re(N_3) \Re(D_b)  + \Im(N_3) \Im(D_b)}{\Re(D_b) ^2 + \Im(D_b)^2} \right), \\
    \Im(\Lambda) &= \frac{|k^*|^2}{2\pi\bar u^2} \left( \omega_0 + \frac{\sigma^2(1-\iota)\omega_0}{D_{\tilde b}} - \frac{\Im(N_3)  \Re(D_b)  - \Re(N_3)  \Im(D_b) }{\Re(D_b) ^2 + \Im(D_b)^2} \right).
\end{align*}

We want to express $\Lambda$ over a strictly positive common denominator. Let us define this denominator as:
\begin{equation*}
    D_{b \tilde b} := 2\pi\bar u^{2} D_{\tilde b} (\Re(D_b)^{2}+\Im(D_b)^{2}).
\end{equation*}
For the real part, we have:
\begin{equation*}
    \operatorname{Re}(\Lambda) = \frac{|k^*|^2}{D_{b \tilde b}} \Big( \bigl(3\sigma D_{\tilde b} - \sigma^3(1-\iota)\bigr)(\Re(D_b)^2 + \Im(D_b)^2) - (\Re(N_3)\Re(D_b) + \Im(N_3)\Im(D_b) )D_{\tilde b}\Big).
\end{equation*}
By pulling out a common factor of $\sigma$ from the numerator, we can write:
\begin{equation*}
    \operatorname{Re}(\Lambda) = \frac{|k^*|^2 \sigma}{D_{b \tilde b}} N_R(\omega_0^2),
\end{equation*}
where $N_R(\omega_0^2)$ is a polynomial in $x = \omega_0^2$. Expanding the terms with the help of the definitions of \(D_b, D_{\tilde b} \) using and grouping by powers of $x$, we obtain the cubic polynomial 
\begin{align*}
  N_R(x) &= 4\iota^{2}\bigl(12\iota^{4}-4\iota^{3}+\iota^{2}-\iota+1\bigr)x^{3} \nonumber \\
         &\quad +\sigma^{2}(\iota-1)\bigl(144\iota^{5}-152\iota^{4}+64\iota^{3}-42\iota^{2}+8\iota-3\bigr)x^{2} \nonumber \\
         &\quad +2\sigma^{4}(\iota-1)^{3}\bigl(72\iota^{3}-56\iota^{2}+22\iota-9\bigr)x \nonumber \\
         &\quad +24\sigma^{6}(\iota-1)^{5}(2\iota-1).
\end{align*}
Similar computations provide 
\begin{equation*}
    \operatorname{Im}(\Lambda) = \frac{|k^*|^2}{D_{b \tilde b}} \Bigl( \bigl(\omega_0 D_{\tilde b} + \sigma^2(1-\iota)\omega_0\bigr)(\Re(D_b)^2 + \Im(D_b)^2) - \bigl(\Im(N_3)\Re(D_b) - \Re(N_3)\Im(D_b)D_{\tilde b} \bigr) \Bigr).
\end{equation*}
Extracting a common factor of $\omega_0$, we write:
\begin{equation*}
    \operatorname{Im}(\Lambda) = \frac{|k^*|^2 \omega_0}{D_{b \tilde b}} N_I(\omega_0^2).
\end{equation*}
Expanding and grouping by powers of $x = \omega_0^2$ yields the second cubic polynomial 
\begin{align*}
  N_I(x) &= \iota^{2}(2\iota-1)(2\iota+1)(4\iota^{2}+1)x^{3} \nonumber \\
         &\quad +\iota\sigma^{2}(2\iota+1)\bigl(24\iota^{4}-56\iota^{3}+46\iota^{2}-16\iota+1\bigr)x^{2} \nonumber \\
         &\quad +\sigma^{4}(\iota-1)^{2}\bigl(48\iota^{4}-80\iota^{3}+32\iota^{2}-2\iota+3\bigr)x \nonumber \\
         &\quad +4\sigma^{6}(\iota-1)^{4}\bigl(4\iota^{2}-6\iota+3\bigr).
\end{align*}
Utilizing similar steps, we find an
expansion of the Denominator $D_{b \tilde b}$ as 
\begin{align*}
    D_{b \tilde b} &= 2\pi\bar u^2 \Big( \iota^2(4\iota^2 - 1)^2 x^3 + \sigma^2(\iota-1)^2(48\iota^4 + 1) x^2 \nonumber \\
    &\qquad\qquad + 8\sigma^4(\iota-1)^4(6\iota^2 + 1) x + 16\sigma^6(\iota-1)^6 \Big).
\end{align*}

Since \(x > 0\) we always have \(D_{b \tilde b} > 0\).

\end{proof}

\begin{remark}[Opposite sign of \(\widehat W(k^*)\)]
\label{rem: sign_of_What}
The proof of Corollary~\ref{cor: lambda_matching} was carried out for \(k^*\in\mathcal T^-\).
In the matching regime \(\widehat W_1=\widehat W_2=\widehat W\) this is equivalent to
\(\widehat W(k^*)<0\), and hence to \(\alpha_{k^*}=-\sigma/(\bar u\widehat W(k^*))>0\); the
complementary case \(k^*\in\mathcal T^+\) is exactly \(\widehat W(k^*)>0\), i.e.\
\(\alpha_{k^*}<0\). We claim that the closed form of Corollary~\ref{cor: lambda_matching} is
valid in both cases, verbatim.
 
\medskip
Where the sign enters: 
Put
\[
  \epsilon:=\operatorname{sgn}\widehat W(k^*)\in\{-1,+1\} .
\]
Since \(\alpha_{k^*}\) is itself proportional to \(1/\widehat W(k^*)\), the self-interaction data are \(\epsilon\)-independent,
\[
  \bar u\,\alpha_{k^*}\widehat W(k^*)=-\sigma,
  \qquad
  \bar u\,\alpha_{k^*}\widehat W(2k^*)=-\sigma\iota,
  \qquad
  \Sigma=\sigma(1-\iota),
\]
whereas the cross-interaction data and the eigenvector carry \(\epsilon\):
\[
  \bar u\,\gamma\,\widehat W(k^*)=\epsilon\,\omega_0,
  \qquad
  \bar u\,\gamma\,\widehat W(2k^*)=\epsilon\,\iota\,\omega_0,
  \qquad
  \chi_{k^*}=\frac{S_2-i\omega_0}{\bar u\gamma\widehat W(k^*)}=-i\epsilon .
\]
In particular \(\chi_{k^*}^{2}=-1\) and \(1-\chi_{k^*}^{2}=2\) for both signs, so Step~1 of the proof holds with \(\chi_{k^*}=-i\epsilon\). The two determinants of Step~2 are even in
\(\epsilon\),
\[
  \det A_{2k^*}=16|k^*|^{4}\bigl(\iota^{2}\omega_0^{2}+\sigma^{2}(1-\iota)^{2}\bigr),
  \qquad
  \det\bigl(2\lambda_+(k^*)-A_{2k^*}\bigr)=4|k^*|^{4}D_b,
\]
because the off-diagonal entry \(\bar u\gamma\widehat W(2k^*)=\epsilon\iota\omega_0\) enters both only through its square. In Steps~3 and~4 the sign \(\epsilon\) enters linearly through \(\bar u\gamma\widehat W(2k^*)\). Carrying \(\epsilon\) through Step~3 shows that \(b_1-b_2\) and \(\tilde b_1+\tilde b_2\) are independent of \(\epsilon\), whereas \(b_1+b_2\) and \(\tilde b_1-\tilde b_2\) are odd in \(\epsilon\). Moreover, the grouping of Step~4 becomes
\[
\mathrm{FP}=(b_1-b_2)\,\frac{i\omega_0-\sigma(1-2\iota)}{\bar u}-\epsilon\,\frac{2\iota\omega_0}{\bar u}\,(b_1+b_2),
\quad
\mathrm{SP}=(\tilde b_1+\tilde b_2)\,\frac{-\sigma(1-2\iota)-i\omega_0}{\bar u}+\epsilon\,\frac{2\iota\omega_0}{\bar u}\,(\tilde b_1-\tilde b_2).
\]
Hence \(\epsilon\) only appears in the products of \(\epsilon\) with odd combinations, so \(\mathrm{FP}\) and \(\mathrm{SP}\) are independent of \(\epsilon\). Steps~5--9 therefore apply verbatim, and the polynomials \(N_R\), \(N_I\) and \(D_{b\tilde b}\) are the same for \(k^*\in\mathcal T^+\) and \(k^*\in\mathcal T^-\).

\end{remark}

Analyzing the roots of the polynomial $N_R$ yields the following corollaries. We omit their proofs, as they rely only on standard arguments for cubic polynomials.

\begin{corollary}
    Let the assumptions of Corollary~\ref{cor: lambda_matching} hold and let \(\iota<\tfrac12\). Then, \(N_R > 0 \) on \((0, \infty)\) and, consequently, the bifurcation is supercritical. 
\end{corollary}
\begin{corollary}[Change of criticality induced by non-reciprocity]
\label{cor: nonreciprocal_transition}
Let the assumptions of Corollary~\ref{cor: lambda_matching} hold and let \(\tfrac12<\iota<1\).
Then \(N_R\) has exactly one positive root \(x^*\), and with
\[
  \gamma^*:=\frac{\sqrt{x^*}}{\bar u\,|\widehat W(k^*)|}
\]
the Hopf bifurcation at \(\alpha_{k^*}\) is subcritical for \(0<\gamma<\gamma^*\) and supercritical
for \(\gamma>\gamma^*\).
\end{corollary}

\begin{remark}[Comparison with the reciprocal system]
\label{rem: matching_criticality}

Corollary~\ref{cor: nonreciprocal_transition} reveals a fundamental difference between reciprocal and non-reciprocal interactions. While the strength of the coupling does not affect the criticality of the bifurcation in the corresponding reciprocal system, it can change the Hopf bifurcation from subcritical to supercritical in the non-reciprocal system.

To see this, consider the reciprocal system~\eqref{eq: 2species_problem_reciprocal} with off-diagonal coupling \((\gamma,\gamma)\). In the matching regime, the homogeneous state undergoes a stationary pitchfork bifurcation at
\(\alpha^*=-\sigma/(\bar u\widehat W(k^*))-\gamma\).
As shown in Appendix~\ref{app: reciprocal_caluclations}, its branch direction satisfies
\begin{equation*}
    \operatorname{sign}\alpha''(0)
    =\operatorname{sign}(1-2\iota),
    \qquad \iota<1.
\end{equation*}
Consequently, the pitchfork bifurcation is supercritical for \(\iota<\tfrac12\) and subcritical for \(\tfrac12<\iota<1\), independently of the coupling strength \(\gamma\).

The situation is different for the non-reciprocal system~\eqref{eq: 2species_problem}. In precisely the regime \(\tfrac12<\iota<1\), where the reciprocal bifurcation is always subcritical, Corollary~\ref{cor: nonreciprocal_transition} shows that increasing \(\gamma\) beyond the threshold \(\gamma^*\) renders the Hopf bifurcation supercritical. Thus, non-reciprocity does not merely replace a stationary bifurcation by an oscillatory one. It also introduces a mechanism through which the interaction strength itself changes the criticality of the emerging branch.

This distinction has an important consequence for the dynamics. If the Hopf bifurcation is the first loss of stability, the Principle of Exchange of Stability implies that the supercritical branch consists of locally orbitally asymptotically stable periodic solutions within the class of even solutions, as shown in Section~\ref{sec: exchange}. In this sense, sufficiently strong non-reciprocity stabilizes the emerging oscillations. The contrast between the two systems is illustrated in Figure~\ref{fig:reciprocal_nonreciprocal_regimes}.
\end{remark}

\subsection{The strong non-reciprocity regime}\label{sec: high_gamma}

In this subsection, we investigate the criticality of the Hopf bifurcation in the regime of strong non-reciprocity, that is, as \(\gamma\to\infty\). In contrast to the matching regime considered previously, we allow for different diffusivities and interaction kernels. Our aim is to understand whether the dependence of the normal-form coefficient \(\Lambda\) on \(\gamma\) simplifies in this limit and to derive explicit conditions for supercriticality.

We first show in Corollary~\ref{cor: gamma_limit} that \(\Re(\Lambda)\) converges to a finite limit with an error of order \(O(\gamma^{-1})\). This limit depends only on the diffusivities and the ratios of the second-harmonic to the critical Fourier coefficients of the self- and cross-interaction kernels. In particular, whenever the limit is nonzero, the criticality becomes independent of \(\gamma\) for sufficiently strong non-reciprocity. Moreover, formula~\eqref{eq: gamma_limit_ReLambda} allows us to identify explicit conditions on the interaction kernels under which this limiting value is positive, and hence the Hopf bifurcation is supercritical for all sufficiently large \(\gamma\).

\begin{corollary} \label{cor: gamma_limit}
Let the assumptions of Proposition~\ref{prop: hopf_bifurcation} hold and \(\hat W(2k^*) \neq 0\) such that \(2 \hat W(2k^*) \neq \pm \hat W(k^*)\). Then, as the non-reciprocity strength
$\gamma\to\infty$,
\begin{equation}\label{eq: gamma_limit_ReLambda}
  \operatorname{Re}\Lambda
   =\frac{|k^*|^2}{\pi\bar u^{2}}\,
     \frac{ (4\iota^2+1)(\sigma_1+\sigma_2) \big( 1- \tfrac{\hat W_1(2k^*) + \hat W_2(2k^*)}{\hat W_1(k^*)+\hat W_2(k^*)}\big)  + 3\iota^2(4\iota^2-1)(\sigma_1+\sigma_2)}
          {(4\iota^2-1)^2}
   \;+\;O\!\left(\gamma^{-1}\right),
\end{equation}
with 
\(
  \iota:=\frac{\hat W(2k^*)}{\hat W(k^*)} , 
 \)
and consequently, 
$$\alpha''(0)=-\dfrac{2}{\bar u\,(\hat W_1(k^*)+\hat W_2(k^*))\,|k^*|^2}\operatorname{Re}\Lambda. $$ 
In particular the limit is finite and depends solely through \(\iota\), \(\sigma_1, \sigma_2\) and the ratio of \((\hat W_1(2k^*)+ \hat W_2(2k^*) )/(\hat W_1(k^*) + \hat W_2(k^*)) \). 

\end{corollary}

\begin{proof}[Proof of Corollary~\ref{cor: gamma_limit}]
For the statement we keep the abbreviations
\[
  \Sigma_i:=\sigma_i+\bar u\,\alpha_{k^*}\hat W_i(2k^*)\ (i=1,2),\qquad
  \iota:=\frac{\hat W(2k^*)}{\hat W(k^*)} ,
\]
At the end we want to represent everything in term of \(\omega_0\). We have the following relations, where the \(\pm\) corresponds to the sign of \(\hat W(k^*)\): 
\begin{center}
\renewcommand{\arraystretch}{1.5}
\begin{tabular}{lll}
\hline
quantity & exact value & leading behaviour as $\omega_0\to\infty$\\
\hline
$2\lambda_+(k^*)$ & $2i|k^*|^2\omega_0$ & $2i|k^*|^2\omega_0$ \\
$\bar u  \gamma \hat W(k^*)$ & $\pm\sqrt{\omega_0^2+S_2^2}$ & $ \pm \omega_0 + O(\omega^{-1}_0) $\\
$ \bar u  \gamma \hat W(2k^*)$ & $\pm \iota\sqrt{\omega_0^2+S_2^2}$ & $ \pm \iota \omega_0 + O(\omega_0^{-1})$\\
$\chi_{k^*}$ & $\dfrac{S_2-i\omega_0}{\pm \sqrt{\omega_0^2+S_2^2}}$ & $ \mp i+ O(\omega^{-1}_0)$\\
$\chi_{k^*}^2 $ &$\dfrac{(S_2-i\omega_0)^2}{\omega_0^2+S_2^2}$  & $ -1- \tfrac{2i S_2}{\omega_0} + O(\omega_0^{-2} ) $ \\
\hline
\end{tabular}
\end{center}

Without loss of generality, we will assume \(\hat W(k^*) > 0\). This does not change the result, since the formula we are going to derive depends only \(\iota\) in a quadratic way. However, the advantage is we can choose just one branch in the above table. In order to handle the case  \(\hat W(k^*) <0 \), we just use the symmetry of our system~\eqref{eq: 2species_problem}. Replacing $W$ by $-W$ leaves $\iota$ unchanged and is equivalent to exchanging the two species, under which~\eqref{eq: gamma_limit_ReLambda} is invariant. 
Our starting point is Proposition~\ref{prop: sub_supercritical}. 
The proof is based on multiple steps. 
We start with rewriting the determinants. 

\emph{Step 1: Formula for determinants.}
We have 
{\small
\begin{align*}
\det(2\lambda_+(k^*) - A_{2k^*}) &= (2\lambda_+(k^*) + 4|k^*|^2 (\sigma_1+\alpha_{k^*} \bar u \hat W_1(2k^*)) ) (2\lambda_+(k^*) + 4|k^*|^2 (\sigma_2+\alpha_{k^*} \bar u \hat W_2(2k^*)) ) \\
&\quad + 16 |k^*|^4 \bar u^2 \gamma^2 \hat W^2(2k^*)
\end{align*}
}
The terms
\begin{equation}
     4|k^*|^2 (\sigma_1+\alpha_{k^*} \bar u \hat W_1(2k^*)), 4|k^*|^2 (\sigma_2+\alpha_{k^*} \bar u \hat W_2(2k^*))
\end{equation}
are real and independent of \(\omega_0\), i.e \(\gamma\).
Additionally, 
\begin{equation*}
    \gamma^2 = \frac{\omega_0^2}{\bar u^2 \hat W^2(k^*)} + O(1) 
\end{equation*}
Hence we find 
\begin{align*}
\det(2\lambda_+(k^*) - A_{2k^*}) 
    &= (2\lambda_+(k^*) + |2k^*|^2\Sigma_1) (2\lambda_+(k^*) + |2k^*|^2\Sigma_2)+ 16 |k^*|^4 \bar u^2 \gamma^2 \hat W^2(2k^*)  \\
    &= -4|k^*|^4  \omega_0^2+16 |k^*|^4 \omega_0^2 \iota ^2+ 8 |k^*|^4 i\omega_0(\Sigma_1+\Sigma_2) +O(1) \\
    &= \frac{4|k^*|^4 \omega_0^2 }{\hat W^2(k^*)} ( 4 \hat W^2(2k^*) - \hat W^2(k^*)) +  8|k^*|^4i\omega_0(\Sigma_1+\Sigma_2) +O(1)
\end{align*}

The computation for \(\det(A_{2k^*})\) follows by setting \(\lambda_+(k^*) =0\) in the above computations. Hence, we obtain 
\begin{equation*}
  \det(A_{2k^*}) =   16 |k^*|^4 \omega_0^2 \iota^2 +O(1) 
\end{equation*}
Here, we require the fact that \(\hat W(2k^*) \neq 0\). Otherwise, the determinant would be of order \(O(1)\) and we would be in a similar situation as in Corollary~\ref{cor: real_part_vanishing_case}.  

\emph{Step 2: Limits for the coefficient \(b_1\).}
We apply~\eqref{eq: chi_quadrat_formula} to find 
\begin{align*}
        b_1&=  \frac{2|k^*|^2 \bigg( \Big( 2 \lambda_+(k^*) + |2k^*|^2 \Sigma_2 \Big) \frac{\sigma_1+i\omega_0}{\bar u} - |2k^*|^2 \gamma \chi_{k^*}^2 \widehat W(2k^*) (\sigma_2+i\omega_0) \bigg) }{\sqrt \pi \det(2\lambda_+(k^*)- A_{2k^*})} \\
        &=  \frac{4|k^*|^4 \bigg( \Big( i \omega_0  + 2 \Sigma_2 \Big) \frac{\sigma_1+i\omega_0}{\bar u} -2 \gamma  \widehat W(2k^*)  \big( -1- \tfrac{2i S_2}{\omega_0} + O(\omega_0^{-2} ) \big) (\sigma_2+i\omega_0) \bigg) }{\sqrt  \pi \big(   \frac{4|k^*|^4 \omega_0^2 }{\hat W^2(k^*)} ( 4 \hat W^2(2k^*) - \hat W^2(k^*)) + 8|k^*|^4 i\omega_0(\Sigma_1+\Sigma_2)+O(1)\big)  }
\end{align*}
Multiplying out the numerator, we find 
\begin{align*}
    &4|k^*|^4 \Big(  \frac{-\omega_0^2+i\sigma_1\omega_0}{\bar u} +  \frac{2i \Sigma_2\omega_0}{\bar u} + O(1) + 2 \big( \frac{\iota }{\bar u } (\omega_0+ 2iS_2 + O(\omega_0^{-1})    \big)  (\sigma_2 +i \omega_0)  \Big) \\
    &\quad = 4|k^*|^4 \Big(  \frac{- \omega_0^2}{ \bar u} + 2 \frac{i \iota  \omega_0^2}{\bar u}+ \Big(  \frac{i (\sigma_1+2\Sigma_2)}{\bar u} -  \frac{2 \iota (-\sigma_2+2S_2)}{ \bar u}\Big) \omega_0 + O(1) \Big)  .
\end{align*}
Recalling that for \(a_1,a_2,d_1,d_2\) we have 
\begin{equation*}
\frac{a_2 \omega_0^2 + a_1 \omega_0}{d_2 \omega_0^2 + d_1 \omega_0} = \frac{a_2}{d_2}\frac{1+\frac{a_1}{a_2\omega_0}}{1+\frac{d_1}{d_2\omega_0}} = \frac{a_2}{d_2} + \frac{1}{\omega_0}\left(\frac{a_1}{d_2} - \frac{a_2 d_1}{d_2^2}\right) + O(\omega_0^{-2}) = \frac{a_2}{d_2} + \frac{a_1 d_2 - a_2 d_1}{d_2^2\,\omega_0} + O(\omega_0^{-2})
\end{equation*}
implies 
\begin{equation*}
    b_1 =  \frac{2 i \iota  -1}{\sqrt{\pi} \bar u \big( 4 \iota^2 -1\big) } + \beta_1 \omega_0^{-1}  + O(\omega_0^{-2})
\end{equation*}
with 
\begin{align*}
    \beta_1 &= \frac{\big( i (\sigma_1+2\Sigma_2)- 2 \iota (-\sigma_2+2S_2)\big)  \big( 4 \iota^2 -1\big) -\big( 2 i \iota  -1\big) 2i(\Sigma_1+\Sigma_2) }{\sqrt{\pi} \bar u \big( 4 \tfrac{\hat W^2(2k^*)}{\hat W^2(k^*)} -1\big)^2} \\
    &=\frac{2\iota \left(\left(4\iota^2-1\right)(\sigma_2-2S_2)+2(\Sigma_1+\Sigma_2)\right) + i\left(\left(4\iota^2 -1\right)(\sigma_1+2\Sigma_2)+2(\Sigma_1+\Sigma_2)\right)}{\sqrt{\pi} \bar u \big( 4\iota^2 -1\big)^2} . 
\end{align*}

\emph{Step 3: Limits for the coefficient \(b_2\).}
Similarly to the previous step, we rewrite the \(b_2\) coefficient as 
\begin{align*}
    b_2 = \frac{4|k^*|^4 \bigg(2 \bar u\gamma \widehat W(2k^*) \frac{\sigma_1+i\omega_0}{\bar u} + \Big( i \omega_0 + 2 \Sigma_1 \Big)  \chi_{k^*}^2 \frac{\sigma_2+i\omega_0}{\bar u} \bigg) }{\sqrt \pi \det(2\lambda_+(k^*)- A_{2k^*})}, \\
\end{align*}
Again let us tackle the numerator first. We concentrate on the expression in the bracket. We find 
\begin{align*}
    &(2 \iota \omega_0 +O(\omega_0^{-1}))\frac{\sigma_1+i\omega_0}{\bar u} + \Big( i \omega_0 + 2 \Sigma_1 \Big)  \big( -1- \tfrac{2i S_2}{\omega_0} + O(\omega_0^{-2} ) \big) \frac{\sigma_2+i\omega_0}{\bar u} \\
    &\quad =  \frac{2 i \iota }{ \bar u  }\omega_0^2 + \frac{\omega_0^2}{\bar u}  + \bigg(  \frac{2 \sigma_1 \iota }{\bar u } 
    - \frac{i(\sigma_2+2\Sigma_1 - 2 S_2) }{\bar u} \bigg) \omega_0 
    +O(1)
    , 
\end{align*}
which implies by the same steps as in the computations of \(b_1\) that 
\begin{equation*}
    b_2 =  \frac{2 i \iota   +1}{\sqrt{\pi} \bar u \big(4 \iota^2 -1\big) } + \beta_2 \omega_0^{-1}+ O(\omega_0^{-2}) 
\end{equation*}
with 
\begin{align*}
    \beta_2 &= \frac{\big( 2\sigma_1 \iota  -i(\sigma_2+2\Sigma_1 - 2 S_2)\big)  \big( 4 \iota^2 -1\big) - \big( 2 i \iota   +1\big) 2i (\Sigma_1+\Sigma_2)}{\sqrt{\pi} \bar u \big( 4\iota^2 -1\big)^2}\\
    &= \frac{2\iota 
\left( \sigma_1\left(
4 \iota^2 -1 \right) +2(\Sigma_1+\Sigma_2) \right)
- i\left( \left( 4\iota^2 -1 \right) (\sigma_2+2\Sigma_1-2S_2) +2(\Sigma_1+\Sigma_2) \right)}{\sqrt{\pi} \bar u \big( 4 \iota^2 -1\big)^2}
\end{align*}

\emph{Step 4: Limits for the coefficient \(\tilde b_1\).}
Using our table we immediately obtain 
\begin{align*}
       \tilde b_1& =  - \frac{8|k^*|^4 }{\sqrt \pi \det(A_{2k^*})}\bigg(  \frac{\Sigma_2 \sigma_1 }{\bar u} -  \bar u\gamma \widehat W(2k^*)  \frac{\sigma_2}{\bar u}\bigg) \\
       &=  \frac{1}{2 \sqrt \pi  \omega_0^2 \iota^2 +O(1)  } \Big( \frac{\iota \sigma_2}{\bar u} \omega_0  + O(\omega_0^{-1})\Big)   \\
       &=  \frac{\sigma_2}{2 \sqrt \pi \bar u\iota  } \omega_0^{-1} + O(\omega_0^{-2}) 
\end{align*}

\emph{Step 5: Limits for the coefficient \(\tilde b_2\).}
Performing similar calculations to \(\tilde b_1\) we obtain

\begin{equation*}
    \tilde b_2 = - \frac{\sigma_1}{2 \sqrt \pi \bar u \iota }   \omega_0^{-1} + O(\omega_0^{-2}) . 
\end{equation*}

\emph{Step 6: \(\Lambda\) expression.}
We notice that \(\Lambda\) is of the form 
\begin{equation*}
    \Lambda= 
-\frac{  |k^*|^2}{2 \sqrt \pi(1-\chi_{k^*}^2)}  (\mathrm{\mathrm{FP}-2 \cdot \mathrm{\mathrm{SP}}} ) . 
\end{equation*}

\emph{Step 7:  \(\mathrm{\mathrm{FP}}\) expression.}
Using our expressions for \(b_1\) and \(b_2\) we find 
\begin{align*}
    \mathrm{\mathrm{FP}}
    &=b_1\big( \frac{-\sigma_1 }{\bar u }- 2\alpha_{k^* }\hat W_1(2k^*) -2 \gamma \hat W(2k^*) ) \\
    &\quad - b_2 \big( \frac{-\sigma_2 }{\bar u }- 2\alpha_{k^*} \hat W_2(2k^*) +2 \gamma \hat W(2k^*) +  \Big)+ \frac{i \omega_0}{\bar u} ( b_1 - b_2)  \\
     &= b_1\big(\frac{\sigma_1-2\Sigma_1}{\bar u } -2 \gamma \hat W(2k^*) \big) - b_2 \big(\frac{\sigma_2- 2\Sigma_2}{\bar u }+2 \gamma \hat W(2k^*)  \big)+ \frac{i \omega_0}{\bar u} ( b_1 - b_2) \\
    &=:(I)-(II) + (III)
\end{align*}
Splitting now \(b_1,b_2\) into the leading order term, the \(\omega_0^{-1}\) term we find  
\begin{align*}
    (I) &= \big(- \frac{ 1}{\sqrt{\pi} \bar u \big( 4 \iota^2-1\big) }+ \omega_0^{-1}  \Re(\beta_1) \big) \big(\frac{\sigma_1-2\Sigma_1}{\bar u } -2 \gamma \hat W(2k^*) \big) 
    \\
    &\quad + i\big(\frac{\sigma_1-2\Sigma_1}{\bar u } -2 \gamma \hat W(2k^*) \big) \big( \frac{2\iota  }{\sqrt{\pi} \bar u \big( 4 \iota^2 -1\big) } +  \omega_0^{-1}  \Im(\beta_1)\big)  + O(\omega_0^{-2}) \\
&= \big(- \frac{ 1}{\sqrt{\pi} \bar u \big( 4 \iota^2 -1\big) }+  \omega_0^{-1}  \Re(\beta_1) \big) \big(\frac{\sigma_1-2\Sigma_1}{\bar u }  -2\frac{\iota }{ \bar u} \omega_0 \big) 
    \\
    &\quad + i\big(\frac{\sigma_1-2\Sigma_1}{\bar u }-2 \frac{\iota }{ \bar u} \omega_0 \ \big) \big( \frac{2\iota  }{\sqrt{\pi} \bar u \big( 4 \iota^2 -1\big) } + \omega_0^{-1}  \Im(\beta_1)\big)  + O(\omega_0^{-1})   
\end{align*}

\begin{align*}
(II) &= \big( \frac{ 1}{\sqrt{\pi} \bar u \big( 4 \iota^2 -1\big) }+ \omega_0^{-1}   \Re(\beta_2) \big) \big(\frac{\sigma_2-2\Sigma_2}{\bar u }   +2 \gamma \hat W(2k^*) \big) 
    \\
    &\quad + i\big(\frac{\sigma_2-2\Sigma_2}{\bar u } +2 \gamma \hat W(2k^*) \big) \big( \frac{2\iota  }{\sqrt{\pi} \bar u \big( 4 \iota^2 -1\big) } + \omega_0^{-1}  \Im(\beta_2)\big)  + O(\omega_0^{-2}) \\
    &= \big( \frac{ 1}{\sqrt{\pi} \bar u \big( 4 \iota^2 -1\big) }+ \omega_0^{-1}   \Re(\beta_2) \big) \big(\frac{\sigma_2-2\Sigma_2}{\bar u }  +2 \frac{\iota }{\bar u} \omega_0 \ \big) 
    \\
    &\quad + i\big(\frac{\sigma_2-2\Sigma_2}{\bar u }  + 2\frac{\iota }{ \bar u} \omega_0 \big) \big( \frac{2\iota   }{\sqrt{\pi} \bar u \big( 4 \iota^2 -1\big) } + \omega_0^{-1}  \Im(\beta_2)\big)  + O(\omega_0^{-2}) 
\end{align*}

\begin{align*}
  (III) =  \frac{ 1}{\bar u } \Big( \Im(\beta_2)-\Im(\beta_1) +  i\big( \Re(\beta_1)-\Re(\beta_2) -  \frac{2\omega_0 }{\sqrt{\pi} \bar u \big( 4 \iota^2 -1\big) } \big)   \Big)    +O(\omega_0^{-1}) 
\end{align*}

This implies 
\begin{align*}
\mathrm{\mathrm{FP}} 
&= \big( \frac{ -(\sigma_1+\sigma_2) + 2(\Sigma_1+\Sigma_2)}{\sqrt{\pi} \bar u^2 \big( 4 \iota^2 -1\big) } \big)  - 2(\Re(\beta_1) + \Re(\beta_2) ) \frac{\iota }{\bar u}  \\
&\quad -i  \frac{8\iota^2}{\sqrt{\pi} \bar u^2 \big( 4 \iota^2 -1\big) }  \omega_0 - 2i\frac{\iota }{ \bar u}  (\Im(\beta_1)+\Im(\beta_2))\\
&\quad +i \big(\sigma_1-\sigma_2 + 2(\Sigma_2-\Sigma_1) \big) \frac{2\iota }{\sqrt{\pi} \bar u^2 \big( 4\iota^2 -1\big) } + (III) +O(\omega^{-1}_0) 
\end{align*}

\emph{Step 8: \( \mathrm{\mathrm{SP}}\) expression.}
For the \(\mathrm{\mathrm{SP}}\) we obtain 
\begin{align*}
    \mathrm{\mathrm{SP}} &=   \tilde  b_1 \Big(  -\frac{\sigma_1 + i \omega_0}{\bar u}    -2  \alpha_{k^*} \hat W_1(2k^*)-  2 \chi_{k^*}^2  \gamma \hat W(2k^*) \Big)  \\
&\quad + \tilde b_2\Big(  \chi_{k^*}^2  \frac{\sigma_2+ i \omega_0}{\bar u} + 2 \chi_{k^*}^2   \alpha_{k^*} \hat W_2(2k^*) -2  \gamma \hat W(2k^*)\Big) \\
&= \tilde  b_1 \Big(  -\frac{i \omega_0}{\bar u}    -  2\big(-1- \tfrac{2i S_2}{\omega_0} \big)  \frac{\iota }{\bar u} \omega_0  \Big)  \\
&\quad + \tilde b_2\Big(  \big(-1- \tfrac{2i S_2}{\omega_0} \big)  \frac{i \omega_0}{\bar u} -2  \frac{\iota }{\bar u} \omega_0  \Big)  + O(\omega_0^{-1}) \\
&= \tilde  b_1 \Big(  -\frac{i \omega_0}{\bar u}   +2  \frac{\iota }{\bar u} \omega_0  \Big)  - \tilde b_2\Big(  \frac{i \omega_0}{\bar u} +2  \frac{\iota }{\bar u} \omega_0  \Big)  + O(\omega_0^{-1}) \\
&=  -\frac{i \omega_0}{\bar u} (\tilde b_1 +\tilde b_2)   +2  \frac{\iota }{\bar u} \omega_0    (\tilde b_1 -\tilde b_2)   + O(\omega_0^{-1}) \\
& =     \frac{2\iota(\sigma_1 + \sigma_2)- i (\sigma_2-\sigma_1) }{2 \sqrt \pi \bar u^2 \iota   }  + O(\omega_0^{-1}) . 
\end{align*}

\emph{Step 9: Combination of \( \mathrm{\mathrm{FP}}\) and \( \mathrm{\mathrm{SP}}\).}
Finally we have
\begin{align*}
\frac{1}{1-\chi_{k^*}^2} \mathrm{\mathrm{FP}}
&= (\frac{1}{2}- \frac{iS_2}{2\omega_0}) \mathrm{\mathrm{FP}} \\
&= \frac{\mathrm{\mathrm{FP}}}{2} -\frac{S_2}{2}\Big( \frac{8 \iota^2 +2 }{\sqrt{\pi} \bar u^2 \big( 4 \iota^2 -1\big) } \Big)  +O(\omega_0^{-1})  
\end{align*}
and 
\begin{equation*}
    \frac{1}{1-\chi_{k^*}^2} \mathrm{\mathrm{SP}}
    =   \frac{2\iota(\sigma_1 + \sigma_2)- i (\sigma_2-\sigma_1) }{4 \sqrt \pi \bar u^2 \iota   }   + O(\omega_0^{-1}). 
\end{equation*}
The real part of this is given by 
\begin{align*}
    \Re\big (\frac{1}{1-\chi_{k^*}^2} \mathrm{\mathrm{FP}} \big) 
    &=   \frac{ -(\sigma_1+\sigma_2) + 2(\Sigma_1+\Sigma_2)}{2\sqrt{\pi} \bar u^2 \big( 4 \iota^2 -1\big) }  - (\Re(\beta_1) + \Re(\beta_2) ) \frac{\iota }{\bar u} \\ 
    &\quad + \frac{\Im(\beta_2)-\Im(\beta_1)}{2\bar u }
    -\frac{S_2}{2}\Big( \frac{8\iota^2+2 }{\sqrt{\pi} \bar u^2 \big( 4\iota^2 -1\big) } \Big)  +O(\omega_0^{-1})
\end{align*}
and 
\begin{equation*}
    \Re\big( \frac{1}{1-\chi_{k^*}^2} \mathrm{\mathrm{SP}} \big)
    =    \frac{\sigma_1 + \sigma_2}{2 \sqrt \pi \bar u^2  }  + O(\omega_0^{-1}). 
\end{equation*}

Substituting the two real parts obtained above into
\[
  \operatorname{Re}\Lambda
   =-\frac{|k^*|^2}{2\sqrt\pi}
     \left[\operatorname{Re}\!\left(\frac{1}{1-\chi_{k^*}^2}\,\mathrm{\mathrm{FP}}\right)
           -2\operatorname{Re}\!\left(\frac{1}{1-\chi_{k^*}^2}\,\mathrm{\mathrm{SP}}\right)\right]
   +O(\omega_0^{-1}),
\]

\emph{Step 10: Simplification. }
To make the final simplification explicit, first combine the real and
imaginary parts of $\beta_1,\beta_2$ from Step~2,3:
\begin{align}
  \Re(\beta_1)+\Re(\beta_2)
   &=\frac{2\iota\Big[(4\iota^2-1)(\sigma_1+\sigma_2-2S_2)+4(\Sigma_1+\Sigma_2)\Big]}
          {\sqrt\pi\,\bar u(4\iota^2-1)^2},
  \label{eq:sum_re_beta}\\[4pt]
  \Im(\beta_2)-\Im(\beta_1)
   &=\frac{-\Big[(4\iota^2-1)(\sigma_1+\sigma_2+2(\Sigma_1+\Sigma_2)-2S_2)+4(\Sigma_1+\Sigma_2)\Big]}
          {\sqrt\pi\,\bar u(4\iota^2-1)^2}.
  \label{eq:diff_im_beta}
\end{align}
Substituting \eqref{eq:sum_re_beta}--\eqref{eq:diff_im_beta} into the expression for
$\Re\big(\tfrac{1}{1-\chi_{k^*}^2}\mathrm{\mathrm{FP}}\big)$ and bringing all
four terms over the common denominator $2\sqrt\pi\,\bar u^2(4\iota^2-1)^2$
gives,
\begin{equation}\label{eq:first_term_step}
  \Re\Big(\frac{1}{1-\chi_{k^*}^2}\,\mathrm{\mathrm{FP}}\Big)
   =\frac{-2\Big[2(\Sigma_1+\Sigma_2)(4\iota^2+1)+(\sigma_1+\sigma_2)\big(8\iota^4+2\iota^2-1\big)\Big]}
          {2\sqrt\pi\,\bar u^2(4\iota^2-1)^2}.
\end{equation}
Every occurrence of $S_2$ cancels identically at this stage: the
$S_2$-contribution of \eqref{eq:sum_re_beta} cancels exactly against that of
\eqref{eq:diff_im_beta} together with the cross term
$-\tfrac{S_2}{2}\big(\tfrac{8\iota^2+2}{\sqrt\pi\bar u^2(4\iota^2-1)}\big)$,
consistent with the $S_2$-independence already observed in
\eqref{eq: gamma_limit_ReLambda}. Moreover the quartic factors as
\begin{equation}\label{eq:quartic_factor}
  8\iota^4+2\iota^2-1=(4\iota^2-1)(2\iota^2+1),
\end{equation}
so that \eqref{eq:first_term_step} reduces to the compact form
\begin{equation}\label{eq:first_term_final}
  \Re\Big(\frac{1}{1-\chi_{k^*}^2}\,\mathrm{\mathrm{FP}}\Big)
   =-\,\frac{2(\Sigma_1+\Sigma_2)(4\iota^2+1)+(\sigma_1+\sigma_2)(4\iota^2-1)(2\iota^2+1)}
          {\sqrt\pi\,\bar u^2(4\iota^2-1)^2}.
\end{equation}
Combining \eqref{eq:first_term_final} with
$\Re\big(\tfrac{1}{1-\chi_{k^*}^2}\mathrm{\mathrm{SP}}\big)
  =\tfrac{\sigma_1+\sigma_2}{2\sqrt\pi\bar u^2}$ over the same denominator, we find 
\begin{align}\label{eq: bracket_final}
\begin{split}
  &\Re\Big(\frac{1}{1-\chi_{k^*}^2}\mathrm{\mathrm{FP}}\Big) 
   -2\Re\Big(\frac{1}{1-\chi_{k^*}^2}\mathrm{\mathrm{SP}}\Big)\\
  &\quad =-\,\frac{2(\Sigma_1+\Sigma_2)(4\iota^2+1)+(\sigma_1+\sigma_2)(4\iota^2-1)\big[(2\iota^2+1)+(4\iota^2-1)\big]}
          {\sqrt\pi\,\bar u^2(4\iota^2-1)^2}
\\[2pt] 
  &\quad =-\,\frac{2(\Sigma_1+\Sigma_2)(4\iota^2+1)+6\iota^2(4\iota^2-1)(\sigma_1+\sigma_2)}
          {\sqrt\pi\,\bar u^2(4\iota^2-1)^2}.
\end{split}
\end{align}
Finally, multiplying \eqref{eq: bracket_final} by the prefactor
$-\dfrac{|k^*|^2}{2\sqrt\pi}$ yields
\begin{equation*}
  \operatorname{Re}\Lambda
   =\frac{|k^*|^2}{\pi\bar u^{2}}\cdot 
     \frac{(4\iota^2+1)(\Sigma_1+\Sigma_2)+3\iota^2(4\iota^2-1)(\sigma_1+\sigma_2)}
          {(4\iota^2-1)^2}
   \;+\;O(\omega_0^{-1})
\end{equation*}
and we obtain~\eqref{eq: gamma_limit_ReLambda} by replacing \(\Sigma_1+\Sigma_2\) by \((\sigma_1+\sigma_2) (1- \tfrac{\hat W_1(2k^*) + \hat W_2(2k^*) }{\hat W_1(k^*) + \hat W_2(k^*) }\big) \).

\end{proof}

Analyzing the numerator~\eqref{eq: gamma_limit_ReLambda} we immediately obtain the following corollary. 
\begin{corollary} \label{cor: high_gamma_cond}
In addition to the assumptions of Corollary~\ref{cor: gamma_limit}, assume that \(\tfrac{\hat W_1(2k^*) + \hat W_2(2k^*)}{\hat W_1(k^*)+\hat W_2(k^*)}<1\), which is equivalent to \(\operatorname{tr}(A(2k^*)) < 0\). Then, the numerator of \(\Re (\Lambda) \) given by~\eqref{eq: gamma_limit_ReLambda} is strictly positive in either of the two cases 
\begin{enumerate}[(i)]
    \item  \(|\iota | > \tfrac{1}{2}\)
    \item   \(0 < |\iota | < \tfrac{1}{2}\) and \( \tfrac{\hat W_1(2k^*) + \hat W_2(2k^*)}{\hat W_1(k^*)+\hat W_2(k^*)}<1- \tfrac{3 \iota^2 (1-4\iota^2)}{4 \iota^2  +1}\)
\end{enumerate}
so that \(\Re(\Lambda) >0\) for large enough \(\gamma\). In the complementary regime \((ii)\) with the reversed inequality we have \(\Re(\Lambda) <0\) for large enough \(\gamma\).
    
\end{corollary}

\begin{remark}
 
Notice that the minimum of \(1- \tfrac{3 \iota^2 (1-4\iota^2)}{4 \iota^2  +1}\) is \((6\sqrt{2}-5)/4 \approx 0.87132\) for \(0 \le |\iota| < 1/2\). 
Hence, if the amplitude of the leading mode \(\hat W_1(k^{*})+\hat W_2(k^{*}) \) is sufficiently large compared to that of the \(\hat W_
1(2k^{*})+\hat W_2(2k^{*}) \) mode (or if they have opposing signs), a supercritical bifurcation is invariably obtained in case \((ii)\) for a sufficiently large non-reciprocity \(\gamma \). 
\end{remark}

\section{Point of critical stability and stability exchange}\label{sec: exchange}

The previous sections establish three ingredients for understanding the dynamics near the homogeneous state. Section~\ref{sec: stability} characterizes the interval of linear stability, Section~\ref{sec: local_bifurcation} establishes the existence of a Hopf branch at the critical parameter \(\alpha_{k^*}\), and Section~\ref{sec: Sub-vs. Supercriticality} determines its criticality through the sign of the cubic normal-form coefficient. However, these results do not yet imply that the system undergoes a Hopf bifurcation from a stable homogeneous state or that the emerging periodic solutions are stable. The missing link is to show that the critical mode \(k^*\) is the first mode to lose stability as \(\alpha\) increases. Otherwise, another mode may already be unstable at \(\alpha_{k^*}\), and the criticality of the Hopf bifurcation alone does not determine the stability of the periodic solutions.

The aim of this section is to establish verifiable conditions under which \(\alpha_{k^*}\) coincides with the first stability threshold \(\alpha_+\). This requires ruling out both possible mechanisms of instability before the critical parameter: a complex conjugate pair of eigenvalues crossing the imaginary axis and a real eigenvalue crossing zero. The maximality condition~\eqref{eq: maximality} ensures that \(k^*\) has the smallest trace threshold, while condition~\eqref{eq: det_positive_interval} excludes a steady instability on \([0,\alpha_{k^*}]\). Lemma~\ref{lem: alpha_plus} combines these conditions to identify the first stability threshold and ensure that the remaining spectrum lies in the left half-plane at the critical point.

Together with the results of the previous sections, this yields the main result of this section, Theorem~\ref{theorem: first_bifurcation}: the homogeneous state is linearly stable for \(\alpha\in[0,\alpha_{k^*})\), and a supercritical Hopf bifurcation at \(\alpha_{k^*}\) gives rise to locally asymptotically stable periodic solutions. We also show that the determinant condition becomes easier to satisfy as the non-reciprocal coupling \(\gamma\) increases, since \(\gamma\) enters the determinant through a non-negative quadratic term. Moreover, Lemma~\ref{lem: pointwise_suffices} provides a convenient sufficient condition under which it is enough to verify~\eqref{eq: det_positive_interval} at \(\alpha_{k^*}\), rather than throughout the interval. This pointwise criterion will be used in Section~\ref{sec: application}. For simplicity, we state the results for \(\alpha\ge0\); the case \(\alpha\le0\) is analogous.

We start by identifying the first bifurcation point. 

\begin{lemma}[Identification of the first bifurcation point]\label{lem: alpha_plus}
Let \(k^*\in\mathcal T^-\) satisfy Assumption~\ref{ass: kernels} together with the maximality
property
\begin{equation}\label{eq: maximality}
  |\hat W_1(k^*)+\hat W_2(k^*)| \;\ge\; |\hat W_1(k)+\hat W_2(k)|
  \qquad\text{for all } k\in\mathcal T^- ,
\end{equation}
and assume
\begin{equation}\label{eq: det_positive_interval}
  \det(A(k))>0 \qquad \text{for all } k\ge1 \text{ and all } \alpha\in[0,\alpha_{k^*}] .
\end{equation}
Then
\[
  \alpha_+ \;=\; T(k^*) \;=\; \alpha_{k^*} \;>\;0 ,
\]
and the infimum \(\inf_{k\in\mathcal T^-}T(k)\) is attained at \(k^*\).

Moreover, there is \(\delta>0\) such that
\[
  \sigma\bigl(\operatorname{D}_uG(\mathbf0,\alpha_{k^*})\bigr)\setminus\{\lambda_\pm(k^*)\}
  \subseteq \{\operatorname{Re}\lambda\le-\delta\};
\]
in particular, \(\lambda_\pm(k^*)\) is the only spectrum of \(\operatorname{D}_uG(\mathbf0,\alpha_{k^*})\)
on the imaginary axis.
\end{lemma}

\begin{proof}
For \(k\in\mathcal T^-\) one has \(\hat W_1(k)+\hat W_2(k)<0\), hence
\(T(k)=(\sigma_1+\sigma_2)/\bigl(\bar u|\hat W_1(k)+\hat W_2(k)|\bigr)>0\), and \(T\) is decreasing
in \(|\hat W_1(k)+\hat W_2(k)|\). By~\eqref{eq: maximality} the infimum over \(\mathcal T^-\) is
therefore attained at \(k^*\), with value \(T(k^*)=\alpha_{k^*}>0\).

By~\eqref{eq: det_positive_interval} no \(D_k\) vanishes on \((0,\alpha_{k^*}]\), so
\(D_+(k)>\alpha_{k^*}\) for every \(k\ge1\) and hence \(\inf_{k\ge1}D_+(k)\ge\alpha_{k^*}\).
Consequently, by Proposition~\ref{prop: linear_stability}
\(\alpha_+=\min\{\inf_{\mathcal T^-}T(k),\ \inf_{k\ge1}D_+(k)\}=\alpha_{k^*}\).

For the last statement we use 
\(\bar u\alpha_{k^*}=-(\sigma_1+\sigma_2)/(\hat W_1(k^*) + \hat W_2(k^*)) \) to obtain 
\[
  \operatorname{tr}A(k)
  =-|k|^2\bigl(\sigma_1+\sigma_2+\bar u\alpha_{k^*} (\hat W_1(k) + \hat W_2(k))\bigr)
  =-|k|^2(\sigma_1+\sigma_2)\Bigl(1-\frac{\hat W_1(k) + \hat W_2(k)}{\hat W_1(k^*) + \hat W_2(k^*)}\Bigr).
\]
If \(\hat W_1(k) + \hat W_2(k) \ge0\) the bracket is greater than one, since \(\hat W_1(k^*) + \hat W_2(k^*)<0\). If \(\hat W_1(k) + \hat W_2(k)<0\), then
\((\hat W_1(k) + \hat W_2(k))/(\hat W_1(k^*) + \hat W_2(k^*)) \le 1\) by~\eqref{eq: maximality}. By~\eqref{eq: cardinality_condition} equality can only hold if \(k=k^*\).

 Hence the bracket is strictly positive for every \(k\neq k^*\), so \(\operatorname{tr}A(k)<0\), while
\(\det A(k)>0\) at \(\alpha_{k^*}\) by~\eqref{eq: det_positive_interval}. Thus every block \(A(k)\) with
\(k\neq k^*\) has both eigenvalues in the open left half-plane. Apply Lemma~\ref{lem: tail} with
\(a:=\alpha_{k^*}\). Every eigenvalue of a mode \(k\ge R(a)\) satisfies \(\operatorname{Re}\lambda\le-c\).
In particular \(k^*<R(a)\), since \(\operatorname{tr}A(k^*)=0\). The remaining modes \(k<R(a)\),
\(k\neq k^*\), are finitely many, so
\[
  \delta:=\min\Bigl\{c,\ \min_{1\le k<R(a),\,k\neq k^*}\bigl|\max\operatorname{Re}\sigma\bigl(A(k)\bigr)\bigr|\Bigr\}>0 .
\]
The claim follows from Lemma~\ref{lemma: block_spectrum}.
\end{proof}

\begin{theorem}[First bifurcation and stability exchange, \(\alpha\ge 0\)] \label{theorem: first_bifurcation}
Suppose there is a mode \(k^*\in \mathcal T^-\) such that~\eqref{eq: maximality} and~\eqref{eq: det_positive_interval} hold. 
Additionally, let the Assumptions~\ref{ass: kernels},~\ref{ass: nonreciprocity} hold for this \(k^*\). 

Then \(\alpha_+>0\) is the first point of critical stability on the positive axis, and the
following hold.
\begin{enumerate}[(a)]
 \item \emph{(Stability window) The homogeneous state \(\boldsymbol{\bar u}\) is linearly stable for every \(\alpha\in[0,\alpha_+)\) and linearly unstable for every \(\alpha>\alpha_+\). }
  \item \emph{(Hopf branch)} At \(\alpha=\alpha_+\) the simple conjugate pair
  \(\lambda_\pm(k^*)=\pm i|k^*|^2\omega_0\), with
  \(\omega_0=\bigl(\bar u^2\gamma^2\hat W(k^*)^2-S_2^2\bigr)^{1/2}>0\),
  is the only spectrum on the imaginary axis and crosses it transversally with speed
  \(\tfrac{d}{d\alpha}\Re\lambda_\pm(k^*)
    =-\tfrac{|k^*|^2}{2}\bar u(\hat W_1(k^*)+\hat W_2(k^*))>0\). A Hopf bifurcation
  occurs in the sense of Proposition~\ref{prop: hopf_bifurcation}: a branch
  \((\boldsymbol u(r),\alpha(r))\) of time-periodic solutions of~\eqref{eq: abstract_problem}
  emanates from \((\boldsymbol{\bar u},\alpha_+)\) with
  \(\alpha(r)=\alpha_++\tfrac12\alpha''(0)\,r^2+O(r^3)\), where \(\alpha''(0)\) is given by equation~\eqref{eq: alpha_second_derivative} and~\eqref{eq: def_Lambda}. Moreover, the frequency is given by \(\kappa(r) = |k^*|^2\omega_0 + \frac{1}{2}\kappa''(0) r^2 + O(r^3)\) with \(\kappa''(0)\) defined by equation~\eqref{eq: def_frequency}. 
  \item \emph{(Exchange of stability)} By the Principle of Exchange of
  Stability~\cite[I.12.2]{kielh2004}, whenever the bifurcation is
  supercritical, i.e.\ \(\Re ( \Lambda ) > 0\), the steady state transfers its stability to the
  emergent orbit: \(\boldsymbol{\bar u}\) loses stability at \(\alpha_+\) and the bifurcating
  time-periodic solution is locally asymptotically stable for
  \(\alpha\in(\alpha_+,\alpha_++\delta_0)\), some \(\delta_0>0\). In particular, if we fix an admissible \(\gamma\), choose \(|\hat W(2k^*)|+ |\hat W_i(2k^*) | \) small enough, and
  \begin{equation*}
     \bigl(\hat W_1(k^*)+\hat W_2(k^*)\bigr)\bigl(B_1\,\hat W_1(k^*)+B_2\,\hat W_2(k^*)\bigr) > 0,
  \end{equation*}
  with \(B_1,B_2\) as in~\eqref{eq: AB_def}, then by
  Corollaries~\ref{cor: sub_subcritical_vanishing_case}--\ref{cor: sub_subcritical_small_case}, we have \(\Re(\Lambda)>0\); the bifurcation is supercritical and the emergent orbit is stable.
\end{enumerate}
If instead \(\hat W_1(k)+\hat W_2(k)\ge 0\) for all \(k\ge 1\) and \(\det A(k)>0\) for all
\(k\ge 1\) and all \(\alpha\ge 0\), then \(\mathcal T^-=\varnothing\), \(\alpha_+=+\infty\),
and \(\boldsymbol{\bar u}\) is stable for every \(\alpha\ge 0\): no bifurcation occurs on the
positive axis.
\end{theorem}
\begin{proof}[Proof of Theorem~\ref{theorem: first_bifurcation}]
    Statement (a) follows from Proposition~\ref{prop: linear_stability}, Lemma~\ref{lem: alpha_plus} and the Principle of Linearised Stability. For this we notice that the Principle of Linearised Stability holds by the same arguments as in~\cite[Theorem~2.7, 2.10]{Carrillo2026} (the part of the proof is independent of the gradient flow structure). Moreover, statement (b) follows from Proposition~\ref{prop: hopf_bifurcation}. Statement (c) follows by Proposition~\ref{prop: sub_supercritical}, Corollary~\ref{cor: sub_subcritical_vanishing_case}, and Corollary~\ref{cor: sub_subcritical_small_case} if we can demonstrate that the Principle of Exchange of Stability holds. But this follows immediately by~\cite[Theorem~I.12.2 and Corollary~I.12.3]{kielh2004}. 
\end{proof}

The next lemma provides easily verifiable conditions that ensure~\eqref{eq: det_positive_interval}. In particular, for our application in Section~\ref{sec: application}, these conditions considerably simplify the verification of~\eqref{eq: det_positive_interval}.
\begin{lemma}\label{lem: pointwise_suffices}
Fix \(k\ge1\). Suppose \(D_k(\alpha_{k^*})>0\) and~\eqref{eq: maximality} for \(k^*\). Then already
\(D_k(\alpha)>0\) for all \(\alpha\in[0,\alpha_{k^*}]\), provided that either
\(\widehat W_1(k)\widehat W_2(k)\le0\) or \(\sigma_1=\sigma_2\).
\end{lemma}

\begin{proof}
If \(\widehat W_1(k)\widehat W_2(k)\le0\), then by Remark~\ref{rem: partition_structure} the set
\(\{D_k>0\}\) is an interval containing the origin, and the claim is immediate.

So let \(\widehat W_1(k)\widehat W_2(k)>0\). If both coefficients are positive, \(D_k\) has no
positive root and there is nothing to prove, so assume \(\widehat W_1(k),\widehat W_2(k)<0\) and
put
\[
  a_i:=-\bar u\widehat W_i(k)>0,
  \qquad
  b_i:=\frac{\sigma_i}{a_i},
\]
so that \(\det(A(k))= D_k(\alpha)=a_1a_2(\alpha-b_1)(\alpha-b_2)+\bar u^2\gamma^2\widehat W(k)^2\). 
Since \(a_1,a_2  > 0\), \(D_k(\alpha)\) is a parabola with vertex at \(\tfrac{1}{2}(b_1+b_2)\). Hence, on the \([0,\tfrac{1}{2}(b_1+b_2)]\) we know \(D_k(\alpha)\) is non-increasing. It is enough to show
\begin{equation}\label{eq: mediant_bound}
  \alpha_{k^*} \le\ \tfrac12(b_1+b_2) .
\end{equation}

By maximality~\eqref{eq: maximality} and \(\sigma_1+\sigma_2=a_1 b_1+a_2b_2\),
\[
  \alpha_{k^*} =- \frac{\sigma_1+\sigma_2}{\bar u (\hat W_1(k^*) + \hat W_2(k^*))}
  \ \le\ \frac{\sigma_1+\sigma_2}{a_1+a_2}
  \ =\ \frac{a_1b_1+a_2 b_2}{a_1+a_2} ,
\]
the weighted average of \(b_1\) and \(b_2\) with weights \(a_1,a_2\). Subtracting the midpoint,
\[
  \frac{a_1 b_1+a_2 b_2}{a_1+a_2}-\frac{b_1+b_2}{2}
  =\frac{(a_1-a_2)(b_1-b_2)}{2(a_1+a_2)} .
\]
For \(\sigma_1=\sigma_2=\sigma\) one has \(b_1-b_2=\sigma(1/a_1-1/a_2)\), which carries the sign \(a_2-a_1\). Hence the right-hand side is \(\le0\) and \eqref{eq: mediant_bound} follows.
\end{proof}

We want to demonstrate that in the symmetric case Lemma~\ref{lem: pointwise_suffices} is trivially satisfied and, therefore, the remaining spectrum lies in the negative half space by Lemma~\ref{lem: alpha_plus}. 

\begin{corollary}[Symmetric case]
Assume that
\(
\sigma_1=\sigma_2=\sigma,  
\)
and the maximal inequality~\eqref{eq: maximality} for \(k^* \in \mathcal T^-  \) as well as~\eqref{eq: cardinality_condition}. 
Then, the condition \(D_k(\alpha_{k^*}) > 0\) of Lemma~\ref{lem: pointwise_suffices} is satisfied if 
\[
\min( \hat W_1(k),\hat W_2(k)) > \frac{\hat W_1(k^*)+ \hat W_2(k^*)}{2}
.\]
In particular, \(D_k(\alpha_{k^*}) > 0\)  if additionally \(W_1=W_2\). 
\end{corollary}
\begin{proof}
Under the assumptions we have \( S_1(k) = \sigma + \alpha_{k^*} \bar u \hat W_1(k), S_2(k) = \sigma + \alpha_{k^*} \bar u \hat W_2(k) \)
and let \(k^*\) denote the mode satisfying~\eqref{eq: maximality}.
At the critical value
\[\alpha^*_k
=-\frac{\sigma_1+\sigma_2}{ \bar u (\hat W_1(k^*)+ \hat W_2(k^*))}
=-\frac{2 \sigma}{ \bar u (\hat W_1(k^*)+ \hat W_2(k^*))},\]
the determinant simplifies to
\[\det(A(k))
= |k|^4
\left(\sigma^2\left(1-\frac{2 \hat W_1(k)}{  (\hat W_1(k^*)+ \hat W_2(k^*))}
\right)\left( 1-\frac{ 2 \hat W_2(k)}{  (\hat W_1(k^*)+ \hat W_2(k^*))} \right)+\bar u^2\gamma^2\widehat W(k)^2
\right). \]

If \[\min(\widehat W_1(k),\widehat W_2(k))>\tfrac12(\widehat W_1(k^*)+\widehat W_2(k^*))\], both factors in the first term are positive, hence \(D_k(\alpha_{k^*})>0\). If \(W_1=W_2\), this holds for every \(k\neq k^*\) by
\eqref{eq: maximality} and~\eqref{eq: cardinality_condition}, while for \(k=k^*\) the first term vanishes and
\(D_{k^*}(\alpha_{k^*})=\bar u^2\gamma^2\widehat W(k^*)^2>0\) by~\eqref{eq: imaginary_condition}.
\end{proof}

\section{Application} \label{sec: application}
In this section our goal is to apply the results of the previous sections to a specific example. We consider the Kuramoto model with non-reciprocal coupling, the two-mode cosine example, as well as the top hat kernel.

\begin{example} \label{example: kuramoto}
Let \(W (x) = W_1(x) = W_2(x) = - \cos(x) \). Then, \(\hat W(k) = 0\) for all \(k \neq  1 \) and \(\hat W(k) = - \pi\) for \(k =1 \). All assumptions of Proposition~\ref{prop: hopf_bifurcation} are satisfied. Indeed, for simplicity assume \(\sigma_1 = \sigma_2  = 1 \) and let \(\gamma >  0\). Then~\eqref{eq: imaginary_condition} holds. Together with the trivial conditions~\eqref{eq: cardinality_condition},~\eqref{eq: determinante_non_vanishing} all assumptions of Proposition~\ref{prop: hopf_bifurcation} hold for the bifurcation point \(\alpha^* = \alpha_+= 2  \). The eigenvalues are given by 
\(\lambda_\pm^* = \pm i \frac{\gamma}{2 } \). Hence, the bifurcation curve \(u(r)\) is given by 
\begin{equation*}
    u(r)(t,x)  = \frac{1}{2 \pi} + \frac{2r}{\sqrt{\pi}} 
    \begin{pmatrix}
        \cos\Big( \frac{\gamma}{2 } t \Big) \cos(x) \\
        -\sin \Big( \frac{\gamma}{2 } t \Big) \cos(x)
    \end{pmatrix}
    + O(r^2).
\end{equation*}
The bifurcation branch is given by
\begin{equation}
    \alpha(r) =  2   + \frac{6 \pi (\gamma^2+8)}{ (\gamma^2+16)} r^2 + O(r^3)
\end{equation} 
with frequency 
\begin{equation*}
    \kappa(r) = \frac{\gamma}{2} - \frac{6 \gamma \pi}{ \gamma^2  + 16} r^2    + O(r^3)
\end{equation*}
in a small interval around the steady state.  
Additionally, Theorem~\ref{theorem: first_bifurcation} holds, i.e. the exchange of stability. 
\end{example}

Similarly, we can investigate a linear combination of \(\cos\) kernels. 
\begin{example} \label{example: cos}
Let \(\sigma_1 =\sigma_2 = 1\), \(a_1, a_2 \in \R\) and \(W(x) = W_1(x) = W_2(x)  = - a_1 \cos (x)  - a_2 \cos(2x) \) with \(a_1 > \max(a_2,0)  \). Then, \(\hat W(1) = - a_1 \pi \), \(\hat W(2) =- a_2 \pi \) and \(\hat W(k) = 0 \) for \(k \ge 3 \). All assumptions of Proposition~\ref{prop: hopf_bifurcation} are satisfied similar to Example~\ref{example: kuramoto}. Then the first bifurcation point is given by \(\alpha^* = \tfrac{2}{a_1}\). The eigenvalues are given by 
\(\lambda_\pm^* = \pm i \frac{a_1 \gamma }{2 } \). The bifurcation curve \(u(r)\) is given by 
\begin{equation*}
    u(r)(t,x) = \frac{1}{2 \pi} + 2r e_1(x)
    \begin{pmatrix}
        \cos\Big( \frac{ a_1 \gamma}{2 } t \Big) \\
        -\sin \Big( \frac{a_1 \gamma}{2 } t \Big) 
    \end{pmatrix}
    + O(r^2).
\end{equation*}

Moreover, if \(a_2/a_1\in (1/2,1)\), there exists \(\gamma^* > 0\) by Corollary~\ref{cor: nonreciprocal_transition} such that the bifurcation is subcritical for \(\gamma < \gamma^*\) and supercritical for \(\gamma > \gamma^*\). 
\end{example}

\begin{example}[Top-hat kernel]
Let \(\sigma_1 =\sigma_2 = 1\), \(R < \pi/5 \) and \(W(x) = W_1(x) = W_2(x)  = \tfrac{1}{2R} \indicator{|x|\le R} \).  Then, \(W \in BV(\T)\) and \(\hat W(k) = \tfrac{1}{kR} \sin(kR)\). 
Moreover, \(\hat W(k)=\dfrac{\sin(kR)}{kR}\); since \(kR<\pi\) for
\(k\le5\), one has \(\hat W(1),\dots,\hat W(5)>0\), and the unique minimiser
\(k^*=\operatorname{arg\,min}_{k\ge1}\hat W(k)\) satisfies \(k^*\ge7\), \(k^*R\in(\pi,2\pi)\) with 
\(\hat W(k^*)<0\).
For every \(\gamma>0\) the assumptions of Proposition~\ref{prop: hopf_bifurcation} hold, and
the first bifurcation on the positive axis is the Hopf point
\[
  \alpha^*=-\frac{1}{\bar u\,\hat W(k^*)}=-\frac{2\pi\,k^*R}{\sin(k^*R)}>0,
  \qquad
  \omega_0=\frac{\gamma\,|\sin(k^*R)|}{2\pi\,k^*R},\qquad
  \lambda_\pm(k^*)=\pm i\,|k^*|^2\omega_0 ,
\]
with eigenvector \((e_{k^*},i\,e_{k^*})\) and curve
\[
  u(r)(t)=\frac{1}{2\pi}+2r\,e_{k^*}(x)
  \begin{pmatrix}\cos(|k^*|^2\omega_0 t)\\ - \sin(|k^*|^2\omega_0 t)\end{pmatrix}+O(r^2).
\]
Setting \(\iota=\hat W(2k^*)/\hat W(k^*)=\cos(k^*R)<\tfrac12\), Corollary~\ref{cor: lambda_matching}
gives \(N_R(\omega_0^2)>0\), hence
\[
  \alpha''(0)=-\frac{2\pi N_R(\omega_0^2)}{\hat W(k^*)\,D_{b\tilde b}(\omega_0^2)}>0 . 
\]
This implies \(\Re(\Lambda) > 0 \) and
the bifurcation is supercritical for every \(\gamma>0\), and by the Principle of Exchange of
Stability the emergent periodic orbit is stable. 
\end{example}

\medskip
\noindent\textbf{AI disclosure.} The mathematical content of this paper was developed entirely by
the authors. AI tools (Claude, Anthropic) were used for proofreading, grammatical and linguistic
editing, checking some of the calculations, and generating the figures.

\appendix

\section{Useful identities} \label{sec:appendix}
In this section we collect some useful identities. Recall the definitions of \(\alpha_{k^*} \) given by~\eqref{eq: alpha_bifurcation_point}, \(S_2\)
given by~\eqref{eq: def_S_2} and \(\chi_{k^*}\) given by~\eqref{eq: def_chi}. 
\begin{table}[h]
\centering
\renewcommand{\arraystretch}{2.0}
\begin{tabular}{@{}c l@{}}
\toprule
\textbf{Label} & \textbf{Relation} \\
\midrule
\refstepcounter{equation}\label{eq: aux_special_case_lyapunov}
\multirow{2}{*}{(\theequation)}
  & $\displaystyle \alpha_{k^*} \hat W_1(k^*) = -\frac{\sigma_1 + S_2}{\bar u}$ \\[6pt]
  & $\displaystyle \gamma \chi_{k^*} \hat W(k^*) = \frac{S_2 - i\omega_0}{\bar u}$ \\[4pt]
\midrule
\refstepcounter{equation}\label{eq: square_root_relation}(\theequation)
  & $\displaystyle \bar u^2 \big(\gamma \hat W(k^*)\big)^2 = \omega_0^2 + S_2^2$ \\[4pt] 
\midrule
\refstepcounter{equation}\label{eq: aux_real_part_chi_gamma1}(\theequation)
  & $\displaystyle \gamma \Re(\chi_{k^*}) \hat W(k^*) = \frac{S_2}{\bar u}$ \\[4pt]
\refstepcounter{equation}\label{eq: aux_real_part_chi_gamma2}(\theequation)
  & $\displaystyle \alpha_{k^*} \hat W_2(k^*) - \gamma \Re(\chi_{k^*}) \hat W(k^*) = -\frac{\sigma_2}{\bar u}$ \\[4pt]
\refstepcounter{equation}\label{eq: aux_real_part_chi_gamma3}(\theequation)
  & $\displaystyle \alpha_{k^*} \hat W_1(k^*) + \gamma \Re(\chi_{k^*}) \hat W(k^*) = -\frac{\sigma_1}{\bar u}$ \\[4pt]
\bottomrule
\end{tabular}
\caption{Auxiliary relations at the critical wave number $k^*$. The quantity
$\chi_{k^*} = \dfrac{S_2 - i\omega_0}{\bar u \gamma \hat W(k^*)}$ enters via
\eqref{eq: aux_special_case_lyapunov}.}
\label{tab: aux_relations}
\end{table}

\section{An abstract inverse theorem for block--diagonal operators}

In this section we justify that we can calculate the inverse of our linearized operator \(\operatorname{D}G\) on each \(2 \times 2 \) block in our formula for \(\Lambda\).  
Throughout this subsection \(\mathcal{H}\) denotes a separable complex Hilbert
space that admits an orthogonal Hilbert--space direct sum decomposition into
finite--dimensional subspaces,
\begin{equation*}
    \mathcal{H} = \bigoplus_{k \in \mathcal{K}} H_k ,
    \qquad
    \mathcal{K} \text{ countable}, \quad \dim H_k < \infty, \quad
    H_k \perp H_j \ (k \neq j).
\end{equation*}
We write \(P_k \colon \mathcal{H} \to H_k\) for the orthogonal projection and
\(x_k := P_k x\), so that \(x = \sum_{k} x_k\) and
\(\norm{x}^2 = \sum_{k} \norm{x_k}^2\). For a linear map \(T\) on a
finite--dimensional space we denote by
\(\sigma_{\min}(T) := \inf_{\norm{v}=1}\norm{Tv}\) its smallest singular value. If \(T\) is invertible then \(\norm{T^{-1}} = \sigma_{\min}(T)^{-1}\).

Given a family \(A_k \in \mathcal{L}(H_k)\), \(k \in \mathcal{K}\), we define the block--diagonal operator
\begin{equation} \label{eq: maximal_block_operator}
    D(A) = \Big\{ x \in \mathcal{H} :
        \sum_{k \in \mathcal{K}} \norm{A_k x_k}^2 < \infty \Big\},
    \qquad
    A x = \sum_{k \in \mathcal{K}} A_k x_k .
\end{equation}

\begin{theorem} \label{theorem: block_inverse}
    Let \(A\) be the operator defined in~\eqref{eq: maximal_block_operator}.
    Then \(A\) is densely defined and closed, and the following are equivalent:
    \begin{enumerate}[(i)]
        \item \(A \colon D(A) \to \mathcal{H}\) is a bijection with bounded inverse;
        \item 
            \(\inf\limits_{k \in \mathcal{K}} \sigma_{\min}(A_k) > 0\).
    \end{enumerate}
    In this case \(A^{-1} \in \mathcal{L}(\mathcal{H})\) is block--diagonal,
    \begin{equation*}
        A^{-1} x = \sum_{k \in \mathcal{K}} A_k^{-1} x_k ,
        \qquad
        \norm{A^{-1}}
        = \sup_{k \in \mathcal{K}} \norm{A_k^{-1}}
        = \frac{1}{\displaystyle \inf_{k \in \mathcal{K}} \sigma_{\min}(A_k)} .
    \end{equation*}
\end{theorem}

\begin{proof}
\emph{Density. } Let \(\mathcal{D}_0 = \operatorname{span}\{x : \; \exists k \in \mathcal K \; \mathrm{such \; that \; } x \in H_k \}\) be the space of
    finitely supported vectors. Since each such vector satisfies
    \(\sum_k \norm{A_k x_k}^2 < \infty\) trivially, \(\mathcal{D}_0 \subseteq D(A)\);
    as the \(H_k\) orthogonally span \(\mathcal{H}\), the domain \(D(A)\) is dense.

    \emph{Closedness.} Let \(x^{(n)} \to x\) and \(A x^{(n)} \to y\) in
    \(\mathcal{H}\). Fix \(k \in \mathcal K\). Continuity of \(P_k\) gives
    \(x_k^{(n)} = P_k x^{(n)} \to x_k\), and since \(A\) is block--diagonal,
    \(P_k A x^{(n)} = (Ax^{(n)})_k=  A_k x_k^{(n)} \to y_k\). As \(A_k\) is continuous on the finite--dimensional space \(H_k\), also \(A_k x_k^{(n)} \to A_k x_k\), where we have  
    \(y_k = P_k y =A_k x_k \) for every \(k\). Therefore
    \[\sum_k \norm{A_k x_k}^2 = \sum_k \norm{y_k}^2 = \norm{y}^2 < \infty,\]
    so \(x \in D(A)\) and \(A x = y\). Thus \(A\) is closed.

    \emph{(ii) \(\Rightarrow\) (i).} Set
    \(C := \sup_k \norm{A_k^{-1}}_{\mathrm{op}} < \infty\) (finite by (ii), since
    \(\norm{A_k^{-1}}_{\mathrm{op}} = \sigma_{\min}(A_k)^{-1}\)) and define
    \(B x := \sum_k A_k^{-1} x_k\). By orthogonality,
    \begin{equation*}
        \norm{Bx}^2 = \sum_k \norm{A_k^{-1} x_k}^2
        \le C^2 \sum_k \norm{x_k}^2 = C^2 \norm{x}^2 ,
    \end{equation*}
    so the series converges and \(B \in \mathcal{L}(\mathcal{H})\) with
    \(\norm{B}_{\mathrm{op}} \le C\). For \(x \in \mathcal{H}\) we have
    \((Bx)_k = A_k^{-1} x_k\), hence
    \[\sum_k \norm{A_k (Bx)_k}^2 = \norm{x}^2 
    < \infty, \] i.e. \(Bx \in D(A)\), and
    \[A(Bx) = \sum_k A_k A_k^{-1} x_k = x.\] For \(x \in D(A)\),
    \[(BAx)_k = A_k^{-1} A_k x_k = x_k\], so \(BAx = x\). Thus \(A\) is a
    bijection with \(A^{-1} = B \in \mathcal{L}(\mathcal{H})\).

    \emph{(i) \(\Rightarrow\) (ii).} Fix \(k\). If \(A_k v = 0\) for
    \(v \in H_k\), then embedding \(v\) in \(\mathcal{H}\) gives \(v \in D(A)\) and
    \(A v = A_k v = 0\). 
    Injectivity of \(A\) forces \(v = 0\), so \(A_k\) is
    injective, hence invertible on the finite--dimensional subspace \(H_k\). 
    
    For \(x_k \in H_k\), the vector \(w := A_k^{-1} x_k \in H_k \subseteq D(A)\)
    satisfies \(A w = A_k w = x_k\), so by uniqueness \(A^{-1} x_k = A_k^{-1} x_k\)
    and \(\norm{A_k^{-1} x_k} = \norm{A^{-1} x_k} \le \norm{A^{-1}}
    \norm{x_k}\). Hence \(\norm{A_k^{-1}} \le
    \norm{A^{-1}}\) for all \(k\), which is (ii).

    \emph{Norm.} The bound \(\norm{B} \le C\) together with the
    lower bound \(\norm{A^{-1}} \ge \norm{A_k^{-1}}\)
    for every \(k\) gives \(\norm{A^{-1}} = \sup_k
    \norm{A_k^{-1}} = (\inf_k \sigma_{\min}(A_k))^{-1}\).
\end{proof}

\begin{remark} \label{rem: singular_values_not_eigenvalues}
    The criterion in Theorem~\ref{theorem: block_inverse} is stated in terms of the
    smallest singular values \(\sigma_{\min}(A_k)\), \emph{not} the eigenvalues of
    \(A_k\). For non--normal \(A_k\) these two quantities can differ arbitrarily,
    so a lower bound on \(\lvert \lambda \rvert\) over the eigenvalues
    \(\lambda \in \sigma(A_k)\) is in general \emph{not} sufficient for bounded
    invertibility. The theorem uses only the orthogonal block structure and
    requires no normality, self--adjointness, or spectral theorem, which is exactly the setting of Section~\ref{sec: stability}. 
\end{remark}

\subsection{Application to the linearized two--species operator}

We specialize Theorem~\ref{theorem: block_inverse} to \(A_0 = D_u G(\boldsymbol{0},
\alpha)\) on the symmetric, mean--zero space. Recall \(d = 1\), the torus
\([-\pi, \pi]\), and the orthonormal basis \(e_k = \sqrt{1/\pi}\,\cos(k \cdot)\),
\(k \ge 1\) and
\begin{equation*}
    \dot L^2_s(\T) = \{ f \in L^2_s(\T) : \hat f(0) = 0 \}, \qquad
    \dot H^2_s(\T) = H^2(\T) \cap \dot L^2_s(\T), \qquad
    \mathcal{H} = \dot L^2_s(\T) \times \dot L^2_s(\T) ,
\end{equation*}
and, for \(k \ge 1\), \(H_k = \operatorname{span}\{(e_k, 0), (0, e_k)\} \)
(an isometric identification). By construction
\(\mathcal{H} = \bigoplus_{k \ge 1} H_k\) and \(A_0\rvert_{H_k} = A(k) = -k^2 M(k)\)
with
\begin{equation*}
    M(k) = 
    \begin{pmatrix}
        \sigma_1 + \bar u \alpha \hat W_1(k) & \bar u \gamma\, \hat W(k) \\[2pt]
        -\bar u \gamma\, \hat W(k) & \sigma_2 + \bar u \alpha \hat W_2(k)
    \end{pmatrix}
\end{equation*}

Throughout we have our standing assumption \(\sigma_1, \sigma_2, \gamma,  \bar u > 0\),
\(\alpha \in \R\), and that \(W, W_1, W_2 \in \mathrm{BV}(\T) \) are even and real--valued. Two consequences are used below: the
uniform bound \(\sup_{k \ge 1} \lvert \hat W_i(k) \rvert \le \norm{W_i}_{L^1} 
< \infty\), and, by integration by parts (the boundary terms vanish since
\(\sin(k\pi) = 0\)), the decay
\begin{equation} \label{eq: BV_decay}
    \lvert \hat W_i(k) \rvert \le \frac{V(W_i)}{k}, \qquad
    \lvert \hat W(k) \rvert \le \frac{V(W)}{k},
\end{equation}
with \(V(\cdot)\) the total variation. In particular
\(M(k) \to M_\infty = \operatorname{diag}(\sigma_1, \sigma_2)\) as \(k \to \infty\)
and \(C_M := \sup_{k \ge 1} \norm{M(k)}_{\mathrm{op}} < \infty\).

\begin{lemma} \label{lem: domain_identification}
    Under the above assumptions, the maximal domain
    \eqref{eq: maximal_block_operator} of \(A_0\) satisfies
    \(D(A_0) = \dot H^2_s \times \dot H^2_s\).
\end{lemma}

\begin{proof}
    In the cosine basis \(x \in \dot H^2_s \times \dot H^2_s\) if and only if
    \(\sum_{k \ge 1} k^4 \norm{x_k}^2 < \infty\), and
    \(\norm{A(k) x_k} = k^2 \norm{M(k) x_k}\). Since
    \(\norm{M(k) x_k} \le C_M \norm{x_k}\), we get
    \(\dot H^2_s \times \dot H^2_s \subseteq D(A_0)\). 
    
    Conversely,
    \(\sigma_{\min}(M(k)) \to \sigma_{\min}(M_\infty) = \min(\sigma_1, \sigma_2) >
    0\), so there are \(c_0 > 0\) and \(R < \infty\) with
    \(\norm{M(k) x_k} \ge c_0 \norm{x_k}\) for \(k \ge R\). Hence
    \[\sum_{k \ge R} k^4 \norm{x_k}^2 \le c_0^{-2} \sum_{k \ge R} k^4
    \norm{M(k) x_k}^2 < \infty\] for \(x \in D(A_0)\), while the finitely many
    terms \(1 \le k < R\) form a finite sum. Thus
    \(D(A_0) \subseteq \dot H^2_s \times \dot H^2_s\).
\end{proof}

\begin{corollary} \label{cor: A0_invertible}
    In addition to the standing assumptions, suppose the non-degeneracy condition~\eqref{eq: determinante_non_vanishing} holds.
    Then \(A_0 \colon \dot H^2_s \times \dot H^2_s \to \dot L^2_s \times
    \dot L^2_s\) is a bijection with bounded inverse, given by
    \begin{equation*}
        A_0^{-1} x = \sum_{k \ge 1} A(k)^{-1} \langle x, e_k  \rangle \, e_k,
        \qquad
        A(k)^{-1}
        = -\frac{1}{k^2 \det M(k)}
        \begin{pmatrix}
            \sigma_2 +\alpha \bar u \hat W_2(k)  & -\bar u \gamma \hat W(k) \\ \bar u \gamma \hat W(k) & \sigma_1 +\alpha \bar u \hat W_1(k)
        \end{pmatrix},
    \end{equation*}
    and
    \begin{equation*}
        \norm{A_0^{-1}}
        = \sup_{k \ge 1} \frac{1}{k^2\, \sigma_{\min}(M(k))} < \infty .
    \end{equation*}
\end{corollary}

\begin{proof}
    By Lemma~\ref{lem: domain_identification}, \(A_0\) is the
    block--diagonal operator of Theorem~\ref{theorem: block_inverse} with
    \(\mathcal{K} = \{ k \in \N : k \ge 1 \}\) and \(A_k = A(k) = -k^2 M(k)\).
    We verify condition (ii), i.e.\ \(\inf_{k \ge 1} \sigma_{\min}(A(k)) > 0\),
    using \(\sigma_{\min}(A(k)) = k^2 \sigma_{\min}(M(k))\) and the
    \(2 \times 2\) identity
    \begin{equation*}
        \sigma_{\min}(M(k))
        = \frac{\lvert \det M(k) \rvert}{\norm{M(k)}_{\mathrm{op}}}
        \ge \frac{\lvert \det M(k) \rvert}{C_M} .
    \end{equation*}
    
    \emph{Tail \(k \ge R\).} Since \(\det M(k) \to \sigma_1 \sigma_2\), we fix \(R\)
    with \(\lvert \det M(k) \rvert \ge \tfrac12 \sigma_1 \sigma_2\) for
    \(k \ge R\) (possible by \eqref{eq: BV_decay}). Then
    \[\sigma_{\min} (M(k)) \ge \sigma_1 \sigma_2 / (2 C_M) =: c_0 > 0,\]
    so
    \(\sigma_{\min}(A(k)) = k^2 \sigma_{\min}(M(k)) \ge c_0\).
    
    \emph{Low modes \(1 \le k < R\).} By \eqref{eq: determinante_non_vanishing}, \(\det M(k) \neq 0\). This implies the minimum 
    \[\tilde c_0 := \min_{1 \le k < R} \sigma_{\min}(M(k)) > 0\] exists and 
    \(\sigma_{\min}(A(k)) \ge k^2 \tilde c_0  \ge \tilde c_0 \).
    
    Therefore, we find 
    \(\inf_{k \ge 1} \sigma_{\min}(A(k)) \ge \min(c_0, \tilde c_0 ) > 0\), and
    Theorem~\ref{theorem: block_inverse} yields the bounded inverse, the block
    formula, and the norm. The explicit inverse follows from
    \(A(k)^{-1} = (-k^2 M(k))^{-1} = -k^{-2} M(k)^{-1}\) and the
    \(2 \times 2\) inverse formula.
\end{proof}

\begin{remark} \label{rem: justifies_modewise}
    Corollary~\ref{cor: A0_invertible} is what licenses the mode--wise inversion
    used in the computations: inverting \(A(k)\) on a single Fourier mode \(e_k\)
    produces the blocks \(A(k)^{-1}\), and the uniform bound
    \(\inf_{k \ge 1} \sigma_{\min}(A(k)) > 0\) guarantees that these assemble
    into a genuine bounded operator \(A_0^{-1}\) on all of \(\dot L^2_s \times
    \dot L^2_s\), mapping into \(\dot H^2_s \times \dot H^2_s\). Without this
    uniform lower bound the mode--wise inverses still exist, but the range of
    \(A_0\) need not be closed and \(A_0^{-1}\) need not be bounded.
\end{remark}

\subsection{Block spectrum}
In this section our goal is to demonstrate that the formal 
calculations of diagonalizing our linear operator by the Fourier basis 
\(\boldsymbol{e_k}\) are indeed justified and the spectrum of  the linear operator is the union of spectra. Notice,
this does not follow immediately, since our \(A(k)\) are generally (besides the matching case) not normal
by the non-reciprocity of our system. 

\begin{lemma} \label{lemma: block_spectrum}
    Let \(A = \bigoplus_{k \ge 1} A(k)\) be the block--diagonal operator
    on \(\mathcal{H} = \bigoplus_{k \ge 1} H_k\), where each \(H_k\) is
    finite--dimensional and each \(A(k) \in \mathcal{L}(H_k)\), and suppose that
    \(A\) has compact resolvent. Then
    \begin{equation*}
        \sigma (A) = \sigma_p(A) = \bigcup_{k \ge 1} \sigma\big(A(k)\big).
    \end{equation*}
    In particular the closure is redundant, since \(\bigcup_{k \ge 1}
    \sigma(A(k))\) is already closed.
\end{lemma}

\begin{proof}
    \emph{Point spectrum.} Let \(\mu \in \C\). For \(x = (x_k)_{k \ge 1} \in D(A)\)
    the block--diagonal action gives \((A x)_k = A(k) x_k\), so
    \begin{equation*}
        A x = \mu x
        \quad \Longleftrightarrow \quad
        A(k) x_k = \mu x_k \ \text{ for every } k \ge 1 .
    \end{equation*}
    If \(x \neq 0\) then some component \(x_{k_0} \neq 0\), and
    \(A(k_0) x_{k_0} = \mu x_{k_0}\) shows \(\mu \in \sigma_p(A(k_0))\).
    Conversely, if \(A(k_0) v = \mu v\) for some \(v \neq 0\), the vector
    supported in block \(k_0\) by \(v\) lies in \(D(A)\) and is an eigenvector of
    \(A\) for \(\mu\). Hence \(\mu \in \sigma_p(A)\) if and only if
    \(\mu \in \sigma_p(A(k))\) for some \(k \ge 1\). Since each \(H_k\) is
    finite--dimensional, \(\sigma_p(A(k)) = \sigma(A(k))\), which yields
    \begin{equation*}
        \sigma_p(A) = \bigcup_{k \ge 1} \sigma\big(A(k)\big) .
    \end{equation*}

    \emph{Full spectrum.} As \(A\) has compact resolvent, the Riesz--Schauder
    theory~\cite[Corollary~4.1.19]{engel2000} gives \(\sigma(A) = \sigma_p(A)\) \footnote{Notice that \(\operatorname{range}(\mu-A)\) is closed for every \(\mu\) by the Fredholm alternative and the compactness of the resolvent. Hence, in the terminology of Engel and Nagel, the approximate point spectrum coincides with the point spectrum in our setting.}, a set of isolated eigenvalues of
    finite algebraic multiplicity with no finite accumulation point. Combined
    with the previous step,
    \begin{equation*}
        \sigma(A) = \sigma_p(A) = \bigcup_{k \ge 1} \sigma\big(A(k)\big) .
    \end{equation*}
    Finally, \(A\) is closed (its resolvent set is non-empty), so \(\sigma(A)\) is
    closed; a set equal to a closed set equals its own closure, whence
    \begin{equation*}
        \sigma(A) = \overline{\bigcup_{k \ge 1} \sigma\big(A(k)\big)} .
    \end{equation*}
\end{proof}

\begin{remark} \label{rem: nonnormal_harmless}
    No normality of the blocks \(A(k)\) is used. Without the compact resolvent
    assumption only the inclusion \(\overline{\bigcup_{k \ge 1} \sigma(A(k))}
    \subseteq \sigma(A)\) holds in general, and it may be strict for non--normal
    blocks, namely at points \(\mu\) with \(\sup_{k} \norm{(A(k) - \mu)^{-1}} =
    \infty\) that are eigenvalues of no single block. The compact resolvent
    removes this possibility by reducing the spectrum to the point spectrum,
    where the block decoupling above is exact.
\end{remark}

\section{Reciprocal system} \label{app: reciprocal_caluclations}

In this subsection we provide a short computation of the second derivative \(\alpha''(0)\) of the bifurcation branch in the reciprocal case. More precisely, we consider the system
\begin{align}\label{eq: 2species_problem_reciprocal}
\begin{cases}
    \partial_t u_{1}
    = \sigma \Delta u_{1}
    + \nabla \cdot \Big(
        u_{1}
        \big(
            \alpha \nabla (W * u_1)
            + \gamma \nabla (W * u_2)
        \big)
      \Big),\\[0.2cm]
    \partial_t u_{2}
    = \sigma \Delta u_{2}
    + \nabla \cdot \Big(
        u_{2}
        \big(
            \alpha \nabla (W * u_2)
            + \gamma \nabla (W * u_1)
        \big)
      \Big),
\end{cases}
\end{align}
for constants \(\sigma,\gamma>0\), a parameter \(\alpha\in\R\), and an even interaction kernel \(W\) on the torus \(\T=[-\pi,\pi]\).
Similar to the main section we write \(F: \dot H_s^2(\T)\times \dot H_s^2(\T) \to L^2_s(\T) \times L_s^2(\T)\) for the abstract problem. 

Our goal is to evaluate the bifurcation formula~\cite[(I.6.11)]{kielh2004} for \(\alpha''(0)\). Diagonalizing the linearized operator with respect to the Fourier basis yields the matrices
\begin{equation*}
A(k)
=
-|k|^2
\begin{pmatrix}
\sigma + \bar u \alpha \widehat W(k)
&
\bar u \gamma \widehat W(k)
\\
\bar u \gamma \widehat W(k)
&
\sigma + \bar u \alpha \widehat W(k)
\end{pmatrix}.
\end{equation*}
Following the computations in~\cite{Carrillo2026}, we obtain the critical parameter value
\begin{equation*}
\alpha^*
=
-\frac{\sigma}{\bar u \widehat W(k^*)}
-\gamma,
\end{equation*}
where we assume that \(\hat W(k^*) \le \hat W(k) \) for all \(k \neq k^*\), \(k \ge 1 \), i.e. \(\alpha^*\) is the first eigenvalue crossing on the positive axis. 
At \(\alpha=\alpha^*\), the linearized operator \(\operatorname{D}_{v}F(0, \alpha^*) \) has a one-dimensional kernel
\[
N(\operatorname{D}_{v}F(0, \alpha^*))
=
\operatorname{span}\{\hat v_0\},
\qquad
\hat v_0
=
\frac1{\sqrt2}
\binom{1}{1}
e_{k^*},
\qquad
\|\hat v_0\|=1.
\]
Since the operator is self-adjoint in the reciprocal case, the adjoint eigenvector coincides with \(\hat v_0\), and therefore
\[
N(L^*)
=
\operatorname{span}\{\hat v_0\}.
\]

We introduce the projection onto the kernel
\begin{equation}\label{eq:app-PQ}
\Pi \mathbf f
:=
\langle \mathbf f,\hat v_0\rangle\,\hat v_0.
\end{equation}
Since
\[
R(L)=N(L^*)^\perp,
\]
the complementary projection \(I-\Pi\) projects onto the range \(R(L)\).

For comparison with the non-reciprocal case, we define in the reciprocal setting the analogue of the real part of the normal-form coefficient in~\eqref{eq: alpha_second_derivative} by
\begin{equation}\label{eq: real_part_lambda_reciprocal}
    \Re(\Lambda)
    :=
    -|k^*|^2 \bar u \hat W(k^*)\,\alpha''(0),
\end{equation}
evaluated at the critical mode \(k^*\).

\begin{theorem}\label{theorem:app}
Assume \(\widehat W(k^*)<0\), so that the crossing eigenvalue moves to the right as \(\alpha\)
increases, and
\(\iota:=\widehat W(2k^*)/\widehat W(k^*)<1\), so that
\(\sigma+\bar u(\alpha^*+\gamma)\widehat W(2k^*)=\sigma(1-\iota)>0\). Then the bifurcation branch
satisfies \(\alpha'(0)=0\) and
\begin{equation}\label{eq:app-result}
\alpha''(0) = -\frac{(\alpha^* +\gamma )^2\big(\widehat W(k^*)-2\widehat W(2k^*)\big)}
{8\pi\,\bar u \bigl(\sigma + \bar u (\alpha^*+\gamma) \hat W(2k^*)\bigr)} .
\end{equation}
In particular \(\operatorname{sign}\alpha''(0)=\operatorname{sign}(1-2\iota)\), so the branch is
supercritical if and only if \(\iota<\tfrac12\) and subcritical if and only if
\(\tfrac12<\iota<1\).
\end{theorem}

\begin{proof}
    We compute formula \cite[(I.6.11)]{kielh2004} for \(\alpha''(0)\) for a pitchfork bifurcation. 
    
    \emph{Step 1: Computing denominator.}
    Differentiating the abstract problem with respect to the variable \(u\) and then with respect to the bifurcation parameter \(\alpha\), we find 
    \begin{equation*}
        \langle \operatorname{D^2}_{x \alpha } F(0, \alpha^*) \hat v_0 , \hat v_0 \rangle = -\bar u\,|k^{*}|^2\widehat W(k^*)
    \end{equation*}

    \emph{Step 2: Computing numerator}
    We need to compute 
\begin{equation*}
    \Pi \operatorname{D^2}_{vv}F(0, \alpha^*) [\hat v_0, (I-\Pi) (\operatorname{D}_x F(0,\alpha^*))^{-1} (I-\Pi) \operatorname{D^2}_{vv}F(0,\alpha^*)[\hat v_0, \hat v_0] ].
\end{equation*}
We have 
\begin{equation*}
     D^2_{\mathbf{vv}}F(0, \alpha^*)[\hat v_0,\hat v_0] = -\frac{(\alpha^*+\gamma)\widehat W(k^*)k^{*2}}{\sqrt\pi}\binom{1}{1}e_{2k^*}
\end{equation*}
and \((I-\Pi) \binom{1}{1}e_{2k^*} = \binom{1}{1}e_{2k^*}\). Additionally, one can check that 
\begin{equation*}
    \operatorname{D}_x F(0,\alpha^*) \boldsymbol{e}_{2k^*} = -4|k^*|^2 (\sigma +\bar u (\alpha^*+\gamma) \hat W(2k^*)) \boldsymbol{e}_{2k^*} ,
\end{equation*}
which implies 
\begin{equation*}
  (\operatorname{D}_x F(0,\alpha^*))^{-1} \bigg(
    -\frac{(\alpha^*+\gamma)\widehat W(k^*)|k^*|^2}{\sqrt\pi}\binom{1}{1}e_{2k^*}
  \bigg)
  = \frac{(\alpha^*+\gamma)\widehat W(k^*)}{4 \sqrt\pi\,\bigl(\sigma + \bar u (\alpha^*+\gamma) \hat W(2k^*)\bigr)}\binom{1}{1} e_{2k^*}.
\end{equation*}
Now, the operator \(I-\Pi  \) acts like the identity on the above expression. 
Hence, we finally obtain at the expression
\begin{equation*}
    \frac{(\alpha^*+\gamma)\widehat W(k^*)}{4 \sqrt\pi\,\bigl(\sigma + \bar u (\alpha^*+\gamma) \hat W(2k^*)\bigr)} \Pi \operatorname{D^2}_{vv}F(0, \alpha^*) [\hat v_0, \boldsymbol{e}_{2k^*}]. 
\end{equation*}
The last computation follows similar to the computations in Proposition~\ref{prop: sub_supercritical} with the identities \(\nabla \cdot (\nabla e_{k^*} e_{2k^*} ) =  -\frac{3}{2\sqrt{\pi}} |k^*|^2 e_{3k^*} + \frac{1}{2\sqrt{\pi}}|k^*|^2 e_{k^*} \) and \(\nabla \cdot (\nabla e_{2k^*} e_{k^*} ) = -\frac{3}{\sqrt{\pi}} |k^*|^2 e_{3k^*} - \frac{1}{\sqrt{\pi}} |k^*|^2 e_{k^*}  \). The projection \(\Pi\) kills the \(3k^*\) mode and we are left with 
\begin{align*}
     &\frac{(\alpha^*+\gamma)\widehat W(k^*)}{4 \sqrt\pi\,\bigl(\sigma + \bar u (\alpha^*+\gamma) \hat W(2k^*)\bigr)} \langle \Pi \operatorname{D^2}_{vv}F(0, \alpha^*) [\hat v_0, \boldsymbol{e}_{2k^*}], \hat v_0 \rangle \\
     &\quad =\frac{(\alpha^*+\gamma)^2 |k^*|^{2}\widehat W(k^*)}{4\pi\,\bigl(\sigma + \bar u (\alpha^*+\gamma) \hat W(2k^*)\bigr)}\Big(\tfrac12\widehat W(k^*)-\widehat W(2k^*)\Big).
\end{align*}

\emph{Step 3: Putting everything together.}
We put the two formulas for the numerator and denominator together to find 
\begin{equation*}
    \alpha''(0) = -\frac{(\alpha^* +\gamma )^2\big(\widehat W(k^*)-2\widehat W(2k^*)\big)}{8\pi\,\bar u \bigl(\sigma + \bar u (\alpha^*+\gamma) \hat W(2k^*)\bigr)}. 
\end{equation*}
\end{proof}

\endappendix

\bibliographystyle{alpha}
\bibliography{quellen}

@article {amann1990,
    AUTHOR = {Amann, Herbert},
     TITLE = {Dynamic theory of quasilinear parabolic equations. {II}.
              {R}eaction-diffusion systems},
   JOURNAL = {Differential Integral Equations},
  FJOURNAL = {Differential and Integral Equations. An International Journal
              for Theory and Applications},
    VOLUME = {3},
      YEAR = {1990},
    NUMBER = {1},
     PAGES = {13--75},
      ISSN = {0893-4983},
   MRCLASS = {35K55 (58E07 92A09)},
  MRNUMBER = {1014726},
MRREVIEWER = {Lumin\ Jiang},
}

@book {engel2000,
    AUTHOR = {Engel, Klaus-Jochen and Nagel, Rainer},
     TITLE = {One-parameter semigroups for linear evolution equations},
    SERIES = {Graduate Texts in Mathematics},
    VOLUME = {194},
      NOTE = {With contributions by S. Brendle, M. Campiti, T. Hahn, G.
              Metafune, G. Nickel, D. Pallara, C. Perazzoli, A. Rhandi, S.
              Romanelli and R. Schnaubelt},
 PUBLISHER = {Springer-Verlag, New York},
      YEAR = {2000},
     PAGES = {xxii+586},
      ISBN = {0-387-98463-1},
   MRCLASS = {47D06 (34G10 35K90 47N20)},
  MRNUMBER = {1721989},
MRREVIEWER = {Charles\ Batty},
}

@book {kielh2004,
    AUTHOR = {Kielh\"ofer, Hansj\"org},
     TITLE = {Bifurcation theory},
    SERIES = {Applied Mathematical Sciences},
    VOLUME = {156},
      NOTE = {An introduction with applications to PDEs},
 PUBLISHER = {Springer-Verlag, New York},
      YEAR = {2004},
     PAGES = {viii+346},
      ISBN = {0-387-40401-5},
   MRCLASS = {47J15 (35B32 37G15 37K50 47N20 58E07)},
  MRNUMBER = {2004250},
MRREVIEWER = {Thomas\ Bartsch},
       DOI = {10.1007/b97365},
       URL = {https://doi.org/10.1007/b97365},
}

@article {Carrillo2026,
    AUTHOR = {Carrillo, Jos\'e{} A. and Salmaniw, Yurij},
     TITLE = {Long-time behaviour and bifurcation analysis of a two-species
              aggregation-diffusion system on the torus},
   JOURNAL = {Calc. Var. Partial Differential Equations},
  FJOURNAL = {Calculus of Variations and Partial Differential Equations},
    VOLUME = {65},
      YEAR = {2026},
    NUMBER = {1},
     PAGES = {Paper No. 19, 61},
      ISSN = {0944-2669,1432-0835},
   MRCLASS = {35B32 (35P05 35Q92 35R05 35R09)},
  MRNUMBER = {4998279},
       DOI = {10.1007/s00526-025-03132-0},
       URL = {https://doi.org/10.1007/s00526-025-03132-0},
}

@article {Carrillo2020,
    AUTHOR = {Carrillo, J. A. and Gvalani, R. S. and Pavliotis, G. A. and
              Schlichting, A.},
     TITLE = {Long-time behaviour and phase transitions for the
              {M}c{K}ean-{V}lasov equation on the torus},
   JOURNAL = {Arch. Ration. Mech. Anal.},
  FJOURNAL = {Archive for Rational Mechanics and Analysis},
    VOLUME = {235},
      YEAR = {2020},
    NUMBER = {1},
     PAGES = {635--690},
      ISSN = {0003-9527,1432-0673},
   MRCLASS = {35R60 (35B40 35K59 60K35)},
  MRNUMBER = {4062483},
MRREVIEWER = {Alp\ O.\ Eden},
       DOI = {10.1007/s00205-019-01430-4},
       URL = {https://doi.org/10.1007/s00205-019-01430-4},
}

@article {Jakub2026,
    AUTHOR = {Carrillo, Jos\'e{} A. and Salmaniw, Yurij and Skrzeczkowski,
              Jakub},
     TITLE = {Well-posedness of aggregation-diffusion systems with irregular
              kernels},
   JOURNAL = {Ann. Inst. H. Poincar\'e{} C Anal. Non Lin\'eaire},
  FJOURNAL = {Annales de l'Institut Henri Poincar\'e{} C. Analyse Non
              Lin\'eaire},
    VOLUME = {43},
      YEAR = {2026},
    NUMBER = {5},
     PAGES = {1201--1258},
      ISSN = {0294-1449,1873-1430},
   MRCLASS = {35K40 (35A01 35A02 35B65 35K55 35Q92)},
  MRNUMBER = {5097436},
       DOI = {10.4171/aihpc/173},
       URL = {https://doi.org/10.4171/aihpc/173},
}

@article {Portisch2022,
    AUTHOR = {J\"ungel, Ansgar and Portisch, Stefan and Zurek, Antoine},
     TITLE = {Nonlocal cross-diffusion systems for multi-species populations
              and networks},
   JOURNAL = {Nonlinear Anal.},
  FJOURNAL = {Nonlinear Analysis. Theory, Methods \& Applications. An
              International Multidisciplinary Journal},
    VOLUME = {219},
      YEAR = {2022},
     PAGES = {Paper No. 112800, 26},
      ISSN = {0362-546X,1873-5215},
   MRCLASS = {35K40 (35K55 35Q92 68T07 92B20 92D25)},
  MRNUMBER = {4379345},
       DOI = {10.1016/j.na.2022.112800},
       URL = {https://doi.org/10.1016/j.na.2022.112800},
}

@book{Boltzman,
author = {Ludwig Boltzmann},
title = {Further Studies on the Thermal Equilibrium of Gas Molecules},
booktitle = {The Kinetic Theory of Gases},
chapter = {},
YEAR = {2003},
pages = {262-349},
doi = {10.1142/9781848161337_0015},
URL = {https://www.worldscientific.com/doi/abs/10.1142/9781848161337_0015},
eprint = {https://www.worldscientific.com/doi/pdf/10.1142/9781848161337_0015}
}

@article{pottslewis2019,
    AUTHOR = {Potts, Jonathan R. and Lewis, Mark A.},
     TITLE = {Spatial memory and taxis-driven pattern formation in model
              ecosystems},
   JOURNAL = {Bull. Math. Biol.},
  FJOURNAL = {Bulletin of Mathematical Biology. A Journal Devoted to
              Research at the Interface of the Life and Mathematical
              Sciences},
    VOLUME = {81},
      YEAR = {2019},
    NUMBER = {7},
     PAGES = {2725--2747},
      ISSN = {0092-8240,1522-9602},
   MRCLASS = {92D25 (92C17)},
  MRNUMBER = {3977631},
       DOI = {10.1007/s11538-019-00626-9},
       URL = {https://doi.org/10.1007/s11538-019-00626-9},
}

@article{giunta2022detecting,
    AUTHOR = {Giunta, Valeria and Hillen, Thomas and Lewis, Mark A. and
              Potts, Jonathan R.},
     TITLE = {Detecting minimum energy states and multi-stability in
              nonlocal advection-diffusion models for interacting species},
   JOURNAL = {J. Math. Biol.},
  FJOURNAL = {Journal of Mathematical Biology},
    VOLUME = {85},
      YEAR = {2022},
    NUMBER = {5},
     PAGES = {Paper No. 56, 44},
      ISSN = {0303-6812,1432-1416},
   MRCLASS = {92D25 (35B38 35Q92)},
  MRNUMBER = {4499096},
MRREVIEWER = {Shanming\ Ji},
       DOI = {10.1007/s00285-022-01824-1},
       URL = {https://doi.org/10.1007/s00285-022-01824-1},
}

@article{giunta2024weakly,
    AUTHOR = {Giunta, Valeria and Hillen, Thomas and Lewis, Mark A. and
              Potts, Jonathan R.},
     TITLE = {Weakly nonlinear analysis of a two-species non-local
              advection-diffusion system},
   JOURNAL = {Nonlinear Anal. Real World Appl.},
  FJOURNAL = {Nonlinear Analysis. Real World Applications. An International
              Multidisciplinary Journal},
    VOLUME = {78},
      YEAR = {2024},
     PAGES = {Paper No. 104086, 19},
      ISSN = {1468-1218,1878-5719},
   MRCLASS = {92D25 (35B32 35B36 35Q56)},
  MRNUMBER = {4708765},
MRREVIEWER = {Yan\ Li},
       DOI = {10.1016/j.nonrwa.2024.104086},
       URL = {https://doi.org/10.1016/j.nonrwa.2024.104086},
}

@article{JKME23,
    AUTHOR = {Jewell, Thomas Jun and Krause, Andrew L. and Maini, Philip K.
              and Gaffney, Eamonn A.},
     TITLE = {Patterning of nonlocal transport models in biology: the impact
              of spatial dimension},
   JOURNAL = {Math. Biosci.},
  FJOURNAL = {Mathematical Biosciences},
    VOLUME = {366},
      YEAR = {2023},
     PAGES = {Paper No. 109093, 16},
      ISSN = {0025-5564,1879-3134},
   MRCLASS = {99-01},
  MRNUMBER = {4663459},
       DOI = {10.1016/j.mbs.2023.109093},
       URL = {https://doi.org/10.1016/j.mbs.2023.109093},
}

@article{FruchartVitelli2021,
  author  = {Fruchart, Michel and Hanai, Ryo and Littlewood, Peter B. and Vitelli, Vincenzo},
  title   = {Non-reciprocal phase transitions},
  journal = {Nature},
  volume  = {592},
  number  = {7854},
  pages   = {363--369},
  year    = {2021},
  doi     = {10.1038/s41586-021-03375-9},
}

@article{YouBaskaranMarchetti2020,
    AUTHOR = {You, Zhihong and Baskaran, Aparna and Marchetti, M. Cristina},
     TITLE = {Nonreciprocity as a generic route to traveling states},
   JOURNAL = {Proc. Natl. Acad. Sci. USA},
  FJOURNAL = {Proceedings of the National Academy of Sciences of the United
              States of America},
    VOLUME = {117},
      YEAR = {2020},
    NUMBER = {33},
     PAGES = {19767--19772},
      ISSN = {0027-8424,1091-6490},
   MRCLASS = {92D10},
  MRNUMBER = {4245726},
       DOI = {10.1073/pnas.2010318117},
       URL = {https://doi.org/10.1073/pnas.2010318117},
}

@article{SahaGolestanian2020,
  title = {Scalar Active Mixtures: The Nonreciprocal Cahn-Hilliard Model},
  author = {Saha, Suropriya and Agudo-Canalejo, Jaime and Golestanian, Ramin},
  journal = {Phys. Rev. X},
  volume = {10},
  issue = {4},
  pages = {041009},
  numpages = {14},
  year = {2020},
  month = {Oct},
  publisher = {American Physical Society},
  doi = {10.1103/PhysRevX.10.041009},
  url = {https://link.aps.org/doi/10.1103/PhysRevX.10.041009}
}

@article{FrohoffThiele2023,
  title = {Nonreciprocal Cahn-Hilliard Model Emerges as a Universal Amplitude Equation},
  author = {Frohoff-H\"ulsmann, Tobias and Thiele, Uwe},
  journal = {Phys. Rev. Lett.},
  volume = {131},
  issue = {10},
  pages = {107201},
  numpages = {7},
  year = {2023},
  month = {Sep},
  publisher = {American Physical Society},
  doi = {10.1103/PhysRevLett.131.107201},
  url = {https://link.aps.org/doi/10.1103/PhysRevLett.131.107201}
}

@article{BraunsMarchetti2024,
  title = {Nonreciprocal Pattern Formation of Conserved Fields},
  author = {Brauns, Fridtjof and Marchetti, M. Cristina},
  journal = {Phys. Rev. X},
  volume = {14},
  issue = {2},
  pages = {021014},
  numpages = {22},
  year = {2024},
  month = {Apr},
  publisher = {American Physical Society},
  doi = {10.1103/PhysRevX.14.021014},
  url = {https://link.aps.org/doi/10.1103/PhysRevX.14.021014}
}

@article{Gerber2026,
doi = {10.1088/1361-6544/ae7ee7},
url = {https://doi.org/10.1088/1361-6544/ae7ee7},
year = {2026},
month = {jul},
publisher = {IOP Publishing},
volume = {39},
number = {7},
pages = {075023},
author = {Gerber, N J and Gvalani, R S and Hairer, M and Pavliotis, G A and Schlichting, A},
title = {Formation of clusters and coarsening in weakly interacting diffusions},
journal = {Nonlinearity}
}

@article{Kreienkamp2022,
doi = {10.1088/1367-2630/ac9cc3},
url = {https://doi.org/10.1088/1367-2630/ac9cc3},
year = {2022},
month = {dec},
publisher = {IOP Publishing},
volume = {24},
number = {12},
pages = {123009},
author = {Kreienkamp, Kim L and Klapp, Sabine H L},
title = {Clustering and flocking of repulsive chiral active particles with non-reciprocal couplings},
journal = {New Journal of Physics}
}

@article{Ivlev2015,
  title = {Statistical Mechanics where Newton's Third Law is Broken},
  author = {Ivlev, A. V. and Bartnick, J. and Heinen, M. and Du, C.-R. and Nosenko, V. and L{\"o}wen, H.},
  journal = {Phys. Rev. X},
  volume = {5},
  issue = {1},
  pages = {011035},
  numpages = {10},
  year = {2015},
  month = {Mar},
  publisher = {American Physical Society},
  doi = {10.1103/PhysRevX.5.011035},
  url = {https://link.aps.org/doi/10.1103/PhysRevX.5.011035}
}

@article{MartinEtAl2025,
  title = {Transition to Collective Motion in Nonreciprocal Active Matter: Coarse Graining Agent-Based Models into Fluctuating Hydrodynamics},
  author = {Martin, David and Seara, Daniel and Avni, Yael and Fruchart, Michel and Vitelli, Vincenzo},
  journal = {Phys. Rev. X},
  volume = {15},
  issue = {4},
  pages = {041015},
  numpages = {43},
  year = {2025},
  month = {Oct},
  publisher = {American Physical Society},
  doi = {10.1103/PhysRevX.15.041015},
  url = {https://link.aps.org/doi/10.1103/PhysRevX.15.041015}
}

@article{SotoGolestanian2014,
  title = {Self-Assembly of Catalytically Active Colloidal Molecules: Tailoring Activity Through Surface Chemistry},
  author = {Soto, Rodrigo and Golestanian, Ramin},
  journal = {Phys. Rev. Lett.},
  volume = {112},
  issue = {6},
  pages = {068301},
  numpages = {5},
  year = {2014},
  month = {Feb},
  publisher = {American Physical Society},
  doi = {10.1103/PhysRevLett.112.068301},
  url = {https://link.aps.org/doi/10.1103/PhysRevLett.112.068301}
}

@article{Meredith2020,
  title={Predator--prey interactions between droplets driven by non-reciprocal oil exchange},
  author={Meredith, Caleb H and Moerman, Pepijn G and Groenewold, Jan and Chiu, Yu-Jen and Kegel, Willem K and van Blaaderen, Alfons and Zarzar, Lauren D},
  journal={Nature Chemistry},
  volume={12},
  number={12},
  pages={1136--1142},
  year={2020},
  publisher={Nature Publishing Group UK London}
}

@article{Lavergne2019,
author = {François A. Lavergne  and Hugo Wendehenne  and Tobias Bäuerle  and Clemens Bechinger },
title = {Group formation and cohesion of active particles with visual perception–dependent motility},
journal = {Science},
volume = {364},
number = {6435},
pages = {70-74},
year = {2019},
doi = {10.1126/science.aau5347},
URL = {https://www.science.org/doi/abs/10.1126/science.aau5347},
eprint = {https://www.science.org/doi/pdf/10.1126/science.aau5347}}

@article{Theveneau2013,
  title={Chase-and-run between adjacent cell populations promotes directional collective migration},
  author={Theveneau, Eric and Steventon, Benjamin and Scarpa, Elena and Garcia, Simon and Trepat, Xavier and Streit, Andrea and Mayor, Roberto},
  journal={Nature cell biology},
  volume={15},
  number={7},
  pages={763--772},
  year={2013},
  publisher={Nature Publishing Group UK London}
}

@article{Tsyganov2003,
  title = {Quasisoliton Interaction of Pursuit-Evasion Waves in a Predator-Prey System},
  author = {Tsyganov, M. A. and Brindley, J. and Holden, A. V. and Biktashev, V. N.},
  journal = {Phys. Rev. Lett.},
  volume = {91},
  issue = {21},
  pages = {218102},
  numpages = {4},
  year = {2003},
  month = {Nov},
  publisher = {American Physical Society},
  doi = {10.1103/PhysRevLett.91.218102},
  url = {https://link.aps.org/doi/10.1103/PhysRevLett.91.218102}
}

@misc{ArjunEtAl2026,
  author       = {Arjun, R. and Patra, Pratyush Prakash and Anil Kumar, A. V.},
  title        = {Nonreciprocal {McKean}--{Vlasov} equations: from stationary instabilities to travelling waves},
  year         = {2026},
  eprint       = {2605.07597},
  archiveprefix = {arXiv},
  primaryclass = {cond-mat.stat-mech},
}

@misc{HuLiu2026,
  author       = {Hu, Junlang and Liu, Zhenxin},
  title        = {Continuity and discontinuity of {McKean}--{Vlasov} phase transitions via bifurcation theory},
  year         = {2026},
  eprint       = {2607.10723},
  archiveprefix = {arXiv},
  primaryclass = {math.AP},
}

\medskip

\end{document}